\documentclass[12pt, reqno]{amsart}
\usepackage{amsmath, amsthm, amscd, amsfonts, amssymb, graphicx, color, float}
\usepackage[bookmarksnumbered, colorlinks, plainpages]{hyperref}
\input{mathrsfs.sty}
\hypersetup{colorlinks=true,linkcolor=red, anchorcolor=green, citecolor=cyan, urlcolor=red, filecolor=magenta, pdftoolbar=true}
\usepackage{soul,cancel}
\usepackage{tikzsymbols}
\usepackage{textcomp}
\usepackage{enumerate}

\newtheorem{theorem}{Theorem}[section]

\newtheorem{proposition}[theorem]{Proposition}
\newtheorem{corollary}[theorem]{Corollary}
\theoremstyle{definition}
\newtheorem{definition}[theorem]{Definition}
\newtheorem{example}[theorem]{Example}

\theoremstyle{remark}

\numberwithin{equation}{section}

\def\11{\textbf{1}}

\begin{document}

\title[Similarities and differences between real and complex spaces]{Similarities and differences between real and complex Banach spaces. An overview and recent developments}

\author[Moslehian]{M. S.~Moslehian}
\address[M.S.~Moslehian]{\mbox{}\newline \indent Department of Pure Mathematics, \newline \indent Center of Excellence in
 Analysis on Algebraic Structures (CEAAS), \newline \indent Ferdowsi University of Mashhad, \newline \indent P.O. Box 1159, Mashhad 91775, Iran}
\email{moslehian@um.ac.ir; moslehian@yahoo.com}
\author[Mu\~{n}oz]{G. A. Mu\~{n}oz-Fern\'{a}ndez}
\address[Gustavo A. Mu\~{n}oz-Fern\'{a}ndez]{\mbox{}\newline \indent Instituto de Matem\'atica Interdisciplinar (IMI), \newline \indent Departamento de An\'{a}lisis Matem\'{a}tico y Matem\'atica Aplicada,\newline \indent Facultad de Ciencias Matem\'aticas, \newline \indent Plaza de Ciencias 3, \newline \indent Universidad Complutense de Madrid,\newline \indent 28040 Madrid, Spain.}
\email{gustavo\_fernandez@mat.ucm.es}
\author[Peralta]{A. M.~Peralta}
\address[A.M.~Peralta]{\mbox{}\newline \indent Departamento de An\'{a}lisis Matem\'{a}tico \newline \indent Facultad de Ciencias\newline \indent Universidad de Granada,
Granada, 18071, Spain.}
\email{aperalta@ugr.es}
\author[Seoane]{J. B.~Seoane-Sep\'{u}lveda}
\address[J.B.~Seoane-Sep\'{u}lveda]{\mbox{}\newline \indent Instituto de Matem\'atica Interdisciplinar (IMI)\newline \indent Departamento de An\'{a}lisis Matem\'{a}tico y Matem\'atica Aplicada \newline \indent
 Facultad de Ciencias Matem\'{a}ticas\newline \indent
 Plaza de Ciencias 3 \newline \indent
 Universidad Complutense de Madrid \newline \indent
 Madrid, 28040, Spain.}
\email{jseoane@mat.ucm.es}

\renewcommand{\subjclassname}{\textup{2020} Mathematics Subject Classification}

\begin{abstract} There are numerous cases of discrepancies between results obtained in the setting of real Banach spaces and those obtained in the complex context. This article is a modern exposition of the subtle differences between key results and theories for complex and real Banach spaces and the corresponding linear operators between them. We deeply discuss some aspects of the complexification of real Banach spaces and give several examples showing how drastically different can be the behavior of real Banach spaces versus their complex counterparts.
\end{abstract}

\maketitle

\tableofcontents

\section{Introduction}

By a complex (respectively, real) linear space, we mean a linear space over the field of complex numbers $\mathbb{C}$ (respectively, real numbers $\mathbb{R}$).\smallskip

In the theory of Banach spaces and operator algebras these objects are usually considered over the field of complex numbers, and a study over the field of real numbers has been systematically studied in recent years (see \cite{LI}). Although $\mathbb{R}$ has very good properties such as the Dedekind completion (i.e., every upper-bounded nonempty subset of $\mathbb{R}$ admits a supremum) and the law of trichotomy (i.e., every nonzero real number is either positive or negative), it fails to satisfy the fundamental theorem of algebra (i.e., there exist nonconstant single-variable polynomials over $\mathbb{R}$ admitting no root in $\mathbb{R}$). These facts entail that functional analysts and operator theorists usually deal with complex linear spaces. Furthermore, there are several pieces of evidence showing that complex linear spaces are more suitable to be used in physics. For example, in quantum mechanics, the state of a system is described as a vector in a complex Hilbert space. On more than one occasion, one has been faced with results that are valid in the setting of complex Banach spaces and algebras but not in the real context, and vice versa, for example, the existence of elements in a unital real Banach algebra whose spectra contain no real numbers, or the impossibility of establishing a version of the Mazur--Ulam theorem for complex normed spaces, which provides a complex affine extension of every surjective isometry in that setting. Furthermore, we know from results by Bourgain and Szarek the existence of two complex Banach spaces which are linearly isometric as real spaces but non-isomorphic as complex spaces \cite{BOU,Sza86,Sza86b}. An outstanding result by Ferenczi shows that there exist two Banach spaces which are isometric as real spaces, but totally in-comparable as complex spaces; where two real or complex Banach spaces are said to be \emph{totally incomparable} if no infinite-dimensional subspace of the one is isomorphic to a subspace of the other (see \cite[Theorem 1]{FER}). Furthermore, this result shows that a theorem of Godefroy and Kalton proving that if a separable real Banach space embeds isometrically into a Banach space, then it embeds linearly isometrically into it (see  \cite{GK}),  cannot be extended to the complex case.\smallskip

In many topics such as real C$^*$-algebras, JB-algebras, real operator spaces, and $KK$-theory, mathematicians study mathematical objects in the setting of real linear spaces; special attention is received by bounded real linear operators acting on real Hilbert spaces. Thus, it is an interesting problem to ask which results of the theory in the complex case still hold for the real case, probably under some different hypotheses or more restricting assumptions, and which facts valid in the complex case do not hold when we restrict ourselves to linear spaces over the real field. Generally, there are some technical difficulties in translating the known results related to the complex case to the real setting. One of the essential tools is the ``complexification'' by means of which one may go from the real to the complex world and prove the new results or employ some known facts and results, returning then to the real setting in order to state the pursued results therein. This idea has been used many times with the aim of extending the inherent results of complex analysis to a real setting. \smallskip

A celebrated example of the use of complexification techniques to tackle a problem in a real setting is the study of real analytic functions on (real) Banach spaces with the aid of the properties of holomorphic functions. Let ${\mathscr X}$ and ${\mathscr Y}$ be Banach spaces over ${\mathbb K}={\mathbb R}$ or ${\mathbb C}$ and let $U\subseteq {\mathscr X}$ be an open set. A mapping $P:\mathscr{X}\rightarrow \mathscr{Y}$ is an \emph{$n$-homogeneous polynomial} if there is an $n$-linear mapping $L:\mathscr{X}^n \rightarrow \mathscr{Y}$ such that $P(x)=L(x,\ldots,x)$ for all $x\in{\mathscr X}$. The fact that a homogeneous polynomial $P:{\mathscr X}\rightarrow{\mathscr Y}$ or an $n$-linear mapping $L:{\mathscr X}^n\rightarrow {\mathscr Y}$ is continuous if and only if it is bounded on the unit ball, $B_{\mathscr X}$, of ${\mathscr X}$, is a standard result in the theory of polynomials between arbitrary Banach spaces (see, for instance, \cite{Dineen1999} for a modern exposition on polynomials). If $P:\mathscr{X}\rightarrow \mathscr{Y}$ and $L:\mathscr{X}^n\rightarrow \mathscr{Y}$ are a continuous $n$-homogeneous polynomial and a continuous $n$-linear mapping, respectively, then we define the norms of $P$ and $L$ by
\begin{align*}
\|P\|&=\sup\{\|P(x)\|:x\in B_{\mathscr X}\}\text{ and }\\
\|L\|&=\sup\{\|L(x_1,\dots,x_n)\|:x_1,\ldots , x_n\in B_{\mathscr X}\}, \hbox{ respectively}.
\end{align*}

A function $f:U\subseteq {\mathscr X}\rightarrow {\mathscr Y}$ is analytic (also called holomorphic in the case when ${\mathbb K}={\mathbb C}$) if $f$ is defined by its Taylor series around every point $a$ in $U$, that is,
 \begin{equation}\label{eq:Taylor_Series}
 f(x)=\sum_{k=0}^\infty P_n(x-a),
 \end{equation}
for all $x$ in the open ball centered at $a$ with radius $\rho$ ($B(a,\rho)$ in short), where $P_n=\frac{1}{n!}{\widehat D}^nf$ (note that here ${\widehat D}^nf$ stands for the $n$-homogeneous polynomial associated with the $n$-th Fr\'echet derivative of $f$) and $\rho>0$ is the radius of convergence of the series in \eqref{eq:Taylor_Series}. Recall that $\rho$ can be calculated by using the Cauchy--Hadamard formula
 $$
 \rho=\frac{1}{\limsup_n \|P_n\|^{\frac1n}}.
 $$
For more details on real and complex analytic mappings between Banach spaces, the reader may consult \cite{Bochnak1970,Bochnak1971bis,Bochnak1971,Dineen1999}.\smallskip

It is well known that in the case where ${\mathscr X}={\mathscr Y}={\mathbb K}$, we have $\rho=\text{dist}(a,\partial U)$ and that \eqref{eq:Taylor_Series} is uniformly convergent to $f$ on $B(a,r)$ for all $r\in(0,\rho)$. This is not necessarily true for arbitrary ${\mathscr X}$ and ${\mathscr Y}$, which motivates the definition of fully analytic functions. The mapping $f:{\mathscr X}\rightarrow {\mathscr Y}$ is \emph{fully analytic} in $U$ if its Taylor series \eqref{eq:Taylor_Series} converges uniformly in every closed ball centered at $a$ contained in $U$ for each $a\in U$. Thus we define the \emph{radius of analyticity}, $\rho_A$, of $f$ at $a$ as the largest $r>0$ such that $f$ is fully analytic in the ball $B(a, r)$. Obviously, $\rho_A\leq \rho$ for every fully analytic mapping $f:{\mathscr X}\rightarrow {\mathscr Y}$. Moreover, if ${\mathbb K}={\mathbb C}$, then it follows from the Cauchy integral formula that $\rho_A=\rho$. It is not known whether $\rho_A=\rho$ is also true for real analytic functions. The interest in this problem can be traced back at least to 1938, when Taylor \cite{Taylor1938} proved that $\rho_A\ge \frac{\rho}{\sqrt{2} e}$. Using the optimal complexification constants of homogeneous polynomials (see \cite{MUN}), the previous estimate can be improved to $\rho_A\geq \frac{\rho}{2}$. A further improvement obtained by Nguyen \cite{Nguyen2009} in 2009 shows that $\rho_A\geq \frac{\rho}{\sqrt{e}}$. The best estimate known nowadays is $\rho_A\geq \frac{\rho}{\sqrt{2}}$ (see \cite{Hajek2014,Sarantopoulos2016}); however, that estimate can be greatly improved for specific spaces. For instance, $\rho_A=\rho$ for any (real or complex) Hilbert space (see the discussion after the proof of \cite[Theorem 1]{Nguyen2009}). Also, it was established in \cite{Boyd2018} that if $\rho_A=\rho$ in $\ell_1({\mathbb R})$, then we also have $\rho_A=\rho$ for any real Banach space.\smallskip

As shown in the previous paragraphs, the importance of complexifications is revealed in the study of real analytic functions, but it manifests too in the study of many other important questions being presented in the next sections. Complexifications have been employed several times in the past, and nowadays are still subjected to study. A unified treatment on complexifications was done in \cite{MUN} (see also \cite {Kirwan}), where, in addition to a number of general results on the construction of several complexification norms, some optimal estimates on the norm of the complex extension of polynomials and multilinear mappings are proved. Other studies on complexifications can be found in \cite{Cuellar, IliKuzLiPoon, RUA, Zyl}.\smallskip

Complexifications and real forms also constitute a key procedure to study real C$^*$-algebras \cite{IsRo, ChuDangRuVen, LI}, real J$^*$B-algebras \cite{Alv86, DaRu94, Pe03}, and real JB$^*$-triples \cite{IsKaRo95}, objects intensively studied in the nineties, and whose topicality is out of any doubt. Actually, in certain problems, real structures are gaining protagonism and topicality. For example, the conclusion of the Mazur--Ulam theorem produces real affine maps, while the recent contributions on Tingley's problem on the extension of isometries between the unit spheres of Banach spaces, C$^*$-algebras, and JB$^*$-triples show that the desired extension is only real linear, and the theory of real structures and morphisms is becoming more useful (see, for example, \cite{Pe2018, YangZhao2014, MoriOza2020, BeCuFerPe2021, KalPe2020, Banakh, BanakhCabello}). It is worth exploring the parallelisms and divergences of celebrated results, like the Gleason--Kahane--\.{Z}elazko theorem, the Kadison--Schwarz inequality, the notions of $n$-positive maps, the Russo--Dye theorem, the Bohnenblust--Karlin theorem, and the Kaup--Banach--Stone theorem, in the real and complex settings. These results are revisited together with a complete presentation of the original sources, the state-of-art of problems, and open questions.\smallskip

We have tried to write this expository article in a self-contained manner. However, a background of basic topics in the theory of Banach spaces is needed for an adequate understanding of the topic.

\section{Real vs complex linear spaces}

This section contains the algebraic tools required to define the complexification of a real linear space. The first subsections are devoted to refreshing the basic notions on complexifications of real linear spaces from a purely algebraic point of view. The analytic ingredients will appear in subsequent subsections.

\subsection{Linear algebra}\ \smallskip

By restricting the scalars to the real numbers, every complex linear space $\mathscr{X}$ can be regarded as a real linear space, denoted by $\mathscr{X}_r$. This process is called \emph{realification}. \smallskip

If $\mathscr{X}$ is a complex linear space of dimension $n$ with a linear basis $\{f_1, \ldots, f_n\}$, then $\{f_1, if_1, \ldots, f_n, if_n\}$ is a linear basis for $\mathscr{X}_r$ and so $\mathscr{X}_r$ is of dimension $2n$.\smallskip

A linear subspace $\mathscr{M}$ of $\mathscr{X}$ is a linear subspace of $\mathscr{X}_r$, but the converse is not true. For example, any nontrivial subspace of $\mathbb{R}^2$ is a line passing $0$, but clearly, it is not a subspace of $\mathbb{C}_r\cong \mathbb{R}^2$.\smallskip

Using one of the determinants' basic properties (i.e., $\det(TS)=\det(T)\det(S)$), we can conclude that for each odd natural $n$ there is no matrix $T$ in $\mathcal{M}_n(\mathbb{R})$ such that $T^{2}=-I$, whilst in $\mathcal{M}_n(\mathbb{C})$, there exist many examples of matrices satisfying this property.\smallskip

For every commuting $n\times n$ complex matrices $T$ and $S$, there exists a unitary matrix $U\in \mathcal{M}_n(\mathbb{C})$ such that both $U^*TU$ and $U^*SU$ are upper-triangular, where $^*$ denotes the conjugate transpose operation. This result does not hold for matrices in $\mathcal{M}_n(\mathbb{R})$. For example, if $T=\left(\begin{matrix} 0 & 1 \\ -1& 0\end{matrix}\right)$ and $S=\left(\begin{matrix} 1 & 1 \\ -1 & 1\end{matrix}\right)$, then there is no matrix $U\in \mathcal{M}_n(\mathbb{R})$ with the required properties (see \cite[p. 76, Problem 3.]{ZHA}).

\subsection{Complexification of real linear spaces}\label{subsec: algebraic complexification and operators}\ \smallskip

The process of producing a complex linear space from a real one is called \emph{complexification}. In such a process, the method follows similar techniques to those employed to construct $\mathbb{C}$ from $\mathbb{R}$.\smallskip

A complex linear space $\mathscr{X}$ is called a complexification of a real linear space ${X}$ if there is a one-to-one/ injective real linear map $\iota: X \to \mathscr{X}$ such that the complex linear span of $\iota({X})$ is $\mathscr{X}$. Obviously, such a complex linear space is unique and will be denoted by ${X}_c$.\smallskip

If ${X}$ is a real linear space, then the direct sum ${X}_c={X}\oplus_{\mathbb{R}} {X}$ as a real linear space can be endowed with a complex structure via the product by complex scalars defined by
\[(\alpha+i\beta)(x,y):=(\alpha x-\beta y, \alpha y+\beta x).\]
Identifying ${X}$ with $\{(x,0): x\in {X}\}$ and denoting $\{(0,x): x\in {X}\}$ by $i{X}$, we can write ${X}_c={X}\oplus_{\mathbb{R}} i {X}$ and denote $(x,y)\in X_c$ by $x+iy$. Note that ${X}_c={X}\oplus_{\mathbb{R}} i {X}$ can be identified with ${X}\otimes_{\mathbb{R}} \mathbb{C}$ in the context of real linear spaces and $({X}\otimes_{\mathbb{R}} \mathbb{C})_r$ is nothing but $X\otimes \ell_2^2$ via  $x \otimes (r+is)\mapsto x\otimes (r,s)$, where $\ell_2^2 \cong \mathbb{C}_{r}$ is the two-dimensional real Hilbert space. In what follows, all tensor products are real.\smallskip

Furthermore, the mappings $\sigma, \tau: {X}_c\to {X}_c$ defined by $\sigma(x,y)=(-y,x)$ and $\tau(x,y)= \overline{(x,y)}=(x,-y)$, respectively, are complex linear and conjugate linear automorphisms on ${X}_c$ satisfying $\sigma^2=-I_{{X}_c}$ and $\tau^2=I_{{X}_c}$. As an example, $\mathcal{M}_n(\mathbb{C})=\mathcal{M}_n(\mathbb{R}) + i \mathcal{M}_n(\mathbb{R})$.\smallskip

Let ${X}$ and ${Y}$ be two real linear spaces and let ${X}_c$ and ${Y}_c$ be their complexifications, respectively. If $T: {X} \to {Y}$ is a real linear mapping, then one can define its complex linear extension $T_c: {X}_c\to {Y}_c$ by $T_c(x+iy)=T(x) + i T(y)$. For each linear mapping $S: {X}_c \to {Y}_c$, we can consider the complex linear operator $\overline{S}: {X}_c\to {Y}_c$ by $\overline{S}(x+iy)=\overline{S(x-iy)}=\overline{S\left(\overline{x+iy}\right)}$. It is easy to check that 
$S$ is of the form $T_c$ for some real linear map $T: {X} \to {Y}$ if and only if $S=\overline{S}$. Let $L(X_c,Y_c)_{sym}$ denote the real linear space of all $S\in L(X_c,Y_c)$ such that $\overline{S} = S$. It follows from the previous arguments that the mapping \begin{equation}\label{eq identification of L(X,Y) is the corresponding complex linear extensions}\begin{aligned} L(X,Y)&\longrightarrow L(X_c,Y_c)_{sym} \\
T& \longmapsto T_c
\end{aligned}
\end{equation} is a real linear bijection, so $L(X,Y)$ and $L(X_c,Y_c)_{sym}$ are indistinguishable as real linear spaces.


\section{Complex structure}\label{subsect: complex structure}

Let ${X}$ be a real linear space. If there is an automorphism $\sigma: {X}\to {X}$ such that $\sigma^2=-I_{X}$, then the product by scalars defined by $(\alpha+i\beta) x:=\alpha x+\beta \sigma(x)$ makes ${X}$ a complex linear space. In this case, we say that ${X}$ admits a complex  structure and denote it by ${X}_\sigma$.\smallskip

Conversely, for each complex linear space $\mathscr{X}$, the mapping $\sigma: \mathscr{X}_{r}\to \mathscr{X}_{r}$ defined by $\sigma(x):= i x$ is an automorphism on $\mathscr{X}_{r}$ with $\sigma^2=-I_{X}$.\smallskip

It is easy to show that for a linear space $X$ of finite dimension $n$, a necessary and sufficient condition for the existence of an automorphism $\sigma$ on $X$ with $\sigma^2=-I_{X}$ is that $n$ is even. In fact, if $\sigma: \ell_2^n \to \ell_2^n$ is such a map, then, by considering it as an $n\times n$ matrix with real entries, we have
$$ \det(\sigma)^2=\det(\sigma^2)=\det(-I_{\ell_2^2})=(-1)^n.$$ \smallskip


\section{Complexifications of Banach spaces}

After reviewing the basic algebraic construction employed to construct the complexification of a real linear space, we try to extend the analytic structure defined by a norm on a real linear space to an appropriate norm on the complexification. In this section, we focus on the complexification of Banach spaces. We will see that, in this case, we have different approaches to extend the norm to the algebraic complexification.\smallskip

First, note that if $({X}, \|\cdot\|)$ is a real Banach space with a continuous automorphism $\sigma$ on $X$ satisfying $\sigma^2=-I_{X}$, then the complex linear space ${X}_{\sigma}$ together with the product by complex scalars defined in Section \ref{subsect: complex structure} and the norm given by
\begin{align}\label{norm by sigma}
\|x\|_0=\sup_{t\in[0,2\pi]}\|x \cos t + \sigma(x) \sin t\|
\end{align}
turns into a complex Banach space and $\|x\| \leq \|x\|_0\leq (1+\|\sigma\|) \|x\|$ for all $x\in X$. Dieudonn\'{e} \cite{DIE} was one of the first authors who found an example of an infinite-dimensional real Banach space ${X}$ (the James space) admitting no automorphism $\sigma$ with the required property.\smallskip

There are interesting questions on the existence, uniqueness, and the number of different complex structures, up to isomorphisms, coexisting in a concrete real Banach space. There exist many examples of Banach spaces admitting no complex structure (see, for example, \cite{DIE}), having a unique complex structure (see \cite{FG}), admitting exactly $n+1$ nonequivalent complex structures (see \cite{FER}), or having infinitely many complex structures up to isomorphisms (see \cite{CUE}, for more examples).


\subsection{Complex conjugation}\label{complex conjugation1}\ \smallskip

Let $\mathscr{X}$ be a complex Banach space. Also, let $\overline{\mathscr{X}}$ be the complex conjugate of $\mathscr{X}$, that is, the same set $\mathscr{X}$ equipped with the same norm and addition but whose product by complex scalars is replaced with $\lambda\cdot x:=\overline{\lambda}x$. Clearly, the spaces $\mathscr{X}$ and $\overline{\mathscr{X}}$ are isometrically isomorphic as real Banach spaces; however, Bourgain \cite{BOU} presented an example of a complex Banach space $\mathscr{X}$, given by an $\ell_2$-sum of finite-dimensional spaces whose distances to their respective conjugate spaces tend to infinity, such that $\mathscr{X}$ and $\overline{\mathscr{X}}$ are not isomorphic as complex Banach spaces. Thus complex Banach spaces may be isomorphic as real Banach spaces while they are not as complex spaces.\smallskip

It is worth noting that for each continuous automorphism $\sigma$ on $X$ satisfying $\sigma^2=-I_{X}$, we have ${X}_{-\sigma}=\overline{{X}_\sigma}$, that is, ${X}_{-\sigma}$ is the complex conjugate space of ${X}_\sigma$.


\subsection{Various norms on the complexification of a real Banach space}\label{subsec: various nroms on Xc}\ \smallskip

Assume that $({X},\|\cdot\|)$ is a real Banach space. There are many ways, in general, to define a (complete) norm $\|\cdot\|_c$ on the algebraic complexification ${X}_c$ whose restriction to $X$ gives the original norm $\|\cdot\|$. Any such a space $({X}_c,\|\cdot\|_c)$ is called a \emph{complexification} of the Banach space $({X},\|\cdot\|)$. A complexification whose norm is \emph{reasonable} (i.e., $\|x-iy\|=\|x+iy\|$) is called a \emph{reasonable complexification} of ${X}$.\smallskip

The minimal reasonable complexification is the \emph{Taylor complexification}\label{label Taylor complexification} defined by
$$\|x+iy\|_T:=\sup_{t\in[0,2\pi]}\| x \cos t -y \sin t\|.$$
In addition, any reasonable norm $\||\cdot |\|$ on ${X}_c$ is equivalent to the Taylor norm\label{eq any reasonable norm is equivalent to the T norm} since, from $\|| x+iy |\|=\||(\cos(t)+i \sin (t))(x+iy)|\|$, one can easily infer (see \cite{MUN}) that
$$\|x+iy\|_T \leq \||x+iy|\|\leq 2 \|x+iy\|_T.$$

Tensor norms provide a systematic way to define reasonable norms on the complexification. The interested reader is invited to consult \cite{DF1993,Ryan2002} for a complete account on tensors. If ${X}$ is a real Banach space, any tensor norm $\alpha$ on ${X}\otimes \ell_2^2$ is a reasonable norm on the complexification (see, for instance, \cite[Proposition 9]{MUN}). As a matter of fact, it can be proved (see, for example, \cite{MUN}) that
 $$
 \|x+iy\|_T=\sup\{\sqrt{\varphi(x)^2+\varphi(y)^2}:\varphi \in {X}^*\text{ and }\|\varphi\|_{{X}^*}\leq 1\}.
 $$
In other words, $\|\cdot\|_T$ can be alternatively described in terms of the injective tensor norm in ${X}\otimes \ell_2^2$ or, equivalently, $({X}_c,\|\cdot\|_T)={X}\otimes_\epsilon \ell_2^2$.\smallskip

Another important reasonable norm on the complexification that is related to tensorial norms, named after Bochnak as the \emph{Bochnak norm}, is defined as
$$\|x+iy\|_B:=\inf\left\{\sum\limits_k|\lambda_k|\,\|x_k\|: x+iy=\sum\limits_k\lambda_k\otimes x_k\in {X}\otimes_{\mathbb{R}} \mathbb{C}\right\}.$$
Observe that the Bochnak norm is nothing but the projective tensor norm in ${X}\otimes \ell_2^2$ or, alternatively, $({X}_c,\|\cdot\|_B)={X}\otimes_\pi \ell_2^2$.\smallskip

Additional reasonable norms on the complexification are obtained by
 $$
 \|x+iy\|_{(p)}=2^{\min\{1/2-1/p,0\}}\sup_{t\in [0,2\pi]}(\| x \cos t -y \sin t\|^p+\|x \sin t +y \cos t\|^p)^\frac{1}{p}
 $$
for each $p\in [1,\infty)$. In the special case in which $p=2$, $\|\cdot\|_{(2)}$ is the so-called \emph{Lindenstrauss--Tzafriri norm}\label{page:LT} (see \cite{LT}) and will be denoted by $\|\cdot\|_{LT}$. Therefore,
$$\|x+iy\|_{LT}:=\sup_{t\in[0,2\pi]}\left(\| x \sin t + y \cos t\|^2+\|x \cos t- y \sin t \|^2\right)^{1/2}.$$
The norm $\|\cdot\|_{(1)}$ was employed by Alexiewicz and Orlicz \cite{AO1953}.\smallskip

The following complexification norms have been considered by Kirwan \cite{Kirwan} (we keep author's original notation):
 \begin{align*}
 \gamma_p(x+iy)&=c_p\left(\frac{1}{2\pi}\int_0^{2\pi}\|x\cot t-y\sin t\|^pdt\right)^\frac{1}{p},\\
 \gamma_{p,q}(x+iy)&=c_{p,q}\left(\frac{1}{2\pi}\int_0^{2\pi}\left(\|x\cos t-y\sin t\|^p+\|x\sin t+y\cos t\|^p\right)^\frac{q}{p}dt\right)^\frac{1}{q},
 \end{align*}
where $1\leq p,q<\infty$ and
 \begin{align*}
 c_p&=\left(\frac{1}{2\pi}\int_0^{2\pi}|\cos t|^p\right)^{-\frac{1}{p}},\\
 c_{p,q}&=\left(\frac{1}{2\pi}\int_0^{2\pi}\left(|\cos t|^p+|\sin t|^p\right)^\frac{q}{p}\right)^{-\frac{1}{q}}.
 \end{align*}
Interestingly, the reasonable complexification norms $\|\cdot\|_T$, $\|\cdot\|_B$, and $\|\cdot\|_{LT}$ complexify, in a natural way, the real versions of the spaces $\ell_\infty$, $\ell_1$, and $\ell_2,$ respectively (see \cite{Kirwan,MUN}), that is, the complexifications of $\ell_\infty({\mathbb R})$, $\ell_1({\mathbb R})$, and $\ell_2({\mathbb R})$ endowed with the norms $\|\cdot\|_T$, $\|\cdot\|_B$, and $\|\cdot\|_{LT}$, are $\ell_\infty({\mathbb C})$, $\ell_1({\mathbb C})$, and $\ell_2({\mathbb C}),$ respectively. \smallskip

Finally, one may observe that if $X$ is a real Banach lattice, then the norm $$\|x+iy\|=\|\,\,|x+iy|\,\,\|$$  makes $X_c$ into a complex Banach lattice, where $|x+iy|:=\sup_{t\in[0,2\pi]}|x\cos t +y\sin t|$ is the extension of the modulus function $|\cdot|$ of $X$ to $X_c$. It is shown in \cite{NEE} that this norm is induced by the so-called $l$-norm on $X\otimes \ell_2^2$.

\subsection{Regular Banach spaces and complex strictly convex complexifications}\ \smallskip

A complex Banach space $\mathscr{X}$ is called \emph{regular} if it is isomorphic to the complexification of a real Banach subspace ${Y}$ of $\mathscr{X}_r$ and is equipped with a reasonable norm.\smallskip

For example, the complexification of any real Banach space $({Y},\|\cdot\|)$ endowed with the norm $$\|y+iz\|:=\sup_{t\in[0,2\pi]}\left(\|y\sin t+z\cos t\|+\|y\cos t-z\sin t\|\right) \quad (y, z\in {Y})$$ is regular. Some other examples are given by the complex spaces $\ell_p(\mathbb{N}, \mathbb{C})$ and $\mathcal{L}_p([0, 1],\mathbb{C}),\,\, 1 \leq p < \infty$ equipped with their usual norms (see \cite{FG}).\smallskip

A complex normed space $\mathscr{X}$ is said to be \emph{complex-strictly convex} if the inequality $$ \frac{1}{2\pi} \int_0^{2 \pi} \| x + e^{i t } y\| dt >1$$ holds for each $x,y\in \mathscr{X}$ with $\|x\|=1$ and $y\neq 0$  (this is formally stronger than the notion of strict c-convexity but it is actually equivalent to it (cf. \cite[Theorem 2]{DowHuMup96}). Recall that a complex normed space $\mathscr{X}$  is  called   \emph{strictly c-convex} if for all $x, y  \in \mathscr{X}$ with $\|x||=1$, $\|x + \lambda y\|<1\,\,(|\lambda|\leq 1)$ implies  $y= 0$; see \cite[Definition 1]{Glob75}). \smallskip

The question for what real Banach spaces there is a complex-strictly convex norm on the corresponding algebraic complexification, remains open. Kadets and Kellerman solved this problem for all separable spaces, and they also showed that the statement holds for many, but not all non-separable ones.

\begin{theorem}\label{thm KadetsKellerman}\cite[Theorem 1]{KadKell2000} Every  separable real normed space $X$ admits a complex strictly convex complexification.
\end{theorem}

The conclusion in the previous theorem holds for some classes of nonseparable spaces; for example, for spaces $X$ with 1-norming separable subspases in $X^*$; however, certain spaces of the form $\ell_{\infty}(\Gamma)$ admit no complex strictly convex complexifications.\smallskip

Some open questions remain open. We first recall that a complex Banach space $\mathscr{X}$ is \emph{complex locally uniformly convex}, if for every $x \in \mathscr{X}$ with $\|x\|= 1$ and every sequence $(y_n) \subset \mathscr{X} \setminus \{0\}$ if $\lim_{n \to \infty}\frac{1}{2\pi}\int_0^{2\pi} \|x + e^{it}y_n\| dt = 1$, then $\lim_{n \to \infty} \|y_n\| = 0$. \smallskip

It would be desirable to characterize those real Banach spaces, on which every equivalent norm can be complexified to a complex-strictly convex one. It is not known, in particular, whether $\ell_\infty$ has this property.\smallskip

Another open question asks whether the statement of Theorem \ref{thm KadetsKellerman} remains true, if one substitutes the complex-strict convexity by complex locally uniform convexity.

\subsection{Complexification of Hilbert spaces}\label{subsec:Complexification of Hilbert spaces} \smallskip

Let ${H}$ be a real Hilbert space. The algebraic complexification, ${H}_c={H}+i {H},$ of $H$ can be equipped with a natural inner product structure via the assignment
	$$
	\langle x+iy,x'+iy'\rangle:=(\langle x,x'\rangle+\langle y,y'\rangle)+i(\langle y,x'\rangle-\langle x,y'\rangle).
	$$
In this case, the identity $$\|x+iy\|^2=\|x\|^2+\|y\|^2$$ holds for all $x, y\in {H}$. However, ${H}$ is not orthogonal to $ i{H}$ in the Hilbert space ${H}_c$.\smallskip

{At this point it is interesting to observe that if $X$ is any real Banach space and $\|\cdot\|_\nu$ is any reasonable complexification norm on $X_c$, then for every $x+iy\in X_c$ and $t\in{\mathbb R}$,
	\begin{align*}
	\|x+iy\|_\nu&=\|e^{it}(x+iy)\|_\nu=\|x\cos t-y\sin t+i(x\sin t+y\cos t)\|_\nu\\
	&\leq \|x\sin t+y\cos t\|+\|x\cos t-y\sin t\|.
	\end{align*}	
Hence, if
	$$
	B(x+iy)=\inf_{t\in[0,2\pi]}\left(\|x\sin t+y\cos t\|+\|x\cos t-y\sin t\|\right),
	$$
then $\|x+iy\|_\nu\leq B(x+iy)$. It is precisely when $X$ is a real Hilbert space and $\|\cdot\|_\nu$ is the Bochnak norm that the latter inequality is in fact an equality.
As a matter of fact, if $H$ is a real Hilbert space} the Bochnack norm on ${H}_c$ can be represented by a simpler formulas as
	$$
	\|x+iy\|_B=\inf_{t\in[0,2\pi]}\left(\|x\sin t+y\cos t\|+\|x\cos t-y\sin t\|\right)
	$$
and
	$$
	\|x+iy\|_B=\left(\|x\|^2+\|y\|^2+2\left(\|x||^2\,\|y\|^2-\langle x,y\rangle^2\right)^{1/2}\right)^{1/2},
	$$
where $x, y\in {H}$. {The previous two formulas were proved in \cite[Proposition 3]{FM}. We provide below an alternative proof of the last formula communicated to the authors of \cite{FM} by the anonymous referee in his/her report. First notice that $(H_c,\|\cdot\|_B)$ is nothing but the projective tensor $H\otimes_\pi \ell_2^2$. Since $\text{span}\{x,y\}\otimes_\pi \ell_2^2$ is 1-complemented in $H\otimes_\pi \ell_2^2$, it suffices to check the formula for $\ell_2^2\otimes_\pi\ell_2^2$. Identifying $\ell_2^2\otimes_\pi\ell_2^2$ with the Schatten 1-class $S_1(\ell_2^2)$, the norm of $x+iy$ with $x=(x_1,x_2)$ and $y=(y_1,y_2)$ is given by $\|T\|_1$ where
	$
	T=\left(\begin{array}{cc}
	x_1 & y_1\\
	x_2 & y_2
	\end{array}
	\right).
	$
It is well-known that $\|T\|_1=\text{trace}(T^*T)^{1/2}$. Since
	$
	T^*T=\left(\begin{array}{cc}
	\|x\|^2 & \langle x,y\rangle\\
	\langle x,y\rangle & \|y\|^2
	\end{array}
	\right),
	$
we have that $\|T\|_1=\sqrt{\lambda_1}+\sqrt{\lambda_2}$ where $\lambda_1$ and $\lambda_2$ are the eigenvalues of $T^*T$. It is elementary to show that
\begin{align*}
\lambda_1&=A+\sqrt{B^2+C^2},\\
\lambda_2&=A-\sqrt{B^2+C^2},
\end{align*}
where $A=\frac{1}{2}(\|x\|^2+\|y\|^2)$, $B=\frac{1}{2}(\|x\|^2-\|y\|^2)$ and $C=\langle x,y\rangle$, concluding the proof.

On the other hand, an argument based on elementary calculus leads to
	$$
	B(x+iy)=\left(\|x\|^2+\|y\|^2+2\left(\|x||^2\,\|y\|^2-\langle x,y\rangle^2\right)^{1/2}\right)^{1/2}.
	$$
We reproduce the proof of the latter formula found in \cite[Proposition 1]{FM}.	
For fixed vectors $x,y\in H$, we define $\phi_{xy}:[0,2\pi]\rightarrow
{\mathbb R}$ by
	$$
	\phi_{xy}(t):=\|x\cos t-y\sin t\|+\|x\sin t+y\cos t\|,
	$$
for every $t\in [0,2\pi]$. Then
	$$
	\phi_{xy}(t)=\sqrt{A+B\cos 2t-C\sin 2t}+\sqrt{A-B\cos 2t+
	C\sin 2t},
	$$
for every $t\in [0,2\pi]$.
Using elementary calculus it is easily seen that the mapping
$\phi_{xy}$ attains its infimum at a point $t_0\in [0,2\pi]$
such that $B\cos 2t_0-C\sin 2t_0=\pm\sqrt{B^2+C^2}$. Therefore
\begin{align*}
B(x+iy)^2&=\left[\inf_{t\in [0,2\pi]}\phi_{xy}(t)\right]^2\nonumber\\
&=\bigg\{\sqrt{A+\sqrt{B^2+C^2}}+\sqrt{A-\sqrt{B^2+C^2}}\bigg\}^2\label{ABC}\\
&=2A+2\sqrt{A^2-B^2-C^2}\nonumber\\
&=\|x\|^2+\|y\|^2+2\big\{\|x\|^2\|y\|^2-\langle x,y\rangle^2\big\}^{1/2}.\nonumber
\end{align*}

To finish this series of comments on the connection between the Bochnak norm and the mapping $B$ in the context of Hilbert spaces it can be proved (see \cite[Theorem 6]{FM}) that if $X$ is a real Banach space and $\|\cdot\|_\mu$ is a reasonable complexificaton norm on $X_c$ such that $\|x+iy\|_\nu=B(x+iy)$ for all $x+iy\in X_c$, then $X$ is a Hilbert space.}\smallskip

In our seeking for similarities and differences between real and complex spaces, we note that in a real Hilbert space ${H}$, the function $f(x)=\|x\|^2,$ ($x\in {H}$) is Fr\'{e}chet differentiable at every point and $f'(x) y = 2\langle x,y\rangle$, while this statement does not hold for complex Hilbert spaces.\smallskip

The conjugate linearity of the inner product of a complex Hilbert spaces in the second variable also produces differences in the adjoint of a bounded linear operator between real or complex Hilbert spaces. For each bounded linear operator $T: \mathscr{H}\to \mathscr{K}$ between Hilbert spaces, there exists a unique bounded linear operator $T^*: \mathscr{K} \to \mathscr{H}$ satisfying the identity $$\langle Tx,y\rangle=\langle x,T^*y\rangle \quad \hbox{ for all } x\in \mathscr{H} \hbox{ and } y\in \mathscr{K}.$$ In the setting of real Hilbert spaces, the assignment $T\mapsto T^*$ is an isometric (real) linear isomorphism, while in the setting of complex Hilbert spaces, it is a conjugate linear isometric isomorphism (i.e., $(\lambda T)^*=\overline{\lambda} T^*$).\smallskip

Along this note, the symbol $\mathcal{B}(\mathscr{X},\mathscr{Y})$ will denote the complex (respectively, real) Banach space of all bounded linear operators $T: \mathscr{X}\to \mathscr{Y}$ between complex (respectively, real) Banach spaces endowed by the operator norm $$\|T\|:=\sup\{\|Tx\|: \|x\|=1\}. $$ We write $\mathcal{B}(\mathscr{X})$ for the space $\mathcal{B}(\mathscr{X},\mathscr{X})$. According to this notation, the symbol $\mathscr{X}^*$ will stand for $\mathcal{B}(\mathscr{X}, \mathbb{C})$ (respectively, $\mathcal{B}(X, \mathbb{R})$) when $\mathscr{X}$ is a complex (respectively, real) Banach space.\smallskip

In the setting of complex Hilbert space operators, given an operator $T\in \mathcal{B}(\mathscr{H})$, the condition $\langle T\xi,\xi\rangle=0,$ for all $\xi\in \mathscr{H}$ implies that $T$ is identically $0$. The proof essentially applies complex scalars; actually, such a result does not hold for bounded linear operators on a real Hilbert space. Take, for example, ${H}=\ell_2^2$ and $T=\left(\begin{matrix} 0 & 1 \\ -1 & 0\end{matrix}\right)$ as the matrix of 90 degrees clockwise rotation.\smallskip

Let $\mathscr{H}$ be a complex Hilbert space. An operator $T\in \mathcal{B}(\mathscr{H})$ is called (\emph{numerically}) \emph{positive} if $\langle T\xi,\xi\rangle \geq 0$ for every $\xi \in \mathscr{H}$. It immediately follows that $T^*=T$, that is, $T$ is self-adjoint. However, if ${H}$ is a real Hilbert space, then the positivity of $T$ in the above terms does not entail that $T$ is self-adjoint. Consider, for example, ${H}=\ell_2^2$ and $T=\left(\begin{matrix} 1 & 1 \\ -1 & 1\end{matrix}\right)$. Thus, it is reasonable that in the definition of positivity for bounded real linear operators on a real Hilbert space, we add the self-adjointness condition of $T$ (see Subsection \ref{subsec: KadisonSchwarz inequ} for a deeper discussion in the setting of real C$^*$-algebras).\smallskip

We continue in the complex setting, and suppose that $T\in \mathcal{B}(\mathscr{H})$ is a bounded linear operator for which there exists an orthonormal basis $\{e_j: j\in \Lambda\}$ of $\mathscr{H}$ consisting of eigenvectors of $T$. Then for each $j$, it follows that $Te_j=\lambda_j e_j$ for some scalar $\lambda_j$. It is easy to see that $T^* e_j =\overline{\lambda_j} e_j$ for every $j\in \Lambda$. If ${H}$ is a real Hilbert space, then each $\lambda_j$ is a real number, and so $T=T^*$. In the complex case, we can only conclude that $T$ is normal (i.e., $T T^* = T^* T$). In other words, a real linear combination of mutually orthogonal projections always gives a self-adjoint operator, while if we consider complex linear combinations, then the result operator is only normal. \smallskip

In the setting of complex Hilbert spaces, a linear functional $\varphi: \mathcal{B}(\mathscr{H})\to \mathbb{C}$ is called \emph{positive} if $\varphi(T^*T)\geq 0$ for every $T\in \mathcal{B}(\mathscr{H})$. Such a functional is always self-adjoint in the sense that $\varphi(T^*)=\overline{\varphi(T)}$ for all $T\in \mathcal{B}(\mathscr{H})$. However, in the framework of real Hilbert spaces, this assertion is not valid, in general. For example, assume that ${H}=\ell_2^2$ and that the linear functional $\varphi: \mathcal{M}_2(\mathbb{R})\to \mathbb{R}$ defined by $\varphi\left(\begin{matrix} \alpha & \beta \\ \gamma & \delta \end{matrix}\right)=\alpha+\beta+\delta$ is positive but not self-adjoint (see \cite{RUA2}). In Subsection \ref{subsec: KadisonSchwarz inequ}, we include a detailed discussion on positive linear functionals on real C$^*$-algebras. \smallskip

The ``complex plank problem''\label{page:plank} asks whether for a finite sequence $(\xi_k)_{k=1}^n$ of unit vectors in a complex Hilbert space $\mathscr{H}$, there exists a unit vector $\xi\in \mathscr{H}$ such that $|\langle \xi,\xi_k\rangle|\geq 1/\sqrt{n}$ for each $k = 1, \ldots, n$ (see \cite{BAA}). This fact is not true, in general, in the setting of real Hilbert spaces. For example, we can consider the usual real Hilbert space $\ell_2^2$ and the unit vectors $\xi_1,\ldots, \xi_{2n}$ uniformly distributed around the circle. Then for each unit vector $\xi\in\ell_2^2$, there exists a vector $\xi_k$ for some $1\leq k\leq 2n$ such that $|\langle \xi,\xi_k\rangle|\leq \sin(\frac{\pi}{2n})\leq \frac{\pi}{2n}$ (see \cite[p. 706]{RT}). This topic will be addressed with more detail in Section \ref{PlankProblem}.\smallskip

We can conclude this subsection with an additional example illustrating the previous notions. For a complex Hilbert space $\mathscr{H}$, the subset $\mathcal{B}(\mathscr{H})_{sa}$ of all self-adjoint or hermitian operators on $\mathscr{H}$ (i.e., all $T\in \mathcal{B}(\mathscr{H})$ with $T^* = T$) is a norm closed real subspace of $\mathcal{B}(\mathscr{H})$. It is well known that the algebraic complexification of $\mathcal{B}(\mathscr{H})_{sa}$ is precisely the whole of $\mathcal{B}(\mathscr{H})$.\smallskip

According to the usual terminology (see, for example, \cite[\S 9, Definition 7]{BDnr}), for an operator $T$ in $\mathcal{B}(\mathscr{H}),$ where $\mathscr{H}$ is a real or complex Hilbert space, the (\emph{spatial}) \emph{numerical range} of $T$ is defined as the set $$W(T)= \{ \langle T\xi ,\xi\rangle : \xi \mbox{~is a unit vector in~} \mathscr{H}\}.$$ The \emph{numerical radius} of $T$ is defined by $$w(T)=\sup\{|\lambda|: \lambda\in W(T)\}.$$ For a complex Hilbert space $\mathscr{H}$, a remarkable result by Sinclair asserts that for each $T\in \mathcal{B}(\mathscr{H})_{sa}$, we have $w(T) = \|T\|$ (see \cite{Sinc71,BDnr}). A classic result in operator theory (see \cite[pp. 116--117]{HAL}) assures that \begin{equation}\label{eq Bohoneblust Karlin Hilbert} \frac{1}{2}\|T\| \leq w(T) \leq \|T\|\quad \hbox{ for all } T\in \mathcal{B}(\mathscr{H}).
\end{equation} As remarked by Ili\v{s}evic et al. \cite[Example 3.15]{IliKuzLiPoon}, the Taylor complexification norm on $\mathcal{B}(\mathscr{H})$ of the restriction of the spectral or operator norm on $\mathcal{B}(\mathscr{H})_{sa}$ is precisely the numerical radius $w(\cdot)$, that is, $$w(T) = \|T\|_{T} = \|H +i K\|\sup_{t\in[0,2\pi]}\| H \cos t - K \sin t\|$$ for all $T = H+i K$ in $\mathcal{B}(\mathscr{H})$ with $H,K\in \mathcal{B}(\mathscr{H})_{sa}$.\smallskip

However, if $H$ is a real Hilbert space, then the left-hand side inequality in \eqref{eq Bohoneblust Karlin Hilbert} may be not true. For example, let $H=\ell_2^2$ and let $T=\left(\begin{matrix} 0 & -1 \\ 1 & 0\end{matrix}\right)$ to get $w(T)=0$ and $\|T\|=1$. In the setting of real Banach spaces, Lumer \cite[Theorem 1]{LUM} proved a deep result showing the existence of positive constants $c_1$ and $c_2$ such that $$\|T\|\leq c_1 w(T) +c_2 w(T^2)^{1/2},$$ for every $T\in \mathcal{B}({X})$. In the case that ${H}$ is a real Hilbert space with dim$(H)>1$, the inequality $$\|T\| \leq 2 w(T) + w(T^2)^{\frac12},$$ holds for all $T\in \mathcal{B}({H})$, and $2$ and $1$ are the best possible constants (see \cite[Theorem 10]{LUM}). Such inequalities were discussed in \cite{LUM}.

\subsection{Schauder Basis}\ \smallskip

Several key results hold only for complex Banach (Hilbert) spaces but not for real Banach (Hilbert) spaces since they depend on complex analysis techniques. Furthermore, some of the facts in the setting of complex Banach (Hilbert) spaces can be proved for real spaces under some extra conditions or changing some hypotheses.\smallskip


Let $\mathscr{X}$ be a real or complex Banach space. A (Schauder) basis of $\mathscr{X}$ is a sequence $(x_n)_{n=1}^\infty$ such that for every vector $x\in \mathscr{X}$, there exists a unique sequence $(\lambda_n)$ of scalars such that $x=\sum\limits_{n=1}^\infty\lambda_nx_n$ in the norm topology on $\mathscr{X}$. A Schauder basis is called a \emph{$1$-unconditional basis} if for every sequence of scalars $(\lambda_n)$ and every sequence of scalars $(\varepsilon_n)$ in the closed unit ball of the corresponding field, it holds that $\left\|\sum\limits_{n=1}^\infty\varepsilon_n\lambda_nx_n\right\|\leq \left\|\sum\limits_{n=1}^\infty\lambda_nx_n\right\|$ (see \cite{LT}).\smallskip

The existence of an unconditional Schauder basis in a real Banach space can be employed to define a norm on its algebraic complexification (see \cite{Kirwan}). Namely, if $\{e_n:n\in{\mathbb N}\}$ is an unconditional Schauder basis of the real Banach space ${X}$, then
 $$
 \widetilde{{X}}=\left\{(\lambda_n)\in{\mathbb C}^{\mathbb N}:\sum\limits_{n=1}^\infty |\lambda_n| e_n\text{ converges in ${X}$}\right\}
 $$
is a complex linear space. Moreover, every $(\lambda_n)\in\widetilde{{X}}$ with $\lambda_n=a_n+ib_n$ can be viewed as a vector $x+iy\in {X}\oplus i {X}$, where
 $$
 x=\sum\limits_{n=1}^\infty a_ne_n\quad\text{and}\quad y=\sum\limits_{n=1}^\infty b_ne_n.
 $$
It is easy to furnish $\widetilde{{X}}$ with a reasonable complexification norm by defining
 $$
 \|(\lambda_n)\|_{unc}=\left\|\sum\limits_{n=1}^\infty |\lambda_n|e_n\right\|_{{X}}.
 $$
It was proved in \cite[Prop 2.13]{Kirwan} that if ${\mathcal B}=\{e_n:n\in{\mathbb N}\}$ is a $1$-unconditional monotone Schauder basis for the real Banach space ${X}$ (i.e., $\left\|\sum\limits_{n=1}^\infty r_ne_n\right\|_{{X}}\leq \left\|\sum\limits_{n=1}^\infty s_ne_n\right\|_{{X}}$ whenever $r_n\leq s_n$ for all $n\in{\mathbb N}$), then ${\mathcal B}$ is also a $1$-unconditional monotone Schauder basis for $\widetilde{{X}}$.\smallskip

The different behavior of real and complex Banach spaces will be again contrasted. Kalton and Wood \cite[Theorem 6.1]{KW} proved that any two $1$-unconditional bases $(x_n)$ and $(y_n)$ in a complex Banach space $\mathscr{X}$ are isometrically equivalent in the sense that there is a permutation $\sigma$ on $\mathbb{N}$ such that $y_{\sigma(n)} = \lambda_n x_n$ for all $n$, where the $\lambda_n$'s are scalars of modulus $1$. However, Lacey and Wojtaszczyk \cite{LW} showed that this conclusion does not hold for real $L_p$-spaces (see also \cite{RAN}).


\subsection{Extension of linear operators to the complexifications}\ \smallskip

Let ${X}$ and ${Y}$ be two real Banach spaces and let ${X}_c$ and ${Y}_c$ be two arbitrary complexifications of $X$ and $Y$, respectively. If $S\in \mathcal{B}({X},{Y})$, then the operator $S_c: X_c\to Y_c$, defined in Subsection \ref{subsec: algebraic complexification and operators} by $S_c(x+iy)=S(x)+i S(y)$, is bounded and there is a constant $m$ such that $\|S\|\leq \|S_c\|\leq m \|S\|$. As seen in \eqref{eq identification of L(X,Y) is the corresponding complex linear extensions}, the real linear subspace $\mathcal{B}(X_c,Y_c)_{sym} = \{ S_c : S\in \mathcal{B}({X},{Y})\}$ is a real linear subspace of $\mathcal{B}({X}_c,{Y}_c)$, which can be algebraically identified with $\mathcal{B}(X,Y)$. Moreover, each $T\in \mathcal{B}(X_c,Y_c)$ can be written in the form $T = T_1 + i T_2$, where $T_1 = \frac{T + \overline{T}}{2}$, $T_2 = \frac{T - \overline{T}}{2 i }\in \mathcal{B}(X_c,Y_c)_{sym}.$ Thus, algebraically $\mathcal{B}(X_c,Y_c) = \mathcal{B}(X,Y) + \mathcal{B}(X,Y)$. But, from the analytic point of view, we have two norms on $\mathcal{B}(X_c,Y_c)_{sym}$, one inherited from $\mathcal{B}({X}_c,{Y}_c)$ and another one obtained when it is identified with $\mathcal{B}(X,Y)$ as a real Banach space. Since, in general, $\|S_c\|\neq \|S\|$ (cf., see Example \ref{example second March}), the Banach space $\mathcal{B}(X_c, Y_c)$ cannot be identified with any complexification of $\mathcal{B}(X_c,Y_c)_{sym} \equiv \mathcal{B}({X},{Y})$ with the operator norm of the latter space.\smallskip

Despite the difficulties, in the setting of real Hilbert spaces, the complex Banach space $\mathcal{B}({H}_c)$ actually is the complexification of $\mathcal{B}({H})$, where $H_c$ is the complex Hilbert space obtained by complexifying $H$ (see Subsection \ref{subsec:Complexification of Hilbert spaces}). Furthermore, $\mathcal{B}({H})_c$ can be identified (completely isometrically) with a subspace of $\mathcal{M}_2(\mathcal{B}({H}))$ via $T+iS\mapsto \left(\begin{matrix} T & -S \\ S & T\end{matrix}\right)$; in particular, $\mathbb{C}$ can be identified with $$\left\{\left(\begin{matrix} t & -s \\s & t \end{matrix}\right)\in \mathcal{B}(\ell_2^2): t, s\in \mathbb{R}\right\}$$
as a real subspace of $\mathcal{M}_2(\mathbb{R})$ (see \cite[page 1051]{RUA}). We will revisit this construction in Subsection \ref{subsec: KadisonSchwarz inequ}.\smallskip

It is not hard to show that if we use the Taylor complexification in both real Banach spaces ${X}$ and ${Y}$, then $\|S_c\| = \|S\|$ for every $\mathcal{B}(X,Y)$ (with respect to the operator norm given by the Taylor complexification norm). The same holds when we employ the Lindenstruass--Tzafriri norm, the Bochnak norm, or the $(p)$ norms (see \cite{LI,MUN}). This particularly interesting property motivates the concept of \emph{natural complexification procedure} (see \cite{MUN}). A natural complexification procedure $\nu$ is a way to assign to each real Banach space $X$, a reasonable complexification norm $\|\cdot\|_\nu$ in such a manner that if ${X}$ and ${Y}$ are arbitrary real Banach spaces and $S\in \mathcal{B}({X},{Y})$, then $\|S_c\|_\nu = \|S\|$, where $\|S_c\|_\nu$ must be understood as the operator norm of $S_c$ as an operator between $({X}_c,\|\cdot\|_\nu)$ and $({Y}_c,\|\cdot\|_\nu)$. That is, $\mathcal{B}({X},{Y})$ and $\mathcal{B}((X_c,\|\cdot\|_\nu),(Y_c,\|\cdot\|_\nu))_{sym}$ are isometrically isomorphic as real Banach spaces.\smallskip


{\subsection{Extension of operators and injectivity: Real vs complex cases}\ \smallskip

A real or complex Banach space $X$ is said to be injective if for every Banach space $Z$ and every subspace $Y$ of $Z$, each operator $T\in {\mathcal B}(Y,X)$ admits an extension ${\overline T}\in{\mathcal B}(Z,X)$. Additionally, if $\lambda>0$ then we say that $X$ is $\lambda$-injective if the extension ${\overline T}$ can be chosen so that $\|{\overline T}\|\leq \lambda \|T\|$. Obviously, if $X$ is 1-injective then any $T\in {\mathcal B}(Y,X)$ can be extended to a $\|{\overline T}\|\leq \lambda \|T\|$ with preservation of its norm, i.e., $\|T\|=\|{\overline T}\|$. The space $\ell_\infty$ is a good example of an injective space. As a matter of fact, $X$ is injective if and only if it is a $C(K)$ space with $K$ being an extremely disconnected compact space as proved by Nachbin \cite{N1950} and Kelley \cite{K1952}. The study of injective Banach spaces has attracted the attention of many researchers since at least the 1940's, generating a vast literature. For a comprehensive global view on the topic we recommend \cite{A2016}, where the interested reader will be able to check that most of the results proved for real injective Banach spaces can be also established in a complex setting without much difficulty. However there is one significant consideration between the real and complex cases that must be underlined. This difference has been detected within the context of separably injective Banach spaces. Recall that a real or complex Banach space $X$ is separable injective if for every separable Banach space $Z$ and every subspace $Y$ of $Z$, each operator $T\in {\mathcal B}(Y,X)$ admits an extension ${\overline T}\in{\mathcal B}(Z,X)$. The concept of $\lambda$-separable injectivity is defined similarly. The spaces $c$ and $c_0$ are 2-separably injective (see \cite{P1940} and \cite{S1941}). Actually $c_0$ is the only separable Banach space that is separably injective (see \cite{Z1977}). The above mentioned difference between real and complex separably injective spaces is found in the following characterization of real 1-separably injective spaces (see for instance\cite[Proposition 2.30]{A2016}).
\begin{proposition}\label{prop:characterization}
	A real Banach space $X$ is 1-separably injective if and only if every countable family of mutually intersecting balls has nonempty intersection.
\end{proposition}
To translate the previous characterization into a complex setting a new property must be defined. We say that a family of balls $\{B(x_\xi; r_\xi)\}_\xi$ in a Banach space $X$ over ${\mathbb K}$ is \emph{weakly intersecting} if for every norm
one $f\in X^*$ the balls $\{B(f(x_\xi),r_\xi)\}_\xi$ have nonempty intersection. The previous property was introduced in \cite{H1973}. Using this terminology, the complex analog of Proposition \ref{prop:characterization} would be
\begin{proposition}
	A complex Banach space $X$ is 1-separably injective if and only if every countable family of weakly intersecting balls has nonempty intersection.
\end{proposition}
}


\subsection{Spectrum of an operator}\label{first results spectrum}\ \smallskip

Let ${X}$ be a real Banach space. The \emph{{\rm(}real{\rm)} spectrum} of an operator $T\in \mathcal{B}({X})$ may be defined as the set $$\{\lambda\in\mathbb{R}: T-\lambda I_{X} \mbox{~is not invertible in ~} \mathcal{B}({X})\},$$ where $I_{X}$ denotes the identity operator on ${X}$. It is well known that this definition has several handicaps. For example, the spectrum of a bounded linear operator $T$ on a real Banach space ${X}$ given by this definition may be empty, such as, the case where ${X}=\ell_2^2$ and $T=\left(\begin{matrix} 1 & 2 \\ -2 & 0\end{matrix}\right)$. Thus matrices with real entries may have complex eigenvalues.\smallskip

Thus it is more appropriate to define the spectrum of $T\in \mathcal{B}({X})$ as the spectrum of $T_c$ in $\mathcal{B}(X_c)$, that is, the set
$${\rm sp}(T) :=\{\lambda\in\mathbb{C}: T_c-\lambda I_{X_c} \mbox{~is not invertible in ~} \mathcal{B}({X}_c)\}.$$

An elementary \emph{spectral theorem} affirming that the identity ${\rm sp}(p(T))=p({\rm sp}(T))$ is true for any bounded linear operator $T$ acting on a complex Hilbert space $\mathscr{H}$ and any polynomial $p$ with complex coefficients, can be obtained because a polynomial with complex coefficients is a product of polynomials of degree $1$. However, the fundamental theorem of algebra fails to be true in the context of polynomials with real coefficients; therefore, the spectral theorem stated above does not hold in this setting for the suggested real spectrum. For example, let ${X}=\ell_2^2$, let $T=\left(\begin{matrix} 0 & 1 \\ -1 & 0\end{matrix}\right)$, and let $p(t)=t^2$.\smallskip

The reader is referred to Section \ref{sec: Banach algebras}, where a more detailed study on the similarities and differences between real and complex Banach algebras is conducted.


\subsection{Invariant subspaces}\ \smallskip

A subspace $\mathscr{M}$ of a real or complex Banach space $\mathscr{X}$ is said to be \emph{invariant under an operator} $T\in \mathcal{B}(\mathscr{X})$ if $T(\mathscr{M})\subseteq \mathscr{M}$. The subspace $\mathscr{M}$ is called \emph{nontrivial} if $\{0\}\neq \mathscr{M}\neq \mathscr{X}$. If $\mathscr{M}$ is invariant under every bounded linear operator commuting with $T$, then it is called \emph{hyperinvariant} for $T$.\smallskip

The problem of whether every bounded operator $T$ on a complex (or real) Banach space $\mathscr{X}$ possesses a nontrivial closed subspace $M$ which is $T$-invariant has been a long standing problem in functional analysis. Enflo provided a counterexample to this question for Banach spaces in 1976, although due to the high complexity of Enflo's construction, his 100 page long paper was not published until 1987, \cite{Enflo1987} (see, also, e.g., \cite{Araujo2020,CP2011,LOM,Read1985}). The problem still remains open for separable Hilbert spaces.\smallskip

Let ${X}$ be an infinite-dimensional separable real Banach space. If $T\in \mathcal{B}({X})$ has a nontrivial invariant closed subspace $\mathscr{M}$, then $\mathscr{M}+i\mathscr{M}$ is a nontrivial closed subspace of ${X}_c$ invariant under $T_c$. If $T$ has no nontrivial closed
invariant subspaces, then it is an interesting question to ask whether the same is true for $T_c\in \mathcal{B}({X}_c)$ (see \cite[Conjecture 3]{ABR}).\smallskip

 Lomonosov \cite{LOM} proved that if a nonscalar bounded linear operator $T$ on a complex Banach space commutes with a nonzero compact linear operator, then $T$ admits a nontrivial hyperinvariant closed subspace. In his proof, Lomonosov used an essential property, that is, bounded linear operators on a finite-dimensional complex space have eigenvalues.\smallskip

 Hooker \cite[p. 132]{HOO} provided, among other results, a counterexample to Lomonosov's result in the real setting. The linear isometry $$T(x_1,y_1,x_2,y_2,\ldots)=(-y_1,x_1,-y_2,x_2,\ldots)$$ on the real Hilbert space $\ell_2$ has no nontrivial closed hyperinvariant subspaces.\smallskip

It is worth noting that for each real Banach space ${X}$ and each nonscalar operator $T\in \mathcal{B}({X})$ commuting with a nonzero compact linear operator on ${X}$, the following statements are equivalent (see \cite{SIR}):
\begin{enumerate}[$\bullet$]
 \item[(i)] $T$ has a nontrivial closed hyperinvariant subspace;
 \item[(ii)] For each pair of real numbers $\alpha$ and $\beta$ with $ \beta\neq 0$, we have $(\alpha-T)^2+\beta^2\ne 0$.
\end{enumerate}

\subsection{Dual}\ \smallskip

Let $\mathscr{X}$ be a complex Banach space. For any $f$ in the dual space, $\mathscr{X}^*$, of $\mathscr{X}$, consider the linear functional ${\rm Re} \varphi : \mathscr{X}_r\to \mathbb{R}$ given by $$({\rm Re} \varphi)(x):={\rm Re}(\varphi(x))\quad (x\in \mathscr{X}_r).$$ A classic result in functional analysis affirms that the assignment $\varphi\mapsto {\rm Re} \varphi$ provides an isometric real linear isomorphism from $\left(\mathscr{X}^*\right)_r$ onto $(\mathscr{X}_r)^*$.\smallskip

On the other hand, if ${X}$ is a real Banach space and $\phi_1,\phi_2\in {X}^*$, then the mapping defined by
 $$
 \widetilde{(\phi_1+i \phi_2)} (x+iy)=\phi_1 (x)-\phi_2(y)+ i \left(\phi_2(x)+\phi_1(y)\right)
 $$
is a linear functional in $({X}_c)^*$. As a matter of fact,
 $$
 \Psi:({X}^*)_c\ni \phi_1+i \phi_2\mapsto \widetilde{(\phi_1+i \phi_2)}\in ({X}_c)^*
 $$
is a natural isomorphism between $(({X}^*)_c,\|\cdot\|_\nu)$ and $(({X}_c)^*,\|\cdot\|_\nu)$ for any natural complexification procedure $\nu$ (see, for instance, \cite{MUN}). Hence if ${X}_c$ is a reasonable complexification of a real Banach space ${X}$, then $({X}_c)^*$ is a reasonable complexification of the real Banach space ${X}^*$. However, the natural isomorphism $\Psi$ is not always an isometry for any 2-dominating natural complexification procedure $\nu$ (i.e., $\|x+iy\|_\nu \geq \sqrt{\|x\|^2+\|y\|^2}$ for all $x,y\in\mathscr{X}$) or any 2-dominated complexification procedure $\nu$ (i.e., $\|x+iy\|_\nu \leq \sqrt{\|x\|^2+\|y\|^2}$ for all $x,y\in\mathscr{X}$), unless $\nu$ is the Lindenstrauss--Tzafriri complexification procedure (see \cite[Proposition 14]{MUN}), in which $\Psi$ is an isometry whenever ${X}$ is a real Hilbert space.\smallskip

According to a well-known property of the projective and injective tensor norms, for any real Banach space ${X}$, it follows that $({X}\otimes_\epsilon \ell_2^2)^*={X}^*\otimes_\pi \ell_2^2$ and $({X}\otimes_\pi\ell_2^2)^*={X}^*\otimes_\epsilon \ell_2^2$, where $\ell_2^2$ is identified with $\mathbb{C}_r$. In terms of complexifications, the duality existing between the injective and projective tensor norms translates into the identities
\begin{align*}
({X}_c,\|\cdot\|_T)^*&=(({X}^*)_c,\|\cdot\|_B)\hbox{ and } ({X}_c,\|\cdot\|_B)^*=(({X}^*)_c,\|\cdot\|_T)
\end{align*} for the Taylor and Bochnak norms on the complexification.\smallskip

In a milestone contribution, Bishop and Phelps \cite{BP} showed that for each real Banach space ${X}$ and each closed bounded convex subset $\mathscr{M}$ of $X$, the set $$\{\phi\in {X}^*: \phi \mbox{~attains its supremum on~} \mathscr{M} \}$$
of linear functionals supported at points of $\mathscr{M}$ is norm-dense in ${X}^*$ (see \cite{Aron, Miguel} and references therein for some generalization in several various directions). Lomonosov \cite{LOM,LOM2000} showed that this statement cannot be extended to general complex Banach spaces by constructing a closed bounded convex set with no support points.


\subsection{Extension of polynomials and multilinear mappings to the complexification}\ \smallskip

It is convenient to recall first the most basic definitions and results about polynomials on Banach spaces. The reader is referred to the excellent monograph \cite{Dineen1999} for a complete and modern exposition on polynomials on Banach spaces. A mapping $P:\mathscr{X}\rightarrow \mathscr{Y}$ between real or complex linear spaces, is an \emph{$n$-homogeneous polynomial} if there is an $n$-linear mapping $L:\mathscr{X}^n \rightarrow \mathscr{Y}$ satisfying
$P(x)={\widehat L}(x):=L(x,\ldots,x)$ for all $x\in \mathscr{X}$. According to a well-known algebraic polarization identity, for each $n$-homogeneous polynomial $P:\mathscr{X}\rightarrow \mathscr{Y}$, there exists a unique symmetric $n$-linear mapping $L:\mathscr{X}^n\rightarrow \mathscr{Y}$ (i.e., $L(x_1,\ldots, x_n) = L (x_{\sigma(1)},\ldots, x_{\sigma(n)})$
for any $(x_1,\ldots, x_n)\in \mathscr{X}^n$ and any permutation $\sigma$ of the first $n$ natural
numbers) such that $P={\widehat L}$. The unique symmetric $n$-linear mapping $L$ is called the \emph{polar} of $P$. The standard notations to represent the linear spaces of all $n$-homogeneous polynomials from $\mathscr{X}$ into $\mathscr{Y}$, the $n$-linear mappings from $\mathscr{X}$ into $\mathscr{Y}$, and the
symmetric $n$-linear mappings from $\mathscr{X}$ into $\mathscr{Y}$ are given by ${\mathcal P}_a(^{n}\mathscr{X};\mathscr{Y})$, ${\mathcal L}_a(^{n}\mathscr{X};\mathscr{Y})$, and
${\mathcal L}_a^{s}(^{n}\mathscr{X};\mathscr{Y})$, respectively. Naturally, a map
$P:\mathscr{X}\rightarrow \mathscr{Y}$ is a \emph{polynomial of degree at most $n$} if
$$
P=P_0 +P_1 +\cdots +P_n,
$$
where $P_k \in {\mathcal P}_a(^{k}\mathscr{X};\mathscr{Y})$ $(1\leq k\leq n)$ and $P_0:\mathscr{X}\rightarrow \mathscr{Y}$ is a
constant function. The polynomials of degree at most $n$ between the normed spaces $\mathscr{X}$ and $\mathscr{Y}$ are denoted by ${\mathcal P}_{n,a}(\mathscr{X};\mathscr{Y})$. If $\mathscr{Y}$ is ${\mathbb K}$ (either ${\mathbb R}$ or ${\mathbb C}$), then ${\mathcal P}_a(^{n}\mathscr{X};{\mathbb K})$, ${\mathcal L}_a(^{n}\mathscr{X};{\mathbb K})$,
${\mathcal L}_a^{s}(^{n}\mathscr{X};{\mathbb K})$, and ${\mathcal P}_{n,a}(\mathscr{X};{\mathbb K})$ are customarily replaced by ${\mathcal P}_a(^{n}\mathscr{X})$, ${\mathcal L}_a(^{n}\mathscr{X})$,
${\mathcal L}_a^{s}(^{n}\mathscr{X})$, and ${\mathcal P}_{n,a}(\mathscr{X})$, respectively.\smallskip

As it happens with linear operators, there are polynomials and multilinear mappings between Banach spaces that are not continuous. Actually, the set of noncontinuous polynomials is extraordinarily large from an algebraic viewpoint (see \cite{GMPS2012}). In any case, the continuity of polynomials and multilinear maps between infinite-dimensional Banach spaces is tightly related to the boundedness. For Banach spaces $\mathscr{X}$ and $\mathscr{Y}$, a polynomial $P\in {\mathcal P}_{n,a}(\mathscr{X};\mathscr{Y})$ or a multilinear mapping $L\in {\mathcal L}_a(^{n}\mathscr{X};\mathscr{Y})$ is continuous if and only if it is bounded on the open or closed unit ball of $\mathscr{X}$, denoted by $B_\mathscr{X}$ and $\overline{B_\mathscr{X}}$, respectively. In that case, the formulas
\begin{align*}
\|P\|&=\sup\{\|P(x)\|_\mathscr{Y}:\|x\|_\mathscr{X}\leq 1\},\\
\|L\|&=\sup\{\|L(x_1,\ldots,x_n)\|_\mathscr{Y}:\|x_k\|_\mathscr{X}\leq 1,\ k=1,\ldots,n\},
\end{align*}
define a complete norm in the spaces of continuous (bounded) $n$-homogeneous polynomials, continuous (bounded) polynomials of degree at most $n$, continuous (bounded) $n$-linear mappings, and continuous (bounded) symmetric $n$-linear mappings between the Banach spaces $\mathscr{X}$ and $\mathscr{Y}$, denoted by ${\mathcal P}(^{n}\mathscr{X};\mathscr{Y})$, ${\mathcal P}_n(\mathscr{X};\mathscr{Y})$, ${\mathcal L}(^{n}\mathscr{X};\mathscr{Y})$, and
${\mathcal L}^{s}(^{n}\mathscr{X};\mathscr{Y})$, respectively. We will rather use ${\mathcal P}(^{n}\mathscr{X})$, ${\mathcal P}_n(\mathscr{X})$, ${\mathcal L}(^{n}\mathscr{X})$, and
${\mathcal L}^{s}(^{n}\mathscr{X})$, whenever $\mathscr{Y}={\mathbb K}$.\smallskip

Throughout the rest of this section, ${X}$ and ${Y}$ will be a pair of real Banach spaces. Multilinear mappings in ${\mathcal L}_a(^{n}{X};{Y})$ admit a unique extension to a multilinear mapping in ${\mathcal L}_a(^{n}{X}_c;{Y}_c)$. Indeed,
if $L\in {\mathcal L}_a(^{n} {X}; {Y})$, then the mapping
 $$
 L_c(x_{1}^{0}+ix_{1}^{1},\ldots ,x_{n}^{0}+ix_{n}^{1})=
\sum\limits_{\epsilon_j=0,1} i^{\sum\limits_{j=1}^{n}{\epsilon}_{j}}
 L(x_{1}^{{\epsilon}_{1}},\ldots ,x_{n}^{{\epsilon}_{n}}) \quad (x_{k}^{0}, x_{k}^{1}\in {X})
 $$
is in ${\mathcal L}_a(^{n} {X}_c; {Y}_c)$ and extends $L$ (see the introduction of \cite{Kirwan}). In addition, if $L$ is bounded, then $L_c$ is bounded too for any pair of reasonable complexification norms in ${X}_c$ and ${Y}_c$. However, the norm of $L_c$ depends strongly on the complexification norms used in ${X}_c$ and ${Y}_c$.\smallskip

Similarly, if $P\in {\mathcal P}_a(^{n} {X}; {Y})$, then $P$ admits a unique complex extension to a homogeneous polynomial $P_c\in {\mathcal P}_a(^{n} {X}_c; {Y}_c)$ given by (see \cite[p. 313]{Taylor1938})
 \begin{equation*}\label{comphom}
 P_c(x+iy)=\sum\limits_{k=0}^{[\frac{n}{2}]}(-1)^{k}{\binom {n}{2k}}
 L(x^{n-2k}y^{2k})
+i\sum\limits_{k=0}^{[\frac{n-1}{2}]}(-1)^{k}{\binom {n}{2k+1}}
 L(x^{n-(2k+1)}y^{2k+1})
 \end{equation*}
for every $x, y$ in ${X}$, where $L\in {\mathcal L}_a^{s}(^{n} {X};{Y})$ is the polar of $P$ and $L(x^ly^m)$ denotes
$L(\underbrace{x,\ldots,x}_{l\ \textrm{times}},
\underbrace{y,\ldots,y}_{m\ \textrm{times}})$ for $l+m=n$.\smallskip

Another useful formula to handle the complexification of any polynomial in ${\mathcal P}_a(^n {X})$ is given by the following identity (see \cite[Theorem 4.12]{Kirwan}):
 $$
 P_c(x+iy)=\frac{2^n}{2\pi}\int_0^{2\pi}P(x\cos\theta+y\sin\theta)e^{in\theta}d\theta .
 $$
Any polynomial $P$ of degree at most $n$ in ${\mathcal P}_{n,a}({X};{Y})$ can be also extended uniquely to a polynomial $P_c \in {\mathcal P}_{n,a}({X}_c;{Y}_c)$. If $P=\sum\limits_{k=0}^{n}P_{k}$ with $P_k\in{\mathcal P}_a(^n {X}; {Y})$, then we just need to set
$P_c=\sum\limits_{k=0}^{n}P_{k,c}$, where $P_{k,c}$ is the complexification of $P_k$ for all $k=1,\ldots,n$. In the special case where ${Y}={\mathbb R},$ a modification of the argument employed in \cite[Theorem 4.12]{Kirwan} (see \cite{M1998}) can be employed to prove that
 $$
 P_{n,c}(x+iy)=\frac{2^n}{2\pi}\int_0^{2\pi}P(x\cos\theta+y\sin\theta)e^{in\theta}d\theta.
 $$

If $P$ is a bounded polynomial in ${\mathcal P}_n({X};{Y})$, then its complex extension, $P_c,$ also is a bounded polynomial in ${\mathcal P}_n({X}_c;{Y}_c)$ for any choice of reasonable complexification norms in ${X}_c$ and ${Y}_c$, although the norm of $P_c$ depends strongly on the complexification norms considered in ${X}_c$ and ${Y}_c$ (see further down).
In the special case of a finite-dimensional space
$({\mathbb R}^N, \|\cdot \|)$, the complexification of a polynomial $P$ on ${\mathbb R}^N$ is the polynomial on ${\mathbb C}^N$ that results by replacing real by complex variables in $P$, that is, the polynomial $P_c$ in $N$ complex variables is defined by
 $$
 P_c(x+iy)=P(x_1+iy_1,\ldots,x_N+iy_N)
 $$
for $x=(x_1,\ldots,x_N)$ and $y=(y_1,\ldots,y_N)$ in ${\mathbb R}^N$.\smallskip

It is simple to prove that if $P\in{\mathcal P}_n({X};{Y})$ or $L\in{\mathcal L}(^n {X};{Y})$, then $P_c$ and $L_c$ are continuous as maps between the complex Banach spaces $({X}_c\|\cdot\|_{{X}_c})$ and $({Y}_c,\|\cdot\|_{{Y}_c})$ for any choice of reasonable complexification norms $\|\cdot\|_{{X}_c}$ and $\|\cdot\|_{{Y}_c}$. It would be desirable to be able to complexify polynomials and multilinear mappings with preservation of their norms. However, that is rarely the case. If $P\in {\mathcal P}_n({X};{Y})$ or $L\in{\mathcal L}(^n {X};{Y})$, then no matter what complexification norms we consider in ${X}_c$ and ${Y}_c$, the complex extensions $P_c$ and $L_c$ of $P$ and $L$ always satisfy $\|P_c\|\ge \|P\|$ and $\|L_c\|\ge \|L\|$. The problem of estimating the size of $\|P_c\|$ has a long standing tradition. Already in 1946, Visser \cite{Visser1946} proved that if $P\in{\mathcal P}_n(\ell_\infty^m({\mathbb R}))$ with $P=P_n+\cdots+P_1+P_0$ and $P_k\in{\mathcal P}(^k\ell_\infty^m({\mathbb R}))$ for $k=1,\ldots,n$ and $P_0\in{\mathbb R}$, then
 $$
 \|P_{n,c}\|_T\leq 2^{n-1}\|P\|,
 $$
where, as usual, $P_{k,c}$ is the complex extension of $P_k$ for $k=1,\ldots,n$. Observe that, as we have commented, the Taylor norm complexifies real $\ell_\infty$-spaces in a ``natural'' way, and therefore
 $$
 \|P_{k,c}\|_T=\sup\{|P_k(z_1,\ldots,z_m)|:(z_1,\ldots,z_m)\in{\mathbb C}^m\text{ and } \|(z_1,\ldots,z_m)\|_\infty\leq 1\}.
 $$
Interestingly, the constant $2^{n-1}$ is optimal and equality is attained for the $n$th Chebysehev polynomials of the first kind $T_n$. Recall that $T_n(x)=\cos(n\arctan x)$ for $x\in [-1,1]$. In a similar fashion, it can be proved (see \cite{RACK,RACK1,REIMER}) that for $n\geq 2$, we have the optimal estimate
 $$
 \|P_{n-1,c}\|_T\leq 2^{n-2}\|P\|
 $$
with equality attained for the Chebyshev polynomial $T_{n-1}$. The following generalization to polynomials on an infinite-dimensional real Banach space can be found in \cite[Propositions 16 and 18]{MUN} (see also \cite{Kirwan}, where a slightly worse estimate is obtained).

\begin{theorem}\label{thm:homopol}{\rm\cite[Propositions 16 and 18]{MUN}} Let ${X}$ be a real Banach space, let $\nu$ be any natural complexification procedure, and let $P\in{\mathcal P}_n({X})$ with $P=P_n+P_{n-1}+\cdots+P_1+P_0$. Then the following estimations hold:
 \begin{align*}
 \|P_{n,c}\|_\nu&\leq 2^{n-1}\|P\|,\\
 \|P_{n-1,c}\|_\nu&\leq 2^{n-2}\|P\| \quad (n\ge 2).
 \end{align*}
 In particular, if $P\in{\mathcal P}(^n {X})$ and $L\in{\mathcal L}(^n {X})$, then
 \begin{align*}
 \|P_c\|_\nu&\leq 2^{n-1}\|P\|,\\
 \|L_c\|_\nu&\leq 2^{n-1}\|L\|.
 \end{align*}
 None of the constants can generally be improved.
\end{theorem}

In the previous result, equality is attained in the first two estimates for the Chebyshev polynomials $T_n$ and $T_{n-1}$, respectively. On the other hand, the $n$-homogeneous polynomial defined on $\ell_2^2$ by
 $$
 P(x,y)=\text{Re}(x+iy)^n,
 $$
for $x,y\in{\mathbb R}$, and its polar $L$ satisfy
\begin{align*}
\|P_c\|_T&=2^{n-1}\|P\|,\\
\|L_c\|_T&=2^{n-1}\|L\|.
\end{align*}

Complexification norm estimates of polynomials and multilinear mappings can be significantly improved when using specific natural complexification procedures. This is the case of the Bochnak norm (see \cite[p. 276]{Bochnak1970} and \cite{Bochnak1971}).

\begin{theorem}\label{equality}
Let ${X}$ be a real Banach space. Then, for every $L\in{\mathcal L}(^{n} {X})$, it follows that
$$ \|L_c\|_B=\|L\|.$$
\end{theorem}

Also, for the $(p)$ norms, we have the following result.

\begin{theorem}\label{complexificationp}{\rm\cite[Proposition 19]{MUN}}
 Let ${X}$ be a real Banach space and let $1\leq p\leq \infty$. Then
 for any $L\in {\mathcal L}(^n {X})$ $(n\geq 2)$, it holds that
 $$ \|L_c\|_{(p)}\leq \begin{cases}
 2^{n/2-1/2}\|L\| & \textrm{if }1 \leq p\leq 4/3,\\
 2^{n/2-2/q}\,\|L\|& \textrm{if }4/3\leq p \leq 2,\\
 2^{n/q -1}\,\|L\|& \textrm{if }2\leq p \leq \infty,
 \end{cases}
 $$
 where $q$ is the conjugate of $p$, that is, $\frac{1}{p}+\frac{1}{q}=1$ and $q=1$ if $p=\infty$.
\end{theorem}

The constant given in the previous proposition is sharp at least for
$p\geq 2$, and equality is achieved for the polar of the polynomial defined on $\ell_2^2$ by
 $$
 P(x,y)=\text{Re}(x+iy)^n,
 $$
where $x,y\in{\mathbb R}$ (see \cite{MUN}).\smallskip

The estimates on the complexification of homogeneous polynomials and multilinear forms appearing in Theorem \ref{thm:homopol} need to be increased by a factor 2 when vector-valued polynomials and multilinear operators are considered.

\begin{theorem}{\rm\cite[Proposition 25]{MUN}}
Let ${X}$ and ${Y}$ be real Banach spaces, let $P\in{\mathcal P}(^n {X};{Y})$, and let $L\in{\mathcal L}(^n {X};Y)$. Then
\begin{align*}
\|P_c\|_{T\rightarrow B}&\leq 2^{n}\|P\|,\\
\|L_c\|_{T\rightarrow B}&\leq 2^{n}\|L\|,
\end{align*}
where $\|P_c\|_{T\rightarrow B}$ {\rm(}respectively, $\|L_c\|_{T\rightarrow B}${\rm)} denotes the norm of $P_c$ {\rm(}respectively, $L_c${\rm)} as a polynomial {\rm(}respectively, multilinear operator{\rm)} between the complex Banach spaces $({X}_c,\|\cdot\|_T)$ and $({Y}_c,\|\cdot\|_B)$.
None of the inequalities can generally be improved.
\end{theorem}

We recall that a natural complexification procedure $\nu$ is 2-dominating if $\|x+iy\|_\nu\ge\sqrt{\|x\|^2+\|y\|^2}$ for all $x$ and $y$ in any real Banach space ${X}$. In the case of homogeneous polynomials and 2-dominating natural complexification procedures, we know the following result.

\begin{theorem}{\rm(\cite[Proposition 20]{MUN} and \cite[Propositions 3.10 and 3.12]{M1998})}
Let ${X}$ be a real Banach space and let $\nu$ be a 2-dominating natural complexification procedure. If $P\in {\mathcal P}(^n {X})$, then
 \begin{align*}
 \|P_c\|_{\nu}&\leq 2^{n-2}\|P\|,\quad \textrm{if $n$ is even,}\\
 \|P_c\|_{\nu}&\leq 2^{n-3/2}\|P\|,\quad \textrm{if $n$ is odd.}
 \end{align*}
If, in addition, ${X}$ is a real Hilbert space, then
 $$\|P_c\|_\nu\leq 2^\frac{n-2}{2}\|P\|.$$
\end{theorem}

The last inequality is optimal, and equality is reached for the polynomial defined on $\ell_2^2$ by $P(x,y)=\text{Re}(x+iy)^n$. Observe that the Linsdenstrauss--Tzafriri norm is 2-dominating. Also, if $P\in{\mathcal P}(^2 {X})$, then any 2-dominating natural complexification procedure satisfies $\|P_c\|_\nu=\|P\|$. This ideal situation never holds when the Taylor complexification is employed (see \cite[Proposition 22]{MUN}).\smallskip

Estimates on the norm of the complexification of nonhomogeneous polynomials have also been studied by several authors in the past. For real polynomials $P$ on the real line with degree at most $n$, Erd\"os \cite{ERDOS} proved that
$$\|P_c\|_{\mathbb D}\leq |T_n(i)|\cdot \|P\|_{[-1,1]},$$
where ${\mathbb D}=\{z\in{\mathbb C}:|z|\leq 1\}$ is the closed unit disk in the complex plane, $\|P_c\|_{\mathbb D}=\max_{z\in{\mathbb D}}|P(z)|$, and $\|P\|_{[-1,1]}=\max_{x\in[-1,1]}|P(z)|$. Obviously, the constant $|T_n (i)|$ cannot generally be improved. For general real Banach spaces, the following result is known.

\begin{theorem}{\rm\cite[Proposition 29]{MUN}}
Let $P$ be a polynomial of degree $\leq n$ on a real Banach space ${X}$ and let $\nu$ be a natural complexification procedure. Then
 $$ \|P_c\|_\nu\leq 2^{n/2}|T_n(i)|\cdot\|P\|.$$
\end{theorem}

\subsection{Zeros of polynomials in Banach spaces}\
\smallskip

To finish this section, we would like to address a topic of study (subsets and subspaces of zeros of polynomials) that has just, recently, started to develop. Thus, although it has rapidly caught the eyes of many researchers in the field, there is still plenty of ongoing work on it. This topic is closely related to that of {\it lineability and spaceability} (which, in a nutshell, consists of the study of existence of large algebraic structures within certain subsets in a topological vector space), we refer the interested reader to the works \cite{01,02,03,04} for a thorough study of the notions of lineability and spaceability.\smallskip

The study of the zeros of polynomials on complex spaces, due to its fundamental nature, has an old origin dating back at least to the 1950's (see, e.g., \cite{aronhajek2006} for references to earlier works). The case of polynomials on $\mathbb{C}^n$ has been widely investigated but the case of polynomials on infinite-dimensional Banach spaces seems to be an even richer source of challenging questions. Let us present here a classical and well-known result due to Plichko and Zagorodnyuk (1998) which is regarded as the starting point for the infinite-dimensional case.

\begin{theorem}\label{plichko-zagorodnyuk2222}\rm{\cite{plichkozagorodnyuk1998}}
	If \,$X$ is an infinite-dimensional complex Banach space and $P$ is an $n$-homogeneous polynomial on $X$, then $P^{-1}(0)$ contains an infinite-dimensional subspace $Y$.
\end{theorem}

On the other hand, if we move to the real scalar setting, the situation is totally different. This can be seen by means of the polynomial $P:\ell_{2}\rightarrow\mathbb{R}$ given by $$P(x)=
\sum\limits_{j=1}^{\infty}x_{j}^{2}.$$ In the finite-dimensional case, the field ($\mathbb{R}$ or $\mathbb{C}$) makes a big difference.

For instance, for the $2$-homogeneous polynomial $$P:\mathbb{C}^n \to \mathbb{C}, \, P(z) = z_1^2 + \cdots + z_n^2,$$ we have that
$P^{-1}(0)$ contains a vector space of dimension $[\frac{n}{2}],$ since
$$\text{span}\{e_1+ ie_2, e_3 + i e_4, e_5 + i e_6,\ldots\} \subset P^{-1}(0),$$ where $i = \sqrt{-1}$, and $e_1=(1,0,0, \ldots,0), \, e_2 = (0,1,0,0, \ldots,0)$, etc.
Nothing important can be said for $P^{-1}(0)$ if $\mathbb{K} = \mathbb{R}$. As the following theorem reveals, this example is, in fact, illustrative of the general situation in the case $\mathbb{K} = \mathbb{C}$ (see, e.g., \cite{AronRueda, arongonzalozagorodnyuk2000, aronboydryanzalduendo2003, ferrerprims2007, ferrer2009}, and the references therein).

\begin{theorem}\label{plichko-zagorodnyuk}\rm{\cite{plichkozagorodnyuk1998, AronRueda, Zagorod2001, arongonzalozagorodnyuk2000}}
	Let $X$ be a complex Banach space. Given positive integers $n$ and $k,$ there is an integer \,$m(n,k) \in \mathbb{N}$
	\,such that, whenever \,$\dim (X) = k$ \,and \,$P:X \to \mathbb{C}$ \,is an $n$-homogeneous polynomial, the set \,$P^{-1}(0)$
	\,contains a subspace of dimension at least $m(n,k).$ Moreover, $m(n,k) \to \infty$ as $k \to \infty.$
\end{theorem}

\begin{corollary}[\cite{plichkozagorodnyuk1998, AronRueda, Zagorod2001, arongonzalozagorodnyuk2000}]
	Let $P:\mathbb{C}^k \to \mathbb{C}$ be an arbitrary {\rm (}not necessarily homogeneous{\rm )} polynomial
	of degree $n.$ Then there is a subspace $V \subset \mathbb{C}^k,$
	whose dimension depends only on $k,$ such that $\dim (V) \to \infty$ as $k \to \infty$,
	satisfying the condition \,$P|_V \equiv P(0).$
\end{corollary}

Let us point out that, when one considers polynomials of the form $\sum x_{j}^{2}$, the case of real polynomials needs a special approach, where odd-homogeneous polynomials and even-homogeneous polynomials are investigated by different fronts (see, e.g., \cite{aronhajek2006}).\smallskip

The following two results also show how different the answer may end up being when comparing the real and complex frameworks.

\begin{theorem}\label{aron-hajek-positivity}\rm{\cite{aronhajek2007}}
	Given any real, separable, infinite-dimensional Banach space $X$ and any odd $n \in \mathbb{N},$ there is an $n$-homogeneous polynomial $P:X\to \mathbb{R}$ such that $P^{-1}(0)$ does not contain an infinite-dimensional subspace.
\end{theorem}

\begin{theorem}\rm{\cite{fernandez2006}}
	Let $E$ be a complex Banach space containing $\ell_\infty.$ For every $n,$ every
	$n$-homogeneous $P:E \to \mathbb{C}$ vanishes on a nonseparable subspace of $E$.
\end{theorem}

Moreover, in \cite{fernandez2006} it is also shown that in the case of real $\ell_\infty$, if $P:\ell_\infty \to \mathbb{R}$ vanishes on a copy of $c_0$, then $P \equiv 0$ on a nonseparable subspace. Furthermore, in 2009, Avil\'es and Todorcevic \cite{avilestodorcevic2009} showed that there exists a $2$-homogeneous polynomial $P:\ell_1(\aleph_1) \to \mathbb{C}$ such that $P^{-1}(0)$ contains no nonseparable subspace.

\vskip .15cm

\begin{theorem}\rm{(Avil\'{e}s, Todorcevic, 2009, \cite{avilestodorcevic2009})}
	There exists a $2$-homogeneous polynomial $P:\ell_1(\mathfrak{c}) \to \mathbb{C}$ such that $P^{-1}(0)$ contains both separable and nonseparable maximal subspaces.
\end{theorem}

Avil\'es and Todorcevic \cite{avilestodorcevic2009} also provide new viewpoints on the research of the zero set of complex polynomials, including new techniques and connections with results related to the existence of certain partitions. \smallskip

Regarding the case of $2$-homogeneous polynomial on a real Banach space $X$, let us recall that a $2$-homogeneous polynomial $P:X \to \mathbb{R}$ is said to be {\em positive
	definite} if $P(x) \geq 0$ for every $x \in X$ and $P(x) = 0$ only for $x = 0.$ The following very recent results by Ferrer are of major importance in this direction (for the case of a compact topological Hausdorff space $K$).

\begin{theorem}\rm{\cite{ferrer2007}} The space $C(K)$ satisfies the following dichotomy. Either
	\begin{enumerate}
		\item[\rm (i)] It admits a positive definite continuous $2$-homogeneous real-valued
		polynomial, or
		\item[\rm (ii)] Every continuous $2$-homogeneous real-valued polynomial vanishes in a
		non-separable closed linear subspace.
	\end{enumerate}
\end{theorem}

When $X=c_{0}\left(\Gamma\right)$ the following result holds for general polynomials (non necessarily homogeneous):

\begin{theorem}\rm{\cite{ferrer2007}} Let \,$\Gamma$ be an uncountable set. If $P:c_{0}\left( \Gamma\right) \rightarrow\mathbb{R}$ is a continuous polynomial, then there is a closed linear subspace $E$ of $c_{0}\left( \Gamma\right) $ such that $E\subset P^{-1}\left( 0\right) $ and $E$ is isometric to $c_{0}\left( \Gamma\right).$
\end{theorem}

Also, Ferrer et al. (2019, \cite{ferrer2019}) proved that whenever $X$ is a real Banach space which cannot be linearly and continuously injected into a Hilbert space, then for any $2$-homogeneous continuous polynomial $P$ on $X$, the zero set $P^{-1}(0)$ is not separable.

\section{Real and complex (classical) polynomial inequalities}

\subsection{Real and complex polarization constants}\ \smallskip

For real or complex Banach spaces $\mathscr{X}$ and $\mathscr{Y}$, it has been already mentioned that any polynomial in ${\mathcal P}_a(^n\mathscr{X};\mathscr{Y})$ is induced by a unique symmetric $n$-linear mapping in ${\mathcal L}_a^s(^n\mathscr{X};\mathscr{Y})$, which we call the polar of $P$. Along this note, the symbol $\widehat{L}$ stands for the polynomial induced by $L$. According to this notation, the mapping
$${\mathcal L}_a^s(^n\mathscr{X};\mathscr{Y})\ni L\mapsto \widehat{L}\in {\mathcal P}_a(^n\mathscr{X};\mathscr{Y})$$
is a natural linear isomorphism, whose inverse is given by the
so-called polarization formula (see \cite[Corollaries 1.6 and 1.7]{Dineen1999}). The following is just one of the many forms in which the polarization formula can be found in the literature:
$$ L(x_1,\ldots,x_n)=\frac{1}{2^nn!}\sum\limits_{\epsilon=\pm 1}\epsilon_1\cdots\epsilon_nP(\epsilon_1x_1+\cdots+\epsilon_nx_n).$$
By restricting our attention to continuous polynomials and continuous symmetric multilinear mappings, the following estimates establish a relationship between the norms of a polynomial $P\in {\mathcal P}(^n\mathscr{X};\mathscr{Y})$ and the norm of its polar $L\in{\mathcal L}^s(^n\mathscr{X};\mathscr{Y})$:
 \begin{equation}\label{eq:Martin}
 \|P\|\leq \|L\|\leq \frac{n^n}{n!}\|P\|.
 \end{equation}
While the first inequality is trivial since $P$ is a restriction of $L$, the second can be derived from the polarization formula (see \cite[Proposition 1.8]{Dineen1999} for a modern proof). The previous estimates show that the natural algebraic isomorphism $L\mapsto \widehat{L}$ is also a topological isomorphism between the spaces ${\mathcal L}_a^s(^n\mathscr{X};\mathscr{Y})$ and ${\mathcal P}(^n\mathscr{X};\mathscr{Y})$ with norm 1, and whose inverse has norm at most $\frac{n^n}{n!}$. It is important to mention that the constant $\frac{n^n}{n!}$ cannot generally be improved since the polynomial $\Phi_n(x_1,\ldots,x_n)=x_1\cdots x_n$ defined on $\ell_1^n({\mathbb K})$ and its polar $L_n$ satisfy $\|L_n\|=\frac{n^n}{n!}\|\Phi_n\|$. All polynomials satisfying the latter identity are called \emph{extremal}.\smallskip

Although $\frac{n^n}{n!}$ is optimal, in general, it might be improved for specific spaces. This serves as a motivation for the definition of the $n$th polarization constant ${\mathbb K}(n,\mathscr{X})$ of a Banach space $\mathscr{X}$ over ${\mathbb K}$:
 $$
 {\mathbb K}(n,\mathscr{X})=\inf\{C>0:\|L\|\leq C\|\widehat{L}\|\text{ for all $L\in{\mathcal L}^s(^n\mathscr{X})$}\}.
 $$
Also, the polarization constant of $\mathscr{X}$ is defined as
 $$
 {\mathbb K}(\mathscr{X})=\limsup_n\sqrt[n]{{\mathbb K}(n,\mathscr{X})}.
 $$
The calculation of ${\mathbb K}(n,\mathscr{X})$ and ${\mathbb K}(\mathscr{X})$ has been studied in the past in several occasions. Depending on whether ${\mathbb K}$ is ${\mathbb R}$ or ${\mathbb C}$, different techniques are used, and sometimes different results are obtained. We present below some remarkable results on this topic, stressing the difference between the real and complex case.\smallskip

It has been pointed out above that ${\mathbb K}(n;\ell_1^n)=\frac{n^n}{n!}$. The constant $\frac{n^n}{n!}$, however, is attained in different ways, depending on whether ${\mathbb K}$ is ${\mathbb R}$ or ${\mathbb C}$.

\begin{theorem}{\rm\cite[Corollary 2]{Sarant1987}}
 An $n$-dimensional complex Banach space is isometrically isomorphic to $\ell_1^m({\mathbb C})$ if and only if ${\mathbb C}(n,E)=\frac{n^n}{n!}$. Also, if ${\mathbb C}(n,E)=\frac{n^n}{n!}$ and $L\in{\mathcal L}^s(^nE)$ is extremal, that is, $\|L\|=\frac{n^n}{n!}\|\widehat{L}\|$, then $\widehat{L}(z_1,\ldots,z_n)=cz_1,\dots,z_n$ for some $c\in{\mathbb C}$.
\end{theorem}
The previous results states that there is, essentially, a unique extremal polynomial in any $n$-dimensional complex Banach space $E$ with ${\mathbb C}(n,E)=\frac{n^n}{n!}$. The same remains true for real Banach spaces with dimension two or three. However, it is no longer true when $n\geq 4$.

\begin{theorem}{\rm\cite[Corollaries 10 and 12]{Sarant1999}}
 If $n\ge 4$ and $|\gamma_{ij}|\le \frac{1}{3\cdot 4^4}$ for $1\leq i<j\leq n$, then the polynomials
 $$
 P(x_1,\ldots,x_n)=
 \begin{cases}
 c\left(x_1x_2x_3x_4+\sum\limits_{1\leq i<j\leq 4} \gamma_{ij}(x_i^2-x_j^2)^2 \right)&\text{if }n=4,\\
 c\left(x_1x_2x_3x_4+\sum\limits_{1\leq i<j\leq n} \gamma_{ij}(x_i^2-x_j^2)^2 \right)x_5\cdots x_n&\text{if }n>4,
 \end{cases}
 $$
 with $c\in{\mathbb R}$, are extremal in ${\mathcal P}(^n\ell_\infty^n)$.
\end{theorem}

Another remarkable difference between the real and complex cases in connection with polarization constants occurs in Hilbert spaces. It is well known that ${\mathbb K}(n;\mathscr{H})=1$ for any real or complex Hilbert space $\mathscr{H}$ and every $n\in{\mathbb N}$. Hence ${\mathcal L}^s(^n \mathscr{H})$ and ${\mathcal P}(^n\mathscr{H})$ are isometrically isomorphic no matter whether $\mathscr{H}$ is a real or complex Hilbert space. The fact ${\mathbb K}(n; \mathscr{H})=1$ was proved by Kellogg \cite{Kellogg} and Van der Corput and Schaake \cite{Vandercorput} when $\mathscr{H}$ is finite-dimensional. Banach \cite{Banach} gave a proof in the case when $H=\ell_2$. For a comprehensive exposition on the topic the reader is referred to \cite{Dineen1999,Harris1975}.\smallskip

The divergence between the conclusions in the real and complex cases can be found in the following result.

\begin{theorem}{\rm\cite[Proposition 2.8]{Sarant1993}} If $X$ is a real Banach space such that ${\mathbb R}(n,X)=1$ for every $n\in{\mathbb N}$, then $X$ is a Hilbert space. Actually, ${\mathbb R}(2,X)=1$ is enough to conclude that $X$ is a Hilbert space.
\end{theorem}

The previous result is not true in the complex setting. If $\mathscr{H}$ is a complex Hilbert space and $\mathscr{H}\oplus_\infty{\mathbb C}$ is the space $\mathscr{H}\times {\mathbb C}$ endowed with the norm
 $$
 \|(x,\lambda)\|_\infty=\max\{\|x\|,|\lambda\},
 $$
then $\mathscr{H}\times {\mathbb C}$ is not a complex Hilbert space, and ${\mathbb C}(n,\mathscr{H}\oplus_\infty{\mathbb C})=1$ (see \cite[p. 94]{Sarant1987BGMS}).\smallskip

The value of the polarization constants of $\ell_\infty$-type spaces is also another issue, where the real and complex cases diverge. It is well known (see, for instance, \cite{Dineen1999} or \cite{Harris1975}) that
 $$
 {\mathbb C}(n,\ell_\infty^m)\leq \frac{n^\frac{n}{2}(n+1)^\frac{n+1}{2}}{2^nn!}.
 $$
However, the same estimate does not hold for ${\mathbb R}(n,\ell_\infty^m)$. Indeed, if $P\in{\mathcal P}(^4\ell_\infty^4({\mathbb R}))$ is defined by
 $$
 P(x_1,x_2,x_3,x_4)=(x_1^2-x_2^2)^2-(x_3^2-x_4^2)^2
 $$
and $L$ is its polar, then
 $$
 \|L\|\ge 3\|P\|,
 $$ which implies that ${\mathbb R}(4,\ell_\infty^4)\geq 3 > \frac{ 25\sqrt{5}}{24}\ge {\mathbb C}(4,\ell_\infty^4)$ (see \cite{Sarant1986}).\smallskip

The polarization constant of finite-dimensional spaces behaves differently in real and complex Banach spaces, as shown recently. If $\mathscr{X}$ is a finite-dimensional complex Banach space, then ${\mathbb C}(\mathscr{X})=1$ (see \cite[Theorem 1.1]{Dimant2021}). However, it was proved in \cite{Dimant2021} that ${\mathbb R}(\ell_1^d)>1$.\smallskip

The last topic we shall deal with in connection to the polarization constants emerges from the following result by Harris \cite[Theorem 1]{Harris1975}:

\begin{theorem}\label{thm:Harris1975}
	Let ${\mathscr X}$ be a complex Banach space, $P\in\mathcal{P}\left(^n {\mathscr X}\right) $ with polar $L\in\mathcal{L}^s\left(^n {\mathscr X}\right)$ and $k_1,\ldots,k_m\in\mathbb{N}\cup\left\{0\right\} $ such that $ k_1+\cdots+k_m=n$. Assume that $ x_1,\ldots,x_m$ are unit vectors in $\mathscr{X}$ satisfying
	\begin{align*}
	\left\|z_1 x_1+\cdots+ z_m x_m\right\|\le\left\|\left(z_1,\ldots,z_m\right)\right\|_p
	\end{align*}
	for all $ \left(z_1,\ldots,z_m\right)\in\mathbb{C}^m $ and for a given $ 1\le p\le\infty $ (here, $\|\cdot\|_p$ denotes the usual $p$-norm). Then,
	\begin{align*}
	\left|L\left(x_1^{k_1},\ldots,x_m^{k_m}\right)\right|\le\frac{k_1 !\cdots k_m ! n^\frac{n}{p}}{k_1^\frac{k_1}{p}\cdots k_m^\frac{k_m}{p}n!}\left\|P\right\| . \label{42}
	\end{align*}
	Moreover, if $ {\mathscr X}=\ell_p^m $, then there exist $ P\in\mathcal{P}\left(^n {\mathscr X}\right) $, $ L\in\mathcal{L}^s\left(^n {\mathscr X}\right) $ with $ L\not\equiv 0$, and unit vectors $x_1,\ldots,x_m\in \mathscr{X}$ for which equality is attained.
\end{theorem}
Under the assumptions of the previous result, considering that
	$$
	\|z_1 x_1+\cdots+ z_m x_m\|\le|z_1|+\cdots+|z_m|\leq\|\left(z_1,\ldots,z_m\right)\|_1
	$$
	for any $m$-tuple of unit vectors $x_1,\ldots,x_m$, it follows that
	$$
	\left|L\left(x_1^{k_1},\ldots,x_m^{k_m}\right)\right|\le\frac{k_1 !\cdots k_m ! n^n}{k_1^{k_1}\cdots k_m^{k_m} n!}\left\|P\right\| .
	$$
	The previous estimate is sharp, however $\frac{k_1 !\cdots k_m ! n^n}{k_1^{k_1}\cdots k_m^{k_m} n!}$ might be replaced by a better (smaller) constant for a specific choice of ${\mathscr X}$. This motivates the definition of the generalized polarization constants:
	\begin{definition}
	If ${\mathscr X}$ is a Banach space over ${\mathbb K}$ and $k_1,\ldots,k_m\in{\mathbb N}\cup\{0\}$, then ${\mathbb K}(k_1,\ldots,k_m,{\mathscr X})$ denotes
	$$
	\inf\{M>0:|L(x_1^{k_1},\ldots,x_m^{k_m})|\leq M\|{\widehat L}\|,\ L\in{\mathcal L}^s(^n{\mathscr X}),\ x_1,\ldots,x_m\in S_{\mathscr X}\}.
	$$
	\end{definition}
The comments made above show clearly that for all complex Banach spaces ${\mathscr X}$, we have
	$$
	1\leq {\mathbb C}(k_1,\ldots,k_m,{\mathscr X})\leq \frac{k_1 !\cdots k_m ! n^n}{k_1^{k_1}\cdots k_m^{k_m} n!},
	$$
	where the second inequality is sharp.\smallskip

The study of ${\mathbb R}(k_1,\ldots,k_m,{\mathscr X})$ for any real Banach space ${\mathscr X}$ is subjected to two differential facts with respect to the complex case:
\begin{enumerate}[$(1)$]
	\item First, it is not true in general that
	$$
	{\mathbb R}(k_1,\ldots,k_m,{\mathscr X})\leq \frac{k_1 !\cdots k_m ! n^n}{k_1^{k_1}\cdots k_m^{k_m} n!}
	$$
	for all real Banach space ${\mathscr X}$.
	\item Second, the best upper bound on ${\mathbb R}(k_1,\ldots,k_m,{\mathscr X})$ for arbitrary real Banach spaces ${\mathscr X}$ is not known.
\end{enumerate}
As for the first issue, it was proved in \cite{Sarant1986} that ${\mathbb R}(2,2,\ell^4_\infty({\mathbb R}))=3$ whereas ${\mathbb C}(2,2,{\mathscr X})\leq \frac{8}{3}<3$ for any complex Banach space ${\mathscr X}$. The second fact is still under study nowadays. Several estimates on ${\mathbb R}(k_1,\ldots,k_m,{\mathscr X})$ are known, but the best fit for ${\mathbb R}(k_1,\ldots,k_m,{\mathscr X})$ is still an open question to our knowledge. Harris (see \cite[Corollary 7]{Harris1997}) proved that for any real Banach space $X$, we have
	$$
	|L(x_1^{k_1},\ldots,x_m^{k_m})|\leq \sqrt{\frac{n^n}{k_1^{k_1}\cdots k_m^{k_m}}}\|P\|
	$$
	for every $k_1,\ldots,k_m\in{\mathbb N}\cup\{0\}$ with $k_1+\cdots+k_m=n$, unit vectors $x_1,\ldots,x_m\in X$ and $P\in{\mathcal P}(^n{X})$ with polar $L\in{\mathcal L}^s(^n{X})$. However, the constant $\sqrt{\frac{n^n}{k_1^{k_1}\cdots k_m^{k_m}}}$ seems to be far from being optimal since, letting $m=n$ and $k_1=\ldots=k_n=1$, we arrive at
	$$
	\|L\|\leq n^\frac{n}{2}\|P\|,
	$$
	which can be substancially improved according to \eqref{eq:Martin}. Another estimate on $\mathbb{R}\left(k_1,\ldots,k_m, X\right)$ can be found in Harris' commentaries to problems 73 and 74 of the Scottish Book (see \cite{ScottishBook}), where it is shown that
	\begin{align*}
	\mathbb{R}\left(k_1,\ldots,k_m, X\right)\le \frac{n^n r^l}{n!}
	\end{align*}
	with $ \displaystyle r=\frac{1+e^{-2}}{2} $ and $ \displaystyle l=\sum_{i=1}^m \left\lceil \frac{k_i}{2}\right\rceil$.\smallskip

Let us mention, to finish this section, that Papadiamantis and Sarantopoulos \cite{PS2016} established a number of analogs of Theorem \ref{thm:Harris1975} in a real setting. For instance, if $X$ is the a real $L_p(\mu)$ with $p\geq 1$, $x_1,\ldots,x_m$ are unit vectors in $X$ with disjoint supports, $k_1,\ldots,k_m\in{\mathbb N}\cup\{0\}$ with $k_1+\cdots+k_m=n$, and $P\in{\mathcal P}(^nX)$ with polar $L\in{\mathcal L}^s(^nX)$, then
	$$
	|L(x_1^{k_1},\ldots,x_m^{k_m})|\leq c(k_1,\ldots,k_m,p)\|P\|,
	$$
	where
	$$
	c(k_1,\ldots,k_m,p)=
	\begin{cases}
	\frac{(k_1^{p-1}+\cdots+k_m^{p-1})^\frac{n}{p}}{n!}&\text{if $p\ge n$,}\\
	\frac{m^\frac{n-p}{p}(k_1^{n-1}+\cdots+k_m^{n-1})}{n!}&\text{if $p< n$.}
	\end{cases}
	$$

\subsection{Real and complex polynomial Bohnenblust--Hille inequality}\ \smallskip

If $\mathscr{X}$ is a Banach space, then the problem of computing the value of the norms
\begin{align*}
 \|P\|: & =\sup\{|P(x)|:x\in{\mathcal{B}}_{\mathscr{X}}\},\\
 \|L\|: & =\sup\{|L(x_{1},\ldots,x_{m})|:x_{1},\ldots x_{m}\in{\mathcal{B}%
 }_{\mathscr{X}}\}
\end{align*}
is usually intractable. For this reason, it would be interesting to obtain reasonably good
estimates on it. In the case that $\mathscr{X}$ is finite-dimensional, the $\ell_{p}$ norm of the coefficients of a given polynomial
on ${\mathbb{K}}^{n}$ (${\mathbb{K}}={\mathbb{R}}$ or ${\mathbb{C}}$) is much easier to handle. Recall that an
$m$-homogeneous polynomial in ${\mathbb{K}}^{n}$ can be written as
\[
P(x)={\sum\limits\limits_{\left\vert \alpha\right\vert =m}}a_{\alpha}x^{\alpha},
\]
where $x=(x_{1},\ldots,x_{n})\in{\mathbb{K}}^{n}$, $\alpha=(\alpha_{1}%
,\ldots,\alpha_{n})\in({\mathbb{N}}\cup\{0\})^{n}$, $|\alpha|=\alpha
_{1}+\cdots+\alpha_{n}$, and $x^{\alpha}=x_{1}^{\alpha_{1}}\cdots x_{n}%
^{\alpha_{n}}$.\smallskip

Thus we define the $\ell_{p}$ norm of $P$, with $p\geq1$, as
$$|P|_{p}=\left( \sum\limits_{|\alpha|=m}|a_{\alpha}|^{p}\right)^\frac{1}{p}.$$
If $E$ has finite dimension $n$, then the polynomial norm $\|\cdot\|$ and the
$\ell_{p}$ norm $|\cdot|_{p}$ ($p\geq1$) are equivalent, and therefore, there
exist constants $k(m,n), K(m,n)>0$ such that
\begin{equation}
 \label{equ:equiv}k(m,n)|P|_{p}\leq\|P\|\leq K(m,n)|P|_{p}
\end{equation}
for all $P\in{\mathcal{P}}(^{m}E)$. The latter inequalities may provide a good
estimate on $\|P\|$ as long as we know the exact value of the best possible
constants $k(m,n)$ and $K(m,n)$ appearing in \eqref{equ:equiv}.\smallskip

The problem presented above is an extension of the well-known polynomial \linebreak
Bohnenblust--Hille inequality. It was
proved in \cite{bh} that there exists a constant $D_{m}\geq1$ such that for
every $P\in{\mathcal{P}}(^{m}\ell_{\infty}^{n})$, we have
\begin{equation}
 |P|_{\frac{2m}{m+1}}\leq D_{m}\Vert P\Vert. \label{equ:BH}%
\end{equation}
Observe that \eqref{equ:BH} coincides with the first inequality in
\eqref{equ:equiv} for $p=\frac{2m}{m+1}$ except for the fact that $D_{m}$ in
\eqref{equ:BH} can be chosen in such a way that it is independent from the
dimension $n$.\smallskip

As a matter of fact, Bohnenblust and Hille \cite{bh} showed that $\frac{2m}{m+1}$ is
optimal in \eqref{equ:BH} in the sense that for $p<\frac{2m}{m+1}$, any
constant $D$ fitting in the inequality
\[
|P|_{p}\leq D\Vert P\Vert,
\]
for all $P\in{\mathcal{P}}(^{m}\ell_{\infty}^{n})$, depends necessarily on $n$.\smallskip

The polynomial and multilinear Bohnenblust--Hille inequalities were only \textit{rediscovered} in the last
few years. These inequalities (or, more precisely, the constants appearing in them) have shown to have quite an impact in several fields of mathematics, such as operator theory, Fourier and harmonic analysis, complex analysis, analytic number theory, and quantum
information theory (see, for example, \cite{bps,RACSAM2021, Boas, boas, ManDom,annals2011, DD, diniz2, monta, psseo} and references therein).\smallskip

The best constants in \eqref{equ:BH} may depend on whether we consider the
real or the complex version of $\ell_{\infty}^{n}$, which motivates the
following definition:
\[
D_{{\mathbb{K}},m}:=\inf\left\{ D>0:\text{$|P|_{\frac{2m}{m+1}}\leq D\Vert
 P\Vert$, for all $n\in{\mathbb{N}}$ and $P\in{\mathcal{P}}(^{m}\ell_{\infty
 }^{n})$}\right\}.
\]
If we restrict our attention to ${\mathcal{P}}(^{m}\ell_{\infty}^{n})$ for some
$n\in{\mathbb{N}}$, then we define
\[
D_{{\mathbb{K}},m}(n):=\inf\left\{ D>0:\text{$|P|_{\frac{2m}{m+1}}\leq D\Vert
 P\Vert$ for all $P\in$}{\mathcal{P}}(^{m}\ell_{\infty}^{n})\right\}.
\]
Note that $D_{{\mathbb{K}},m}(n)\leq D_{{\mathbb{K}},m}$ for all
$n\in{\mathbb{N}}$.\smallskip

It was recently shown in \cite{bps} that the complex polynomial
Bohnenblust--Hille inequality is, at most, subexponential, that is, for any
$\varepsilon>0$, there is a constant $C_{\varepsilon}>0$ such that
$$D_{{\mathbb{C}},m}\leq C_{\varepsilon}\left( 1+\varepsilon\right)^{m}$$ for
all positive integers $m$. However, in the real case, the behavior is quite different, more precisely, $$\lim\sup_{m}D_{\mathbb{R},m}^{1/m}=2.$$



\subsection{Bernstein and Markov type inequalities in Banach spaces}\ \smallskip

Estimates on the derivatives of polynomials are known as Bernstein and Markov inequalities. The classical estimates on the norm of the first and successive derivatives of a polynomial in one real variable proved by the brothers Markov in the late 19th century have been generalized in several forms to the case of polynomials in an arbitrary Banach space. V.A. Markov \cite{MAR} proved in 1892 that
 $$
 \|P^{(k)}\|\leq \frac{n^2(n^2-1^2)\cdots (n^2-(k-1)^2)}{1\cdot 3\cdots (2k-1)}\|P\|
 $$
for every polynomial in ${\mathcal P}_n({\mathbb R})$. The norms are calculated as the supremum of the absolute value over the unit interval $[-1,1]$. Equality is attained for the $n$th Chebyshev polynomial of the first kind, namely, $T_n(x)=\cos(n\arccos x)$ for $x\in[-1,1]$. The results had been previously proved by Markov in 1889 for the first derivative, motivated by a question of Mendeleiev, author of the periodic table, who was interested in estimating the maximum value of the derivative of a quadratic polynomial. Markov's estimate on the $k$th derivative was generalized in 2002 to polynomials on a real Hilbert space \cite{Munoz2002} and for polynomials on an arbitrary real Banach space in 2010 \cite{Harris2010} (see also \cite{Skalyga1,Skalyga2}). In fact, if ${X}$ is a real Banach space, then
 $$
 \|\widehat{D}^{(k)}P(x)\|\leq \frac{n^2(n^2-1^2)\cdots (n^2-(k-1)^2)}{1\cdot 3\cdots (2k-1)}\|P\|
 $$
for every $P\in{\mathcal P}_n(X)$ and every $x\in X$ with $\|x\|\leq 1$. In the last inequality $D^{(k)}P$ stands for the $k$-th Fr\'echet derivative of $P$ and accordingly, $\widehat{D}^{(k)}P$ is the $k$-homogeneneous polynomial induced by $D^{(k)}P$.\smallskip

The situation is completely different in the complex setting. First, the well-known \emph{Bernstein's inequality of trigonometric polynomials} states that
 $$
 |T'(\theta)|\leq n\|T\|
 $$
for all $\theta\in{\mathbb R}$ and every trigonometric polynomial $T$ of degree $n$. As a consequence of Bernstein's inequality, complex polynomials in ${\mathcal P}_n({\mathbb C})$ satisfy
 $$
 \|P^{(k)}\|\leq \frac{n!}{(n-k)!}\|P\|,
 $$
where now the norms are calculated as the supremum of the modulus over the unit disk. Equality is attained for $P(z)=z^n$. This divergence between the real and complex cases in one variable is translated to the infinite-dimensional case.\smallskip

It is interesting to observe that in any real Hilbert space ${H}$, homogeneous polynomials satisfy the following estimate:
 $$
 \|DP(x)\|\leq n\|P\|
 $$
for all $x\in {H}$ with $\|x\|\leq 1$ and every $P\in{\mathcal P}(^n {H})$. Hence homogeneous polynomials on a real Hilbert space satisfy Bernstein's inequality. As a matter of fact, this is a characteristic property of real Hilbert spaces, that is, a real Banach space ${X}$ is an inner product space if and only if $\|DP(x)\|\leq n\|P\|$ for all $x\in {X}$ with $\|x\|\leq 1$ and every $P\in{\mathcal P}(^n {X})$ (see, for instance, \cite{Dineen1999}). Here we find another worth mentioning difference between the real and complex settings since Bernstein's inequality does not characterize complex inner product spaces. Indeed, in \cite{Harris1975}, it was shown that $\ell_\infty^2({\mathbb C})$ satisfies Bernstein's inequality, although it is not a Hilbert space.\smallskip

There is yet one more significant difference between real and complex Bernstein-Markov type inequalities in Banach spaces. In the rest of the section, we restrict our attention to homogeneous polynomials on Banach spaces. Sarantopoulos found in 1991 a good Markov estimate for the polynomial associated to the $k$th Fr\'echet derivative on a homogeneous polynomial on a complex Banach space (see \cite{Sarantopoulos1991}). Sarantopoulos results in this line were improved by Harris in 1997. If ${\mathscr X}$ is a complex Banach space, then (see \cite[Corollary 3]{Harris1997})
$$
\|{\widehat D}^{k}P(x)\|\leq
\begin{cases}
\frac{n^n k!}{k^k(n-k)^{n-k}}\|P\|&\text{if $\|x\|\leq 1$,}\\
\frac{n^n k!}{k^k(n-k)^{n-k}}\|P\|\|x\|^{n-k}&\text{if $\|x\|\geq 1$,}
\end{cases}
$$
for all $P\in{\mathcal P}(^n{\mathscr X})$. In particular,
$$
\|{\widehat D}^{k}P\|\leq \frac{n^n k!}{k^k(n-k)^{n-k}}\|P\|.
$$
The latter estimate had already been established by Harris in 1975 \cite[Corollary 1]{Harris1975}. The constant $\frac{n^n k!}{k^k(n-k)^{n-k}}$ cannot generally be improved since equality is attained in ${\mathscr X}=\ell_1^2({\mathbb C})$. \smallskip

Markov's inequalities for homogeneous polynomials on real Banach spaces provide different estimaes. In general, these type of problems in a real setting are more difficult to tackle. In addition, the results that are known for real Banach spaces are not so explicit and clear as Harris' estimates. Let us see what we know for real homogeneous polynomials. Harris \cite{Harris1997} proved that there exist constants $c_{n,k}>0$ such that
	$$
	\|{\widehat D}^{k}P(x)\|\leq c_{n,k}\|P\|
	$$
where $X$ is any real Banach space, $x\in B_X$, $P\in{\mathcal P}(^nX)$ and the optimal choice for $c_{n,k}$ can be obtained as a solution
to an extremal problem for polynomials of one real variable. The following bounds on $c_{n,k}$ follow from \cite{Sarantopoulos1991}:
	$$
	\frac{n^nk!}{k^k(n-k)^{n-k}}\leq c_{n,k}\leq {n\choose k}\frac{n^\frac{n}{2}k!}{k^\frac{k}{2}(n-k)^{\frac{n-k}{2}}},
	$$
for $1\leq k\leq n$. The upper bound can be improved for large values of $n$ as follows: there exists a constant $M>0$ such that
	$$
	c_{n,k}\leq (Mn\log n)^k.
	$$
The latter bound is due to Nevai and Totik \cite{NT1986} and seems to provide the exact asymptotic growth of the $c_{n,k}$'s at least for $k=1$ and $k=2$ (see \cite{RS2003}). Unfortunately no closed formula is known to provide the exact value of the $c_{n,k}$'s. However, the method developed by Harris \cite{Harris1997} can be applied to approach not only the exact value of $c_{n,k}$ for a given choice of $n$ and $k$, but also a procedure to construct a homogeneous polynomial for which $\|{\widehat D}^{k}P(x)\|= c_{n,k}\|P\|$ (see \cite[Table I]{Harris1997}).



\subsection{Linear polarization constants}\ \smallskip

In the literature, linear polarization constants represent the ratio between the norm of the product of linear forms and the product of the norms of the linear forms. More specifically, if ${\mathscr X}$ is a Banach space, real or complex, and $L_1,\ldots,L_n$ are $n$ bounded linear functionals in ${\mathscr X}^{*}$, then the $n$-homogeneous polynomial $P\in{\mathcal P}(^n{\mathscr X})$ defined by
	$$
	P(x)=L_1(x)\cdots L_n(x)
	$$
obviously satisfies
	$$
	\|P\|\leq \|L_1\|\cdots \|L_n\|.
	$$
On the other hand, it can be proved (\cite{BST1998}) that there exists a universal constant $K_n$, depending only on $n$, such that
	$$
	\|L_1\|\cdots \|L_n\|\leq K_n\|P\|.
	$$
The authors of \cite{BST1998} showed that whenever ${\mathscr X}$ is a complex Banach space then
	\begin{equation}\label{eq:lin_pol_const}
	\|L_1\|\cdots\|L_n\|\leq n^n\|L_1\cdots L_n\|
	\end{equation}
for every choice of bounded linear functionals $L_1,\ldots,L_n\in{\mathscr X}^{*}$. Moreover, if ${\mathscr X}=\ell_1^n({\mathbb C})$ and $L_k(z_1,\ldots,z_n)=z_k$, then
	$$
	\|L_1\|\cdots\|L_n\|= n^n\|L_1\cdots L_n\|,
	$$
proving that, at least in a complex setting, $K_n=n^n$ is the smallest possible constant in the inequality
	$$
	\|L_1\|\cdots\|L_n\|= K_n\|L_1\cdots L_n\|
	$$
for all complex Banach spaces ${\mathscr X}$ and all $L_1,\ldots,L_n\in{\mathscr X}^{*}$.\smallskip

For real Banach spaces, using a complexification argument, it can also be proved that there exists a universal constant $K_n$ depending only on $n$ such that 	
	\begin{equation}\label{eq:lin_pol_const_real}
	\|L_1\|\cdots\|L_n\|\leq K_n\|L_1\cdots L_n\|
	\end{equation}
for every real Banach space $X$ and every $L_1,\ldots,L_n\in X^{*}$. However, the best (smallest) possible choice for $K_n$ in \eqref{eq:lin_pol_const_real} does not need to be $n^n$, as it happens in the complex case. The question of whether the best fit for $K_n$ in \eqref{eq:lin_pol_const_real} is $n^n$ or not remained as an open problem for some time. It was already proved in \cite{BST1998} that \eqref{eq:lin_pol_const_real} holds with $K_n=n^n$ at least for $n=1,2,3,4,5$, but it was not until 2004 that a full answer to the problem was found. As a matter of fact, $K_n$ can be replaced by $n^n$ in \eqref{eq:lin_pol_const_real} for all $n\in{\mathbb N}$ (see \cite{RS2004}).\label{page:nton} Moreover, if ${\mathscr X}=\ell_1^n({\mathbb R})$ and $L_k(x_1,\ldots,x_n)=x_k$, then
$$
\|L_1\|\cdots\|L_n\|= n^n\|L_1\cdots L_n\|.
$$

The estimates \eqref{eq:lin_pol_const} and \eqref{eq:lin_pol_const_real} motivate the definition of linear polarization constants. Although $n^n$ is optimal in general in \eqref{eq:lin_pol_const} and \eqref{eq:lin_pol_const_real}, it might be improved for specific choices of ${\mathscr X}$.
\begin{definition}[\cite{BST1998}]
	We define the $n$-th linear polarization constant of the (real or complex) Banach space ${\mathscr X}$ as
	$$
	c_n({\mathscr X})=\inf\{M>0:\|L_1\|\cdots \|L_n\|\leq M\|L_1\cdots L_n\|:L_1,\ldots,L_n\in {\mathscr X}^{*}\}.
	$$
Alternatively, we also have
	$$
	c_n({\mathscr X})=1/\inf_{L_1,\ldots,L_n\in S_{{\mathscr X}^{*}}}\sup_{\|x\|=1}|L_1(x)\cdots L_n(x)|.
	$$
The linear polarization constat of ${\mathscr X}$ is defined as
	$$
	c({\mathscr X})=\limsup_n\sqrt[n]{c_n({\mathscr X})}.
	$$
\end{definition}
Interestingly, $\limsup_n$ can be replaced by $\lim_n$ in the definition of $c({\mathscr X})$ (see \cite[Proposition 4]{RS2004}). It is also worth to mention that $c({\mathscr X})=\infty$ if and only if $\dim({\mathscr X})=\infty$ (see \cite[Proposition 12]{RS2004}).\smallskip

The calculation of the constants $c_n({\mathscr X})$ for specific choices of ${\mathscr X}$ is, in most cases, a winding struggle. Only a selected bunch of linear polarization constants are known with precision. For example, we have already shown that $c_n(\ell_1^n({\mathbb K}))=n^n$. Further, if $L_1(\mu)$ is any real or complex $L_1$-space with $\dim(L_1(\mu))\geq n$, then $c_n(L_1(\mu)))=n^n$. In general, the results and the techniques required to study linear polarization constants depend strongly on whether we consider real or complex Banach spaces. The study of linear polarization constants in Hilbert spaces is a paradigmatic example of the dichotomy existing between the real and complex cases, for which reason we will devote special attention to Hilbert spaces.\smallskip

The calculation of $c_n(\ell_2^n({\mathbb K}))$ plays a central role in the theory of polarization constants since $c_n(\ell_2^n({\mathbb K}))$ is a lower bound for $c_n({\mathscr X})$ whenever ${\mathscr X}$ is an infinite-dimensional Banach space over ${\mathbb K}$. In other words	
\begin{theorem}\label{thm_AdR}\rm{\cite{RS2004}}
If ${\mathscr X}$ is an infinite-dimensional Banach space, then 	
	$$
	c_n(\ell_2^n({\mathbb K}))\le c_n({\mathscr X})\le n^\frac{n}{2}c_n(\ell_2^n({\mathbb K})),
	$$
for all $n\in{\mathbb N}$.
\end{theorem}
In 1998, Arias-de-Reyna made the following significant advancement:
\begin{theorem}\label{thm:AdR}\rm{\cite{AdR1998}}
	If $x_1,\ldots,x_n$ are unit vectors in a complex Hilbert space ${\mathscr H}$ endowed with the inner product $\langle\cdot,\cdot\rangle$, then
	\begin{equation}\label{eq:AdR}
	\sup_{\|x\|=1}|\langle x,x_1\rangle\cdots \langle x,x_n\rangle|\ge n^{-\frac{n}{2}}.
	\end{equation}
	In other words, $c_n({\mathscr H})\leq n^\frac{n}{2}$. Further,
	$c_n(\ell_2^n({\mathbb C}))=n^\frac{n}{2}$ for all $n\in{\mathbb N}$ and therefore $c_n({\mathscr H})=n^\frac{n}{2}$ whenever $\dim({\mathscr H})\ge n$.
\end{theorem}
The proof of Arias-de-Reyna remarkable result relies on complex Gaussian variables and cannot be adapted to $\ell_2^n({\mathbb R})$. As a matter of fact, the question of whether $c_n(\ell^n_2({\mathbb R}))=n^\frac{n}{2}$ remains as an open problem nowadays despide the efforts of many mathematicians. It is important to observe that \eqref{eq:AdR} follows from the so-called complex plank problem described on page~\pageref{page:plank}.\smallskip

Several works have been devoted to establish a real version of Theorem \ref{thm:AdR}, but no complete success have been achieved so far. We present below some results related to the generalization of Theorem \ref{thm:AdR} to the real case. \smallskip

First we have to say that $c_n(\ell^n_2({\mathbb R}))=n^\frac{n}{2}$ at lest for $n\leq 5$ (see \cite[Proposition 15]{RS2004} and \cite[Theorem 2]{PR2004}). In fact \cite[Theorem 2]{PR2004} actually shows that $c_n(\ell^d_2({\mathbb R}))=n^\frac{n}{2}$ whenever $n\leq \min\{d,5\}$.\smallskip

Now, if $x_1,\ldots,x_n$ are unit vectors in $\ell_2^n({\mathbb R})$ and $0\leq \lambda_1\leq \ldots\leq \lambda_n$ are the eigenvalues of the Gram matrix, i.e., the positive definite Hermitian matrix $A=[\langle x_j,x_k\rangle]_{1\leq j,k\leq n}$, then
\begin{enumerate}[$\bullet$]
	\item Marcus, 1997 (see \cite{Marcus1997} and \cite{RS2004}):
	$$
	\sup_{\|x\|=1}|\langle x,x_1\rangle\cdots \langle x,x_n\rangle|\ge (\lambda_1/n)^{n/2}.
	$$
	\item Matolcsi, 2005 (see \cite{Matolcsi2005_1}): If $x_1,\ldots,x_n$ are linearly independent then
	$$
	\sup_{\|x\|=1}\left|\langle{x}, x_{1}\rangle \cdots \langle{x},
	x_{n}\rangle\right|\geq
	\left(\frac{n}{\frac{1}{\lambda_1}+\cdots+\frac{1}{\lambda_n}}\right)^{n/2}n^{-n/2}.
	$$
	\item Matolcsi, 2005 (see \cite{Matolcsi2005_2}):
	$$
	\sup_{\|x\|=1}\left|\langle{x}, x_{1}\rangle \cdots \langle{x},
	x_{n}\rangle\right|\geq \sqrt{\lambda_{1}\cdots\lambda_n}n^{-n/2}.
	$$
	\item Mu\~noz, Sarantopoulos and Seaone, 2010 (see \cite{MUS}):
	 $$
		\sup_{\|x\|=1}\left|\langle{x}, x_{1}\rangle \cdots \langle{x},
		x_{n}\rangle\right|\geq
		\max\left\{\left(\lambda_{1}/n\right)^{n/2},\,
		\left(1/\lambda_{n}n\right)^{n/2}\right\}.
	$$
\end{enumerate}
On the other hand, the reader may find of interest the following evolution of the bounds known on $c_n(\ell_2^n)$:
\begin{enumerate}[$\bullet$]
	\item Litvak, Milman and Schechtman, 1998 (see \cite{LMS1998}):
	$$
	c_n(\ell_2^n({\mathbb R}))\leq(4e^{2\gamma} n)^\frac{n}{2},
	$$
	where $\gamma$ is the Euler-Mascheroni constant and $4e^{2\gamma}\approx 12.6892$.
	\item Garc\'{\i}a-V\'{a}zquez and Villa, 1999 (see \cite{GV1999}):
	$$
	c_n(\ell_2^n({\mathbb R}))\leq	(2e^\gamma n)^\frac{n}{2},
	$$
	where $2e^\gamma\approx 3.5622$.	
	\item A straightforward use of complexifications yields:
	$$
	c_n(\ell_2^n({\mathbb R}))\leq	(2n)^{\frac{n}{2}}/4< (2n)^{\frac{n}{2}}.
	$$
	\item Frenkel, 2008 (see \cite{Frenkel}):
	$$
	c_n(\ell_2^n({\mathbb R}))\leq	\left(\frac{3\sqrt{3}}{e} n\right)^\frac{n}{2},
	$$
	where	$\frac{3\sqrt{3}}{e}\approx 1.9115$.
	\item Mu\~noz, Sarantopoulos and Seoane, 2010 (see \cite{MUS}):
	$$
	c_n(\ell_2^n({\mathbb R}))\leq	n(\sqrt{2}n)^\frac{n}{2},
	$$
	for sufficiently lage $n$'s.
\end{enumerate}

Linear polarization constants have also been estimated for other Banach spaces, producing different results in real and complex settings. We have already mentioned that
	\begin{enumerate}[$\bullet$]
		\item $c_n({\mathscr H})=n^\frac{n}{2}$ for any complex Hilbert space with $\dim({\mathscr H})\ge n$ (see \cite{AdR1998}).
		\item $c_n(L_1(\mu))=n^n$ for any real or complex $L_1$ space with $\dim(L_1(\mu))\ge n$ (see \cite{BST1998}).
	\end{enumerate}
Other estimates and exact values of various linear polarization constants are listed below:
\begin{enumerate}[$\bullet$]
	\item Kro\'{o} and Pritsker, 1999 (see \cite{KP1999}):
	$$
	c_n(\ell_\infty^2({\mathbb C}))=2^{n-1}.
	$$
	\item R\'{e}v\'{e}sz and Sarantopoulos, 2004 (see \cite{RS2004}): If $p,q\geq 1$ with $1/p+1/q=1$, the complex $L_p(\mu)$ satisfies
	$$
	c_n(L_p(\mu))\leq\begin{cases}
	n^{n/p}&\text{if $1\leq p\leq 2$,}\\
	n^{n/q}&\text{if $p\geq 2$.}
	\end{cases}
	$$
	If, in addition, $\dim(L_p(\mu))\geq n$ and $1\leq p\leq 2$, then
	$$
	c_n(L_p(\mu))=n^{n/p}.
	$$
	\item R\'{e}v\'{e}sz and Sarantopoulos, 2004 (see \cite{RS2004}): If $n,d\in{\mathbb N}$, then
	$$
	c_n(\ell_{1}^d({\mathbb K}))= \max_{\substack{k_1+ \cdots +k_d=n \\ k_i\geq0}}
	\frac{n^n}{k_1^{k_1} \cdots k_d^{k_d}}=\prod_{l=0}^{d-1}{\Big(
		\frac {n} {[\frac{n+l}{d}]} \Big)}^{[\frac{n+l}{d}]}.
	$$
	In particular, if $n=m \cdot d$, then
	$$
	c_n(\ell_{1}^d({\mathbb K}))= d^n.
	$$
	\item Anagnostopoulos and R\'{e}v\'{e}sz, 2006 (see \cite{AR2006} and \cite{MUS}):
	$$
	c_n(\ell_2^2({\mathbb R}))=2^{n-1}.
	$$
	\item Anagnostopoulos and R\'{e}v\'{e}sz, 2006 (see \cite{AR2006}):
	$$
	K e^\frac{n}{2}\leq c_n(\ell_2^2({\mathbb C}))\leq M e^\frac{n}{2}
	$$
	where $0<K<M$.
	\item From the previous two results it is easily seen that
	$$
	c(\ell_2^2({\mathbb R}))=2\quad\text{whereas}\quad c(\ell_2^2({\mathbb C}))=\frac{e}{2}.
	$$
	\item To describe the asymptotic growth of a sequence, in particular $c(\ell_p^d({\mathbb K}))$, we shall use the standard symbols $\asymp$ and $\prec$. Observe that for any two sequences of real numbers $(a_d)$ and $(b_d)$, $a_d\prec b_d$ means that $a_d\leq Lb_d$ for some positive constant $L$, whereas $a_d\asymp b_d$ means that $a_d\prec b_d$ and $b_d\prec a_d$. Then
	$$
	c(\ell_p^d({\mathbb K}))\asymp
	\begin{cases}
	\sqrt[p]{d},&\text{if $1\leq p\leq 2$},\\
	\sqrt{d},&\text{if $p\geq 2$},
	\end{cases}
	$$
	and
	$$
	\sqrt{d}\prec c_n(\ell_\infty^d({\mathbb K}))\prec d^{\frac{1}{2}+\varepsilon},\quad\text{for all $\varepsilon>0$}.
	$$
	The previous asymptotic estimates are due to Carando, Pinasco and Rodr\'iguez (see \cite{CPR2017}).
\end{enumerate}



\subsection[Norm of products of polynomials]{Norm of products of polynomials: The factor problem for homogeneous polynomnials}	\ \smallskip

If ${\mathscr X}$ is a Banach space and $P_j\in {\mathcal P}(^{k_j}{\mathscr X})$ with $1\leq j\leq n$ are $n$ homogeneous polynomials on ${\mathscr X}$, then it is straightforward to show that
	$$
	\|P_1\cdots P_n\|\leq \|P_1\|\cdots\|P_n\|.
	$$
At the other end of the scale, it was proved in \cite{BST1998} that there exists a universal constant $M_{k_1,\ldots,k_n}>0$ depending only on the degrees $k_1,\ldots,k_n$ such that
	$$
	\|P_1\|\cdots\|P_n\|\leq M_{k_1,\ldots,k_n}\|P_1\cdots P_n\|.
	$$
Estimates of this type fall within the so-called factor problem. This problem had already been studied for products of polynomials in one (complex) variable many decades earlier than the question was stated for Banach spaces by Ben\'itez, Sarantopoulos and Tonge \cite{BST1998}. \smallskip

The calculation of the best (smallest) possible value of $M_{n_1,\ldots,n_k}$ for all or specific Banach spaces has been a fruitful field of interest for many great mathematicians in the last 20 years. It was proved in \cite{BST1998} that
	$$
	M_{k_1,\ldots,k_n}=\frac{(k_1+\cdots+ k_n)^{k_1+\cdots+k_n}}{k_1^{k_1}\cdots k_n^{k_n}}
	$$
whenever ${\mathscr X}$ is any complex Banach space. For this choice of the constant $M_{k_1,\ldots,k_n}$, equality is attained in $\|P_1\|\cdots\|P_n\|\leq M_{k_1,\ldots,k_n}\|P_1\cdots P_n\|$ (see \cite[Example 1]{BST1998}) for ${\mathscr X}=\ell_1({\mathbb K})$ and
	$$
	P_j((x_i)_{i=1}^\infty)=x_{k_1+\cdots+k_{j-1}+1}\cdots x_{k_1+\cdots+k_j}.
	$$
Observe that letting $k_1=\ldots=k_n=1$, we recover the estimate
	$$
	\|L_1\|\cdots\|L_n\|\leq n^n\|L_1\cdots L_n\|
	$$
for all $L_1,\ldots,L_n\in{\mathscr X}^{*}$, showing that the factor problem for homogeneous polynomials generalizes the linear polarization problem.\smallskip

In the following we shall present other estimates related to the factor problem for homogeneous polynomials on Banach spaces. The fisrt thing that should be pointed out is the lack of known sharp estimates for arbitrary real Banach spaces. Some of the most succesful attempts to improve the constant $M_{k_1,\ldots,k_n}=\frac{(k_1+\cdots+ k_n)^{k_1+\cdots+k_n}}{k_1^{k_1}\cdots k_n^{k_n}}$ have been focused on Hilbert spaces. In 1998, Boyd and Ryan (see \cite{BR1998}) proved that
	$$
	\|P_1\|\cdots\|P_n\|\leq \frac{(k_1+\cdots+k_n)!}{k_1!\cdots k_n!}\|P_1\cdots P_n\|
	$$
for $P_j\in{\mathcal P}(^{k_j}{\mathscr H})$ $(1\leq j\leq n)$, where ${\mathscr H}$ is a complex Banach sapace. The constant, however, is not optimal. In 2012, Pinasco \cite{P2012} found the following improvement of the previous estimate
	$$
	\|P_1\|\cdots\|P_n\|\leq\sqrt{\frac{(k_1+\cdots+ k_n)^{k_1+\cdots+k_n}}{k_1^{k_1}\cdots k_n^{k_n}}}\|P_1\cdots P_n\|.
	$$
Moreover, the constant is sharp whenever the complex Hilbert space ${\mathscr H}$ satisfies $\dim({\mathscr H})\geq n$, and equality is attained for the polynomials defined by $P_i(z)=z_i^{k_i}$ for $1\leq i\leq n$, where $(z_1,\ldots,z_n)$ are the first $n$ coordinates of $z$ with respect to an orthonormal basis of ${\mathscr H}$. As for real Hilbert spaces, it is easy to derive an estimate using the last inequality in combination with the Lindenstrauss--Tzafriri complexification norm (see page \pageref{page:LT}), however, that estimate can be greatly enhanced. Actually, Malicet et al. (see \cite{Malicet}) proved in 2016 that the inequality
	$$
	\|P_1\|\cdots\|P_n\|\leq\sqrt{\frac{2^{k_1+\cdots+ k_n}\Gamma\left(k_1+\cdots+k_n+\frac{d}{2}\right)}{\Gamma\left(\frac{d}{2}\right)k_1^{k_1}\cdots k_n^{k_n}}}\|P_1\cdots P_n\|,
	$$
holds whenever ${\mathscr H}$ is a $d$-dimensional real Hilbert space and $P_j\in{\mathcal P}(^{k_j}{\mathscr H})$ $(1\leq j\leq n)$.

For complex $L_p(\mu)$ spaces we have (see \cite{CPR2013})
	$$
	\|P_1\|\cdots\|P_n\|\leq\sqrt[p]{\frac{(k_1+\cdots+ k_n)^{k_1+\cdots+k_n}}{k_1^{k_1}\cdots k_n^{k_n}}}\|P_1\cdots P_n\|.
	$$
Interestingly, equality is attained in the previous estimate for any choice of polynomials $P_j$ of degree $n_j$ ($1\leq j\leq n$) such that $P_j$ and $P_k$ do not share any common variables for $1\leq j\ne k\leq n$ (see \cite{CPR2013}).\smallskip

Some of the estimates that have appeared in this section have been recently improved in finite-dimensional spaces, real or complex. For instance, if ${\mathscr X}$ is a $d$-dimensional Banach space over ${\mathbb K}$, then (see \cite[Theorem 2.1]{CPR2017_2})
	$$
	\|P_1\|\cdots\|P_n\|\leq \frac{(C_{\mathbb K}4ed)^{\sum_{j=1}^n k_j}}{2^\frac{2}{C_{\mathbb K}}}\|P_1\cdots P_n\|,
	$$
where $C_{\mathbb R}=1$, $C_{\mathbb C}=2$ and, as usual, $P_j\in{\mathcal P}(^{k_j}{\mathscr X})$ $(1\leq j\leq n)$.\smallskip

Now, if ${\mathscr H}$ is a $d$-dimensional Hilbert space over ${\mathbb K}$, then (see \cite[Proposition 2.2]{CPR2017_2})
	$$
	\|P_1\|\cdots\|P_n\|\leq \left(\frac{e^{H_d}C_{\mathbb K}}{4}\right)^{\sum_{j=1}^n k_j}\|P_1\cdots P_n\|,
	$$
where $C_{\mathbb R}=1$, $C_{\mathbb C}=2$, $H_d=\sum_{j=1}^n\frac{1}{j}$ and $P_j\in{\mathcal P}(^{k_j}{\mathscr H})$ $(1\leq j\leq n)$.

\subsection[Real and complex plank problems]{Real and complex plank problems and their relationship with linear polarization constants}\label{PlankProblem}\ \smallskip

\begin{figure}[h!]
	\centering
	\includegraphics[height=.5\textwidth,keepaspectratio=true]{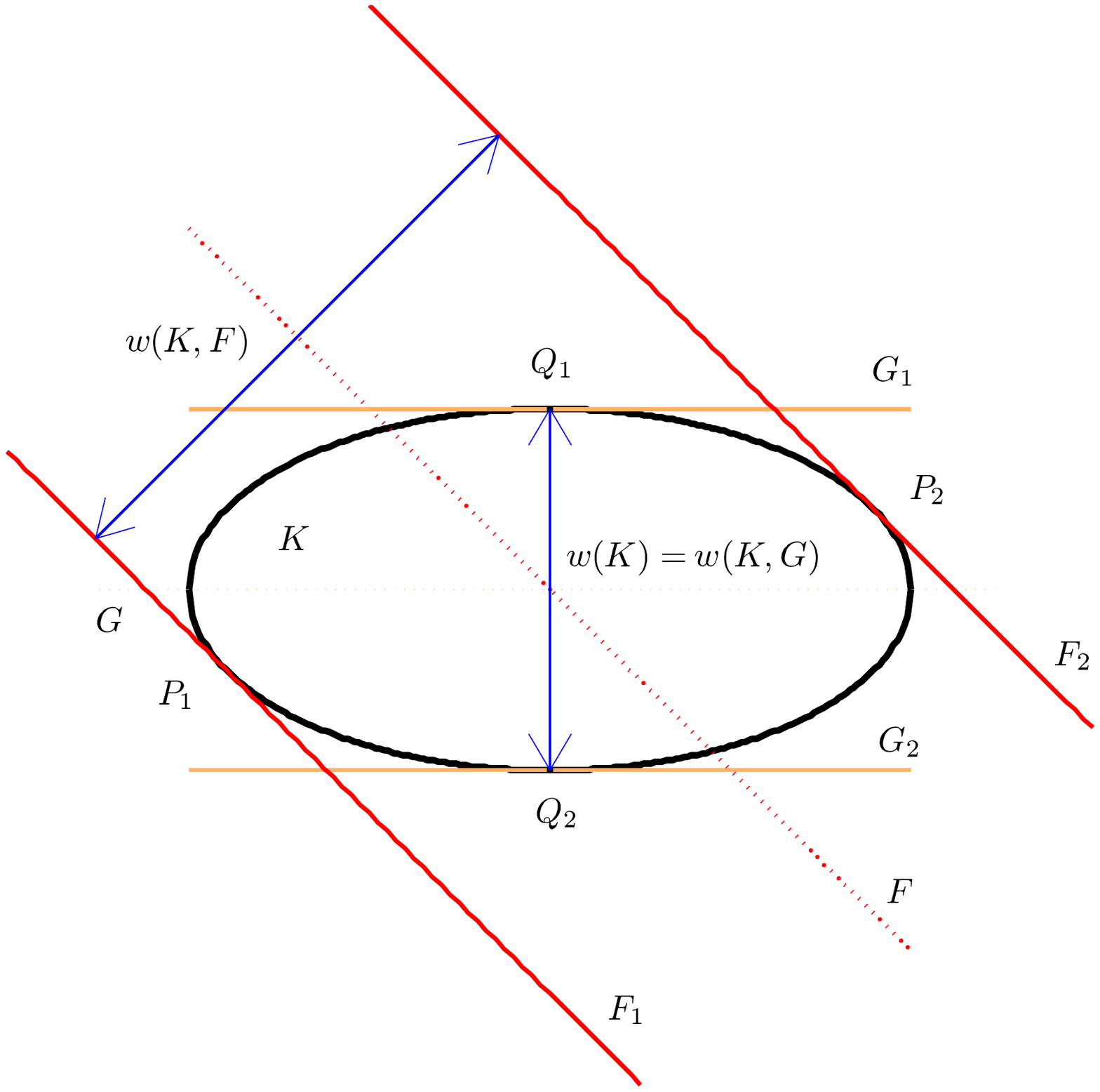}
	\caption{Width $w(K,F)$ of $K$ parallel to $F$ and minimum width of $K$}\label{fig:1}
	\end{figure}

To understand the essence of plank problems we need to introduce a few concepts.
\begin{definition}
	Let ${\mathscr X}$ be a Banach space over ${\mathbb K}$ and $K\subset {\mathscr X}$ a convex body (convex, compact set with nonempty interior). Then
	\begin{enumerate}[$\bullet$]
		\item A plank $P$ in ${\mathscr X}$ is a set of points laying between two parallel hyperplanes, i.e.,
		$$
		P=\{x\in{\mathscr X}:|f(x)-f(a)|\leq \delta\},
		$$
		where $f\in{\mathscr X}^{*}$, $\delta>0$ and $a\in{\mathscr X}$. If $f$ has norm 1, then the width of $P$ is $w(P)=2\delta$.
		\item If $F$ is a hyperplane of ${\mathscr X}$, the width $w(K,F)$ of $K$ parallel to $F$ is the distance between two supporting hyperplanes to $K$ parallel to $F$ (see Figure \ref{fig:1}).
		\item The minimum width of $K$ is $w(K)=\inf_F w(K,F)$.
		\item If $P$ is a plank parallel to the hyperplane $F$, the width of $P$ relative to $K$ is	 $w_K(P)=w(P)/w(K,F)$.
	\end{enumerate}
\end{definition}

The study of plank problems goes back to the 1930's when Tarski posed the question:
	\begin{quote}
	If $K$ is a convex body (in ${\mathbb R}^n$) covered by $n$ planks of widths $w_1,\ldots,w_n$, is it true that $w_1+\cdots+w_n\geq w(K)$?
	\end{quote}
Intuition tells us that the answer to the latter question is yes, however a formal proof of Tarski's problem is not easy. Tarski gave his own proof for a disc in ${\mathbb R}^2$ in the 1932 (see \cite{Tarski} for Tarski's original solution or \cite{K1994} for a modern exposition). Tarski's plank problem was proved in general by Bang in 1951 \cite{Bang}. At the end of his paper, Bang also posed the following strengthened version of Tarski's plank problem:
\begin{quote}
	If $K$ is a convex body covered by the planks $P_1,\ldots,P_n$, is it true that
	$w_K(P_1)+\cdots+w_K(P_n)\geq 1$?
\end{quote}
A positive answer to Bang's plank problem was found by Ball in 1991 (see \cite{Ball1991}) for convex bodies with central symmetry. Ball's solution is formulated in terms of real Banach spaces.
\begin{theorem}\label{thm:BallPlankTheorem}\rm{\cite{Ball1991}}
	If ${\mathscr X}$ is a real Banach space, $f_1,\ldots,f_n\in {\mathscr X}^{*}$ have norm $1$, and $t_1,\ldots,t_n\geq 0$ with $t_1+\cdots+t_n=1$, then 	there exists a unit vector $x\in {\mathscr X}$ such that $|f_k(x)|\geq t_k$ for
	$1\leq k\leq n$.
\end{theorem}
The previous result will be called from now on Ball's real plank theorem.
There exists a very close connection between Ball's plank theorem and linear polarization constants. This relationship is revealed by letting $t_k=\frac{1}{n}$ $(1\leq k\leq n)$ in Theorem \ref{thm:BallPlankTheorem}. Then for any $f_1,\ldots,f_n\in S_{{\mathscr X}^{*}}$ there exists $x\in S_{\mathscr X}$ such that
	$$
	|f_k(x)|\geq \frac{1}{n}.
	$$
Hence
	$$
	\inf_{f_1,\ldots f_n\in S_{{\mathscr X}^{*}}}\|f_1\cdots f_n\|\geq \frac{1}{n^n},
	$$
from which the $n$th linear polarization constant of ${\mathscr X}$ satisfies $c_n({\mathscr X})\leq n^n$ for all real Banach spaces. This estimate, which is optimal, was already mentioned in page \pageref{page:nton}, and was proved in \cite{RS2004}.\smallskip

\begin{figure}[h!]
	\centering
	\includegraphics[height=.4\textwidth,keepaspectratio=true]{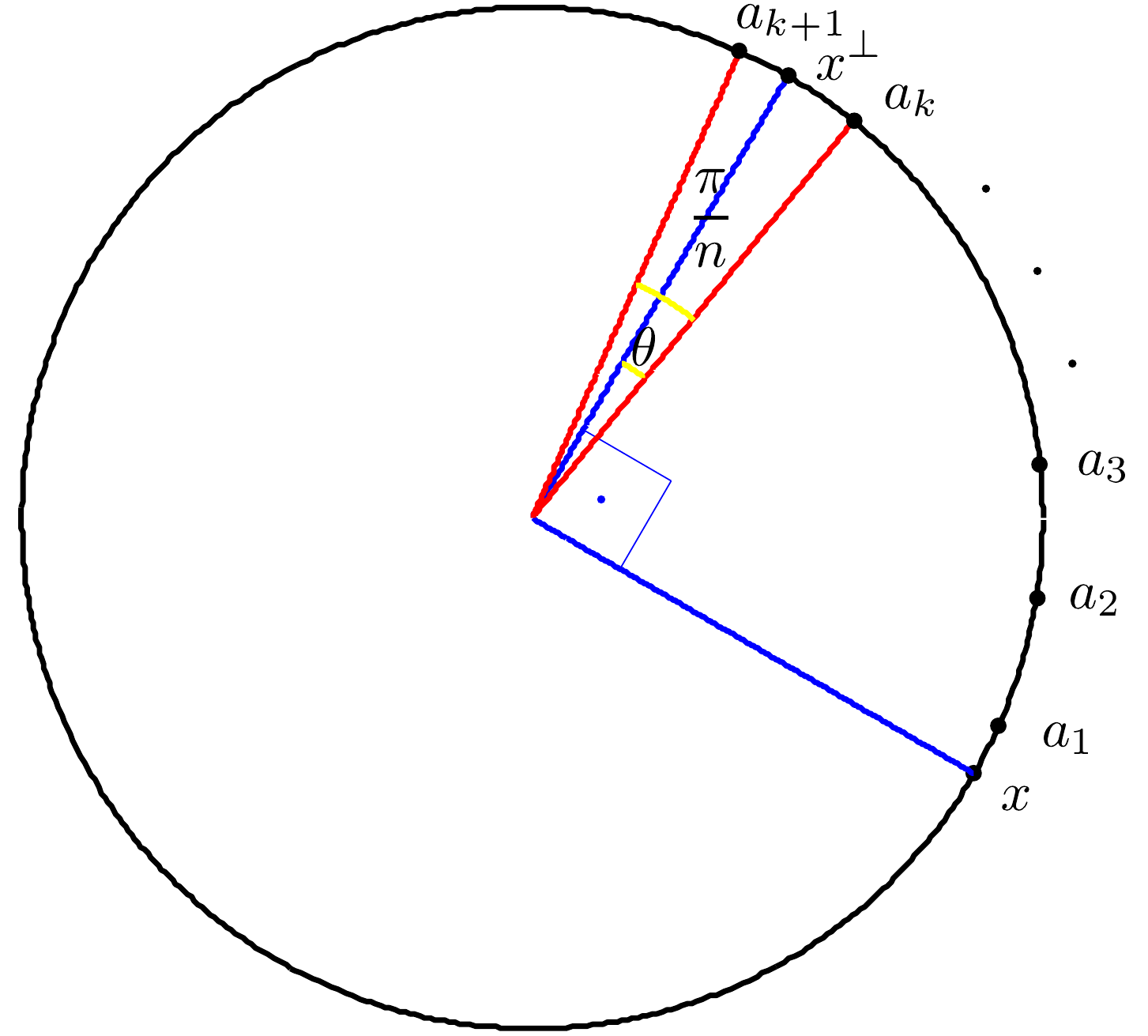}
	\caption{}\label{fig:Kirwan}
	\end{figure}

Ball's plank theorem for real Banach spaces (Theorem \ref{thm:BallPlankTheorem}) admits an analog for complex Hilbert spaces, although with a slightly different statement.
\begin{theorem}\label{thm:complexBallThm}\rm{\cite{BAA}}
	Let $\left({\mathscr H},\, \langle \cdot,\cdot \rangle\right)$ be a complex
	Hilbert space, $a_1,\ldots,a_n$ unit vectors in ${\mathscr H}$, and
	$t_1,\ldots,t_n\geq 0$ with $\sum_{k=1}^{n}{t_{k}^2}=1$.
	Then there exists a unit vector $x\in {\mathscr H}$ such that $|\langle{x},
	a_{k}\rangle|\geq t_{k}$ for $1\leq k\leq n$. In particular
	$$
	|\langle x,a_1\rangle\cdots \langle x,a_n\rangle|\geq t_1\cdots t_n.
	$$
\end{theorem}
The previous result will be named Ball's complex plank theorem from now on. Observe that putting $t_k=\frac{1}{\sqrt{n}}$ in Theorem \ref{thm:complexBallThm}, for every unit vectors $a_1,\ldots,a_n$ there exists $x$ with $\|x\|=1$ such that
	$$
	|\langle x,a_1\rangle\cdots \langle x,a_n\rangle|\geq n^{-\frac{n}{2}}.
	$$
From the previous fact we can infer straightforwardly that the $n$-th linear polarization constant of a complex Hilbert space ${\mathscr H}$ is at most $n^{\frac{n}{2}}$, or equivalently $c_n({\mathscr H})\leq n^{\frac{n}{2}}$. Equality is attained whenever $\dim({\mathscr H})\ge n$, providing an alternative proof of Theorem \ref{thm_AdR}. Unfortunately Ball's complex plank theorem is not true for real Hilbert spaces in general. Indeed, as pointed out by Kirwan \cite[p. $706$]{RT}, if we distribute $2n$ points $a_1,\ldots,a_{2n}$ within the unit circle in ${\mathbb R}^2$, then for any unit vector $x$ in the plane there is some $a_{k}$ with $1 \leq k \leq 2n$ for which
	\[
	 |\langle x,a_k\rangle|=\cos\left(\frac{\pi}{2}-\theta\right)=\sin\theta\leq \sin\frac{\pi}{n}\leq
	\frac{\pi}{n}<\frac{1}{\sqrt{2n}},
	\] (see Figure \ref{fig:Kirwan}).



\section{Real Banach algebras, real C$^*$-algebras, real J$^*$B-algebras, and real JB$^*$-triples}\label{sec: Banach algebras}

Banach algebras have been among the most studied objects in functional analysis since the beginning of the theory. A \emph{real or complex Banach algebra} is a real or complex Banach space $(\mathscr{A}, \|\cdot\|)$ equipped with an associative (and bilinear) product $\mathscr{A}\times \mathscr{A}\to \mathscr{A}$, $(a,b) \mapsto a b,$ satisfying $$\| a b\| \leq \|a\| \ \|b\|\quad \hbox{ for all } a,b\in \mathscr{A}.$$ The latter condition is clearly a link between the algebraic and the analytic structures assuring the continuity of the norm. Different substructures are obtained by adding extra hypotheses on the Banach algebra. For example, a real or complex Banach algebra $\mathscr{A}$ is called \emph{commutative} if its product enjoys the property that $a b = b a $ for all $a,b\in \mathscr{A}$. We say that $\mathscr{A}$ is unital if there exists a necessarily unique element $\11\in \mathscr{A}$ satisfying $\11 a = a = a \11$ for all $a\in \mathscr{A}$. Clearly, every complex Banach algebra is a real Banach algebra by just restricting the product by scalars to the real field. {The center of a real or complex Banach algebra $A$ (denoted by $Z(A)$) will consist in all elements $a\in A$ such that $ a z = z a $ for all $a\in A$.} \smallskip

We will see a good list of examples along with the paper. For the moment, we begin with the best-known models. For each compact Hausdorff space $K$, the spaces $C(K, \mathbb{R})$ and $C(K) = C(K, \mathbb{C})$ of all real-valued and complex-valued functions on $K$, respectively, are examples of commutative real and complex Banach algebras with respect to the supremum norm and the pointwise product. {Let $F$ stand for a closed subset of $K$. The set $C_F^{\mathbb{R}}(K) :=\{f\in C(K) : f(F)\subseteq \mathbb{R}\}$ is a closed real subalgebra of $C(K)$.} These Banach algebras always admit a unit, namely, the constant function $1$. For a locally compact Hausdorff space $\Omega$, the Banach spaces $C_0(\Omega, \mathbb{R})$ and $C_0(\Omega)$, respectively, of all real-valued and complex-valued continuous functions on $\Omega$ vanishing at infinity are examples of nonunital commutative real or complex Banach algebras when equipped with the supremum norm and the pointwise product. We observe that for every real or complex Banach space $\mathscr{X}$, the space $\mathcal{B}(\mathscr{X}),$ of all bounded linear operators on $\mathscr{X},$ is a real or complex Banach algebra with respect to the composition and the operator norm. In particular, the Banach spaces $\mathcal{M}_n(\mathbb{R})$ and $\mathcal{M}_n(\mathbb{C})$ are real and complex Banach algebras with respect to the matricial product and the operator norm. These latter examples are unital but noncommutative Banach algebras if dim$(\mathscr{X})\geq 2$.\smallskip

Despite the fact that the general strategy in this paper consists in extending the norm from a real structure to its algebraic complexification, with the unique condition that the corresponding extension preserves the same algebraic and analytic structures assumed on the real object, there exists another procedure consisting in assuming that we already have an appropriate extension and considering a suitable real subspace whose complexification is the structure from which we began. In analogy with what has been considered in the previous sections, we deal with real forms of complex Banach spaces. Suppose that $\tau: \mathscr{X}\to \mathscr{X}$ is a conjugate-linear isometry of period 2 (i.e., $\tau^2 = Id_{\mathscr{X}}$) on a complex Banach space $\mathscr{X}$. The set $$\mathscr{X}^{\tau} :=\left\{ x\in \mathscr{X} : \tau (x) =x \right\}$$ of all $\tau$-fixed points in $\mathscr{X}$ is a closed real subspace and hence a real Banach space when equipped with the restricted norm. The real Banach space $\mathscr{X}^{\tau}$ is called a \emph{real form} of the space $\mathscr{X}$. A mapping $\tau$ satisfying the above properties is called a \emph{conjugation} on $\mathscr{X}$. If $\mathscr{X}$ is regarded as the complexification of $\mathscr{X}^{\tau}$, then $\|x-i y \| = \|\tau (x +i y)\| = \|x +i y\|$ for all $x,y\in \mathscr{X}^{\tau}$, that is, $\mathscr{X}$ with its original norm is a reasonable complexification in the sense we employed in Subsection \ref{subsec: various nroms on Xc}.\smallskip

Let us illustrate this construction with an example. Suppose that $K$ is a compact Hausdorff space. By the celebrated Banach--Stone theorem, every surjective conjugate-linear isometry $\Phi : C(K)\to C(K)$ is of the form $$\Phi (f) (t) = u(t) \overline{f(\sigma(t))} \quad { (f\in C(K), t\in K),}$$ where $\sigma: K\to K$ is a homeomorphism and $u\in C(K)$ with $|u(t)| = 1$ for all $t\in K$. Hence every conjugation $\tau$ on $C(K)$ must be of the form $\tau(f) (t) = u(t) \overline{f(\sigma(t))}$ with $ u(\sigma(t)) =u(t)$ and $\sigma^2(t) =t$ for all $t\in K$. The real form $C(K)^{\tau}$ is a real Banach space whose complexification is $C(K),$ and its norm admits an extension to $C(K)$. However, the structure of the real form is, in principle, different from a $C(K)$-space. Namely, $C(K)^{\tau}$ need not be a subalgebra of $C(K)$ ---we will see later that it admits a concrete geometric Jordan structure---. Assuming $u=\11$ in $C(K)$, the real form $C(K)^{\tau}$ is a real Banach subalgebra of $C(K)$.\smallskip

If $\tau$ is a conjugation on a complex Banach algebra $\mathscr{A}$ and $\tau$ is multiplicative, then the real form $\mathscr{A}^{\tau}$ is a real Banach algebra.\smallskip

Clearly, the unit element in a unital real or complex Banach algebra $\mathscr{A}$ satisfies $\|\11\| \geq 1$. It is well known that we can renorm $\mathscr{A}$ with another Banach algebra norm in such a way that the unit element has norm one. For this purpose, we shall consider the representation of $\mathscr{A}$ into $\mathcal{B}(\mathscr{A})$ through the left and right multiplication operators. We recall that a \emph{homomorphism} (respectively, \emph{isomorphism}) between two real or complex Banach algebras $\mathscr{A}$ and $\mathscr{B}$ is a linear (respectively, bijective linear) mapping $\Psi: \mathscr{A}\to \mathscr{B}$ preserving the associative product, that is, $\Psi (a b ) = \Psi(a) \Psi (b)$ for all $a,b\in \mathscr{A}$. We can consider the linear maps $$L : \mathscr{A}\to \mathcal{B}(\mathscr{A}) \ \hbox{ and } R : A\to \mathcal{B}(\mathscr{A})$$ defined by $L(a) = L_a: \mathscr{A}\to \mathscr{A}$, $L_a(x) :=a x$ and $R(a) = R_a: \mathscr{A}\to \mathscr{A}$, $R_a(x) :=x a$, respectively. It is well known that $L$ and $R$ are two homomorphisms. These maps are called the left and right regular representations of $\mathscr{A}$ into $\mathcal{B}(\mathscr{A}),$ respectively. One of the advantages of the left (respectively, right) regular representation is that, assuming that $\mathscr{A}$ is unital, by renorming it via the norm $\|| a |\| :=\|L_a\|_{\mathcal{B}(\mathscr{A})}$ (respectively, $\|| a |\| :=\|R_a\|_{\mathcal{B}(\mathscr{A})}$), we find an equivalent algebra norm on $\mathscr{A}$ satisfying $\||\11 |\| = 1$. Henceforth, we shall assume that for each unital (real or complex) Banach algebra $\mathscr{A}$, we have $\|\11\| =1$.\smallskip

A conjugate-linear multiplicative mapping between two complex Banach algebras will be called a \emph{conjugate-linear homomorphism}. A \emph{conjugate-linear isomorphism} is a conjugate-linear bijection that is also multiplicative.\smallskip

A real or complex Banach algebra $\mathscr{A}$ without unit can be always regarded as a norm closed subalgebra of a unital Banach algebra. It suffices to consider the \emph{unitization} $\mathscr{A}_1 = \mathscr{A}\oplus \mathbb{K} \11$ with the obvious extension of the product and the norm $\|a +\lambda \11\| := \|a\| + |\lambda|$ (see \cite[Definition I.3.1]{BD}). We can also consider the left regular representation and the subalgebra of $\mathcal{B}(\mathscr{A})$ generated by $L(\mathscr{A})$ and the identity on $\mathscr{A}$. In the latter case, we have $$\|a+\lambda \11\| :=\|L_a + \lambda Id_{\mathscr{A}}\|_{\mathcal{B}(\mathscr{A})} = \sup_{\|x\|\leq 1} \| a x +\lambda x\| \quad (a\in \mathscr{A}, \lambda\in \mathbb{K}).$$

\subsection{Standard complexification of a Banach algebra}\label{subsec: standard complexification Banach algebras}\ \smallskip

From a strictly algebraic point of view, given a real Banach algebra $A$, there is only one natural extension of its product to an associative product on its complexification $A_c=A+i A$, which is defined by $$ (a+i b) (c+ i d) = a c - bd + i (a d + b c) \quad (a,b,c,d\in A).$$ Clearly, $A_c$ is commutative whenever $A$ is, and if $A$ admits a unit $\11$, the same element is a unit in $A_c$. When we regard $A$ merely as a Banach space, we can consider its Taylor complexification given in Subsection \ref{subsec: various nroms on Xc} (see page \pageref{label Taylor complexification}) whose norm is given by $$\|x+iy\|_T:=\sup_{t\in[0,2\pi]}\|x\cos t-y\sin t\| \quad (x+iy \in A_c).$$

By considering the left regular representation of ${A}_c$ as a subalgebra of $\mathcal{B}(A_c, \|\cdot\|_{T})$ (with the operator norm given by the Taylor complexification), we define a Banach algebra norm on $A_c$ given by this representation, that is, $$\|a+i b\|_{T,a} := \|L_{a+ib}\|_{\mathcal{B}(A_c, \|\cdot\|_T)} = \sup_{\|(x+i y)\|_T\leq 1}\|(a+i b) (x+i y)\|_T.$$ This complex Banach algebra $(A,\|\cdot\|_{T,a})$ (respectively, this norm $\|\cdot\|_{T,a}$) is called the \emph{standard complexification} of ${A}$ (respectively, the standard norm) in references like \cite{LUM}. If $\|\cdot\|_r$ is any reasonable complete norm on the complex space $A_c$, for example, $\|x+i y \|_p^p :=\|x\|^p + \|y\|^p$ with $1\leq p<\infty,$ and $\|x+i y \|_{\infty} :=\max\{\|x\|, \|y\|\}$ ---We recall that all reasonable norms on $A_c$ are equivalent to the Taylor norm (see page \pageref{eq any reasonable norm is equivalent to the T norm})---, then we can reproduce the above procedure to obtain a complex Banach algebra norm $\|\cdot\|_{r,a}$ on $A_c$. All these complex Banach algebra norms are reasonable and equivalent to $\|\cdot\|_{T,a}$.\smallskip

As narrated in the monographs \cite[\S 13]{BD} and \cite[\S 2.1]{LI}, there is another method to extend the norm of a real Banach algebra $(A,\|\cdot\|)$ to a norm on the complex Banach algebra $A_c$. Namely, let $\mathcal{B}_{A}$ denote the closed unit ball of $A$ and let $V$ denote the absolutely convex hull of the set $\mathcal{B}_{A}\times\{0\}$ in $A_c$, that is,
$$V= |co| \left( \mathcal{B}_{A}\times\{0\} \right)= \left\{ \sum\limits_{j} \alpha_j a_j : a_j \in \mathcal{B}_{A}\times\{0\}, \ \alpha_j\in \mathbb{C} \hbox{ with } \sum\limits_j |\alpha_j|\leq 1 \right\},$$
which is an absorbent set in $A_c$. The Minkowski functional associated with $V$ defines a reasonable, complete algebra norm $\|\cdot\|_m$ on $A_c$, whose open unit ball is precisely $V$ and its restriction to $A$ coincides with $\|\cdot\|$ and satisfies $$\max\{\|a\|,\|b\|\}\leq \|a+i b\|_m \leq 2 \max\{\|a\|,\|b\|\}$$ for all $a,b\in A$ (see \cite[Proposition I.13.3]{BD} or \cite{LI}). The reader interested in additional results on the complexification of a normed real algebra can also consult \cite[\S 1.1.5]{Cabrera-Rodriguez-vol1}. As remarked in the just quoted monograph, ``\emph{Due to the power of complex methods, the possibility of regarding (isometrically) any real normed algebra as a real subalgebra	of a normed complex algebra becomes a relevant fact.}'' The complexification method in \cite{Cabrera-Rodriguez-vol1} by means of the projective tensor norm, is precisely the one by Bonsall and Duncan \cite{BD}, and has the advantage that it works without problems in	the non-associative setting.\smallskip

All the above procedures define equivalent reasonable algebra norms on the complexification of a real Banach algebra $A$. It should be noted that the processes of unitization and complexification on a real Banach algebra can be interchanged, and the resulting algebra does not change. If $A_c$ is the complexification of $A$ equipped with a reasonable complete Banach algebra norm $\|\cdot\|$, then we can define a conjugation $\tau$ on $A_c$ given by $\tau (a+i b)= \overline{a+i b} = a- i b $, which is clearly a period-2 isometry, because the norm on the complexification is reasonable. Furthermore, the conjugation $\tau$ is a conjugate-linear homomorphism on $A_c$, and $A = A_c^{\tau}$ is a real Banach subalgebra of $A_c$.\smallskip

{Each complex Banach algebra $\mathcal{A}$ can be always regarded as a real Banach algebra, $\mathcal{A}_{\mathbb{R}},$ by just restricting the product by scalars to the real field. Conversely, it is interesting to have tools to determine if a real Banach algebra is obtained from a complex one in this way. According to the standard terminology, a real (normed) algebra is said to be of \emph{complex type} if it is possible to extend the scalar multiplication to complex scalars so that the algebra becomes a complex (normed) algebra under an equivalent norm (cf. \cite[Definition 6.1]{Inglestam64}). The following technical characterization of real normed algebras of complex type was established by L. Ingelstam with tools developed by I. Kaplansky \cite{Kap49} and a complex norm given by a formula close to the Taylor complexification.
	
\begin{theorem}\label{t Ingelstam  characterization of complex type}{\rm\cite[Proposition 6.2 and Corollary 6.3]{Inglestam64}} A real normed algebra $A$ is of complex type if and only if there exists a continuous linear operator $J$ on $A$ satisfying: \begin{enumerate}[$(a)$] \item $J$ is an $A$-module homomorphism, that is, $$J(ab) =J(a) b = a J(b), \hbox{ for all } a,b\in A;$$
\item $-J^2$ is the identity map on $A$.
\end{enumerate} Furthermore, the equivalent Banach algebra complex norm is given by $$\|| x |\|= \max_{\theta\in \mathbb{R}} \|\cos(\theta) x + \sin(\theta) J(x)\|\  \ (x\in A).$$
		
\noindent Consequently, a real (normed) algebra $A$ with identity $\mathbf{1}$ is of complex type if and only if there exists an element $\boldsymbol{\iota}$ in the center of $A$, satisfying $\boldsymbol{\iota}^2 = -\mathbf{1}$.
\end{theorem}
	
The mapping $J$ in the previous theorem is called a complex multiplication. Propositions 2.1 and 2.2 in \cite{Inglestam66} prove that every complex multiplication on a real normed algebra $A$ is automatically continuous in any of the following cases:\begin{enumerate}[$(a)$] \item $A$ is a real Banach algebra with an approximate identity;
\item the set of (left or right) topological divisors of $0$ is not all of $A$. 	
\end{enumerate}

Whether the original norm of a real Banach algebra is not only equivalent to a complex Banach algebra norm, but it is itself a complex norm is another type of question. We shall add some answer.

\begin{proposition}\label{p original real norm is a complex norm} Let $A$ be a real normed algebra with norm $\|\cdot\|$. Then we can define a product by complex scalars on $A$ making the latter a complex normed algebra for its original norm if and only if there exists a continuous linear operator $J$ on $A$ satisfying: \begin{enumerate}[$(a)$] \item $J$ is an $A$-module homomorphism, that is, $$J(ab) =J(a) b = a J(b) \hbox{ for all } a,b\in A;$$
\item $-J^2$ is the identity map on $A$;
\item For each real $\theta$ the mapping $\cos(\theta) Id_{A} + \sin(\theta) J$ is a non-expansive mapping on $A$.
\end{enumerate}
\end{proposition}

\begin{proof} For the ``only if'' part we observe that if there exists a product by complex scalars making $(A,\|\cdot\|)$ a complex normed algebra, by defining $J(a)= i a$ the first two proerties are clear, and for the last one $$\| (\cos(\theta) Id_{A} + \sin(\theta) J) (a)\| = \| (\cos(\theta) + i \sin(\theta)) a \| = \|a\| \ \ (a\in A, \theta \in \mathbb{R}).$$
	
For the sufficient implication, it is clear that defining $(\alpha + i \beta) a = \alpha a + \beta J(a)$ ($\alpha+i \beta\in \mathbb{C}$, $a\in A$), we get a structure of complex algebra on $A$. It remains to prove that the original norm is a complex norm. Since for each $\theta \in \mathbb{R}$, by hypothesis, we have $\|\cos(\theta) Id_{A} + \sin(\theta) J\|\leq 1,$ the linear mapping $\cos(\theta) Id_{A} + \sin(\theta) J$ is a bijection with inverse $ \cos(\theta) Id_{A} - \sin(\theta) J,$ which is also non-expansive, we deduce that $\cos(\theta) Id_{A} + \sin(\theta) J$ is a linear isometry for all real $\theta$. Therefore, for a non-zero complex number $\alpha + i \beta$ we write $\alpha + i \beta = |\alpha + i \beta| (\cos(\theta) + i \sin(\theta))$ to get  $$\begin{aligned} \|(\alpha + i \beta ) a \| &= \| (\alpha + i \beta ) a \| = |\alpha + i \beta| \| \cos(\theta) a + \sin(\theta) J (a)\|\\
&= |\alpha + i \beta| \| (\cos(\theta) Id_{A}  + \sin(\theta) J ) (a)\|  = |\alpha + i \beta| \|a\|.	 
\end{aligned} $$   	
\end{proof}

We shall see later (see Theorem \ref{t Inglestam real Cstar algebra} below) that in the setting of real C$^*$-algebras, the Gelfand-Naimark axiom is a powerful geometric tool to simplify the conclusion of Theorem \ref{t Ingelstam  characterization of complex type} and Proposition \ref{p original real norm is a complex norm}.}

An element $a$ in a real or complex unital algebra $\mathscr{A}$ is called \emph{invertible} if there exists $b$ in $\mathscr{A}$ with $a b = b a =\11$. This element $b$ is unique, it is called the \emph{inverse} of $a$ in $\mathscr{A}$, and it will be denoted by $a^{-1}$. If ${A}$ is a unital real Banach algebra and $A_c$ denotes its complexification, then it is easy to check that the set $A^{-1}$ of all invertible elements in $A$ coincides with the intersection of $A$ with the set $A_c^{-1}$ of all invertible elements in $A_c$. Therefore, the usual topological properties of $A^{-1}$ and of the inverse mapping hold in the real setting, too.\smallskip

If $A$ is a unital real Banach algebra, then the \emph{spectrum} of an element $a\in A,$ $\sigma_{A}(a)$, is defined as the spectrum of $a$ in the complexification of $A$, that is, \begin{equation}\label{eq def of spectrum in the real setting} \sigma_A (a) := \sigma_{A_c} (a) =\{\lambda\in \mathbb{C} : a-\lambda \11 \notin A_c^{-1}\}.
\end{equation} If $A$ is not unital, then the spectrum of an element $a\in A$ is defined as the spectrum of this element in the unitization of $A$, which is precisely the spectrum of the element in the complexification of $A$. Therefore, by the celebrated Gelfand theorem, the spectrum of each element is a nonempty compact subset of the complex plane, bounded by the norm of the element. As in the case of operators (see page \pageref{first results spectrum}), if we define the spectrum of an element in a unital real Banach algebra $A$ in terms of real numbers and invertible elements in $A$, then we might find an empty set. So, the natural definition for the spectrum in the setting of real Banach algebras is the one given in \eqref{eq def of spectrum in the real setting}.\smallskip

Since the natural conjugation $a+i b \mapsto \overline{a+i b} = a- ib $ on the complexification, $A_c$, of a real Banach algebra $A$ is a conjugate-linear unital homomorphism, it is not hard to see that $$\sigma_{A_c} \left(\overline{a+ i b }\right) = \sigma_{A_c} (a- i b ) = \overline{\sigma_{A_c} (a+i b)}\quad (\hbox{for all } a+i b \in A_c),$$ and consequently \begin{equation}\label{eq the spectrum in the real setting is self adjoint} \sigma_{A} (a) = \sigma_{A_c} \left( a \right) = \overline{\sigma_{A_c} (a)} = \overline{\sigma_{A} (a)}\quad (\hbox{for all } a \in A).
\end{equation}
By \cite[Proposition 1.1.100]{Cabrera-Rodriguez-vol1} we also know that \begin{equation}\label{eq spectrum in complexification} \sigma_{A_c} (a)= \sigma_{A} (a) = \{\alpha +i\beta : \alpha,\beta \in \mathbb{R} \hbox{ such that } (a-\alpha \textbf{1})^2+\beta^2 \textbf{1} \notin A_c^{-1}\},
\end{equation} for all $a\in A$.\smallskip

The \emph{spectral radius} of an element $a\in A$ is defined as the corresponding spectral radius in the complexification, that is,
$$r(a) = r_A(a) = r_{A_c}(a) = \max\{ |\lambda| : \lambda \in {\rm sp}(a)\}.$$ Since the famous Gelfand--Beurling formula holds for every complex Banach algebra (see \cite[Theorem 3.2.8]{Aupe91}), we conclude that the same identity is also true for real Banach algebras, that is, $$r(a) = \lim_{n\to \infty} \|a^n\|^{\frac1n} =\max\{|\alpha + i \beta| : \alpha,\beta \in \mathbb{R} \hbox{ such that } (a-\alpha \textbf{1})^2+\beta^2 \textbf{1} \notin A_c^{-1}\} .$$

Before dealing with more concrete structures, like real and complex C$^*$-algebras, we revisit some results requiring a simple background. A fascinating achievement in the theory of Banach algebras is the so-called Gleason--Kahane--\.{Z}elazko theorem.

\begin{theorem}\label{t GKZ}{\rm(Gleason--Kahane--\.{Z}elazko theorem \cite{Gle,KaZe,Ze68})} Let $F:\mathscr{A} \to \mathbb{C}$ be a nonzero linear mapping, where $\mathscr{A}$ is a complex Banach algebra. Then the following statements are equivalent:\begin{enumerate}[$(a)$]\item $F(a)\in {\rm sp}(a)$ for every $a\in \mathscr{A}$;
\item $F$ is unital if $\mathscr{A}$ is unital or admits a unital extension to the unitization of $\mathscr{A}$ and maps invertible elements to invertible elements;
\item $F$ is multiplicative.
\end{enumerate} The mapping $F$ is continuous if it satisfies any of the previous equivalent conditions.
\end{theorem}

This is nowadays one of the fundamental contributions in functional analysis and the theory of complex Banach algebras, and it is contained in most reference books (see, for example, \cite[Theorem III.10.9]{Rudin73}, \cite[Theorem 2.4.13]{Palmer}, or \cite[Theorem II.17.7]{BD}). The Gleason--Kahane--\.{Z}elazko theorem still is a pole of attraction (see, for example, \cite{MashreRansRans2015,MashreRansRans2018,RoiStern81,TouBrits2020,TouSchuBrits2017,TouSchuBrits2018}). However, during its early years, its importance was disputed. For example, in \cite[p. 25]{Rudin73}, it was affirmed that ``This striking result has apparently found no interesting applications as yet'' ---nothing farthest from its real role in mathematics---. Subsequent years have witnessed a whole explosion of new ideas and applications coming out induced by this important result. Indeed, the Gleason--Kahane--\.{Z}elazko theorem was applied by Cabello S{\'a}nchez and Moln{\'a}r \cite{CabMol2002} while studying the reflexivity of the isometry group and the automorphism group of uniform algebras and topological algebras of holomorphic functions, by Cabello S{\'a}nchez \cite{CabSan2004} for investigation of the Banach algebras $L_{\infty}(\mu)$ for various measures $\mu$ , and by Jim{\'e}nez-Vargas, Morales Campoy, and Villegas-Vallecillos \cite{JimVMorVill2010} in exploration of the algebraic reflexivity of the isometry group of some spaces of Lipschitz functions.\smallskip

It is worth noting that Choda and Nakamura \cite{ChodaNakamura71} gave two short proofs of the Gleason--Kahane--\.{Z}elazko theorem in the special case in which $\mathscr{A}$ is a C$^*$-algebra, while a simple proof for complex Banach algebra with a hermitian involution was established by Ch\={o} \cite{Cho}.\smallskip

The Gleason--Kahane--\.{Z}elazko theorem is not valid for real Banach algebras. For example, let ${A}= C([0,1],\mathbb{R})$ be the real algebra of all continuous real-valued functions on $[0,1]$ and let $F: A\to \mathbb{R}$, $ F(f):=\frac12 (f(0)+f(1))$ ($f\in A$). Since $$\min\{f(0),f(1)\}\leq F(f)=\frac12 (f(0)+f(1)) \leq \max\{f(0),f(1)\} \quad (f\in {A}),$$ the intermediate value theorem implies that $F(f)\in {\rm sp}(f)$ for all $f\in {A}$, but it can be easily checked that $F$ is not multiplicative (see \cite{SS}). Another example can be given by the mapping $G: A\to \mathbb{R}$, $G(f) = \int_0^1 f(t) dt$. By the mean value theorem, $G(f)$ lies in ${\rm sp}(f)$ for all $f\in A$, and clearly, $G$ is not multiplicative.\smallskip

Despite the obstacles in the real setting, Kulkarni established the following version of the Gleason--Kahane--\.{Z}elazko theorem for real Banach algebras, which was originally proved by an ingenious application of functional calculus and Hadamard's factorization theorem.

\begin{theorem}\label{t GKZ Kulkarni}{\rm(Kulkarni--Gleason--Kahane--\.{Z}elazko theorem, \cite{Kul84})} Let $F: A\to \mathbb{C}$ be a nonzero linear map, where $A$ is a unital real Banach algebra. Then the following statements are equivalent: \begin{enumerate}[$(i)$]\item $F$ is multiplicative;
\item $F(\11) = 1$ and $F(a)^2 + F(b)^2 $ lies in ${\rm sp}(a^2 + b^2)$ for all $a,b\in A$ with $ab = ba$;
\item $F(\11) = 1$ and $F(a)^2 + F(b)^2 \neq 0$ for all $a,b\in A$ with $ab=ba$ and $a^2+b^2$ invertible.
\end{enumerate}
\end{theorem}

The original Gleason--Kahane--\.{Z}elazko theorem can be derived from the previous result via the following ingenious idea: Let $\mathscr{A}$ be a complex Banach algebra and let $F :A \to \mathbb{C}$ a linear mapping satisfying statement $(b)$ in Theorem \ref{t GKZ}. Given two elements $a$ and $b$ in $\mathscr{A}$ such that $a b = ba $ and $a^2 + b^2$ is invertible, the identity $a^2 + b^2 = (a+i b) (a-i b)$ implies that $(a+i b)$ and $(a-i b)$ are invertible, and hence $F(a)^2 + F(b)^2 = F(a+i b) F(a-i b)$ must be a nonzero complex number.\smallskip

If instead of studying the algebraic reflexivity of the isometry group and local isometries and automorphisms, we are interested in 2-local isometries and automorphisms (or their weak versions), in the way introduced by \v{S}emrl \cite{Semrl97} and Larson and Sourour \cite{LarSou}, like in the studies conducted by Hatori et al. \cite{HaMiOkTak07} on 2-local isometries and 2-local automorphisms between uniform algebras, on weak-2-local isometries between uniform and Lipschitz algebras by Li et al. \cite{LiPeWangWang} and by Jim{\'e}nez Vargas and Villegas-Vallecillos \cite{JimVVill20}, then we realize that the appropriate tool is the following theorem due to Kowalski and S{\l}odkowski. We omit additional details for the sake of brevity.

\begin{theorem}\label{t KS}{\rm(Kowalski--S{\l}odkowski theorem \cite{KoSlod})} Let $\mathscr{A}$ be a complex Banach algebra and let $\Delta : \mathscr{A}\to \mathbb{C}$ be a mapping satisfying $\Delta (0)=0$ and $$\Delta (x) - \Delta (y) \in \sigma (x-y)$$ for every $x,y\in \mathscr{A}$. Then $\Delta$ is linear and multiplicative.
\end{theorem}

The following spherical versions of the Gleason--Kahane--\.{Z}elazko and Kowalski--S{\l}odkowski theorems, which are suitable tools to study weak-2-local isometries can be found in \cite{LiPeWangWang}. From now on, we write $\mathbb{T}$ for the unit sphere of $\mathbb{C}$.

\begin{theorem}\label{t GKZ sphere}{\rm(Spherical Gleason--Kahane--\.{Z}elazko theorem, \cite[Proposition 2.2]{LiPeWangWang})} Let $F:\mathscr{A} \to \mathbb{C}$ be a linear mapping, where $\mathscr{A}$ is a unital complex Banach algebra. Suppose that $F(a)\in \mathbb{T} \ {\rm sp}(a)$ for every $a\in \mathscr{A}$. Then the mapping $\overline{F(1)} F$ is multiplicative.
\end{theorem}

\begin{theorem}\label{t KS sphere}{\rm(Spherical Kowalski--S{\l}odkowski theorem, \cite[Proposition 3.2]{LiPeWangWang})} Let $\mathscr{A}$ be a unital complex Banach algebra and let $\Delta : \mathscr{A}\to \mathbb{C}$ be a mapping satisfying the following properties: \begin{enumerate}[$(a)$]\item $\Delta$ is 1-homogeneous;
\item $\Delta (x) - \Delta (y) \in \mathbb{T} \ \sigma (x-y)$ for every $x,y\in \mathscr{A}$.
\end{enumerate} Then $\Delta$ is linear, and there exists $\lambda_0\in \mathbb{T}$ such that $\lambda_0 \Delta$ is multiplicative.
\end{theorem}

An interesting contribution due to Oi (see \cite{Oi19}) shows that by replacing hypothesis $(a)$ in the previous theorem by the condition $\Delta(0)=0$ we can get a similar conclusion to that in the Kowalski--S{\l}odkowski theorem.

\begin{theorem}\label{t KS sphere 0}{\rm\cite{Oi19}} Let $\mathscr{A}$ be a unital complex Banach algebra and let $\Delta : \mathscr{A}\to \mathbb{C}$ be a mapping satisfying the following properties: \begin{enumerate}[$(a)$]\item $\Delta (0) =0;$
		\item $\Delta (x) - \Delta (y) \in \mathbb{T} \ \sigma (x-y)$ for every $x,y\in \mathscr{A}$.
	\end{enumerate} Then $\Delta$ is is complex-linear or conjugate-linear and $\overline{\Delta(\textbf{1})} \Delta$ is multiplicative.
\end{theorem}

The Gleason--Kahane--\.{Z}elazko and Kowalski--S{\l}odkowski theorems are now influencing on the developing of new problems in the fruitful line of preservers. We shall see some related results after presenting the basic background on C$^*$-algebras.

\subsection{Division real Banach algebras}\ \smallskip

This is an appropriate moment to introduce another example of a real Banach algebra. We refer to one of the few mathematical models about which we know the exact date and place in which they were invented. We are speaking about Hamilton's quaternions, whose origins were explicitly dated in a letter by Hamilton to his friend and fellow mathematician Graves, in which he wrote ``\emph{And here (at Brougham Bridge) there dawned on me the notion that we must admit, in some sense, a fourth dimension of space for the purpose of calculating with triples ... An electric circuit seemed to close, and a spark flashed forth.}'' (Dublin, 16th of October 1843).\smallskip

The algebra of \emph{quaternions}, $\mathbb{H}$,\label{def quaternions} is the four-dimensional real linear space with basis $\{\11,i,j,k\}$ and associative multiplication defined by $$\11 \hbox{ is the identity and } i^2=j^2= k^2= -1 = ijk.$$ All the other possible products follow from these identities, for example, $i j = (i j k) (-k) = -(-k) =k$, $j k = (-i) (ijk) = i,$ $j i = j (ijk)(kj)= - (j k) j = -i j =-k.$ The algebra $\mathbb{H}$ is noncommutative. When equipped with the Euclidean norm $$\|\alpha+\beta i+\gamma j+\delta k\|:= (\alpha^2+\beta^2+\gamma^2+\delta^2)^{\frac 12},$$ the quaternions become a real Banach algebra, and this norm actually satisfies the identity $$\| h_1 h_2 \| = \|h_1\| \ \|h_2\|\quad \hbox{ for all } h_1,h_2\in \mathbb{H}$$ (see \cite[Definition I.14.3]{BD}). There is a matricial identification of $\mathbb{H}$ in terms of $4\times 4$ matrices with real entries in which $\mathbb{H}$ embeds in $\mathcal{M}_4(\mathbb{R})$ as a real subalgebra via the assignment
$$\alpha+\beta i+\gamma j+\delta k \mapsto
\left( \begin{matrix}
\alpha & -\beta &-\gamma&-\delta \\
\beta & \alpha &-\delta &\gamma\\
\gamma&\delta&\alpha&-\beta\\
\delta&-\gamma&\beta&\alpha
\end{matrix}\right).$$

Each nonzero quaternion $h = a + b i + c j + d k$ has a unique inverse given by $$h^{-1} =\frac{ 1}{ a^2 + b^2 + c^2 + d^2} ( a - b i - c j - d k ).$$

A real or complex Banach algebra $\mathscr{A}$ is called a \emph{division algebra} if every nonzero element in $\mathscr{A}$ is invertible. By the celebrated Gelfand--Mazur theorem, each complex normed division algebra is isometrically isomorphic to $\mathbb{C}$ (see \cite[Theorem I.14.2]{BD}). The real setting is completely different. Clearly, $\mathbb{R}$ and $\mathbb{C}$ are real division Banach algebras, and as we have seen before, $\mathbb{H}$ also enjoys this property. This list exhausts all possibilities, because for each real normed division algebra $A$, there
exists an isomorphism $\Phi$ of $A$ onto $\mathbb{R}$, $\mathbb{C}$, or $\mathbb{H}$ such that $\|\Phi(x)\| = r(x)$ ($x\in A$) (see \cite[Theorem I.14.7]{BD}).\smallskip

Here we have another difference between real and complex Banach algebras. Any unital complex Banach algebra $\mathscr{A}$ is commutative if for some $\kappa>0$, the inequality $$\|a\|^2\leq \kappa\|a^2\|$$ holds for all $a\in \mathscr{A}$ (see \cite[Corollary II.16.8]{BD}). The same conclusion does not hold for real Banach algebras. For example, the real Banach algebra $\mathbb{H}$ actually satisfies $\|a^2\|=\|a\|^2$ for all $a\in \mathbb{H}$.\smallskip

{Let us return to the problem of determining whether a real Banach algebra admits an structure of complex Banach algebra for the same product and a subtle equivalent norm. A necessary condition on a unital real Banach algebra to admit a complex structure is to contain $\mathbb{C}\mathbf{1}$ in its center, and consequently its center must be at least two dimensional. Since the center of $\mathbb{H}$ is $\mathbb{R}\mathbf{1}$, we can immediately deduce that $\mathbb{H}$ does not admit a complex structure as an algebra.}


\subsection{Complexification of Banach $^*$-algebras}\ \smallskip

An \emph{algebra involution} on a real (respectively, complex) Banach algebra $\mathscr{A}$ is a real linear (respectively, conjugate-linear) mapping $^*: \mathscr{A}\to \mathscr{A}$ satisfying
\begin{enumerate}[$(a)$]\item $(a b)^* = b^* a^*$ for all $a,b\in \mathscr{A}$;
\item $(a^*)^* = a$ for all $a\in \mathscr{A}$.
\end{enumerate} The \emph{self-adjoint} or \emph{hermitian part} of $\mathscr{A}$ is the set $$\mathscr{A}_{sa} =\{ a\in \mathscr{A} : a^* = a\},$$ while the \emph{skew symmetric} part of $\mathscr{A}$ is defined as $$\mathscr{A}_{skew} =\{ a\in \mathscr{A} : a^* = - a\}.$$ The sets $\mathscr{A}_{sa}$ and $\mathscr{A}_{skew}$ are real subspaces of $A$. If $\mathscr{A}$ is a complex Banach algebra, then also $\mathscr{A}_{skew} = i \mathscr{A}_{sa}$. In any case, we have $$\mathscr{A}= \mathscr{A}_{sa} \oplus_{\mathbb{R}} \mathscr{A}_{skew}.$$ A real or complex Banach algebra equipped with an algebra involution is called a real or complex \emph{Banach $^*$-algebra}. In some references, like in \cite[Definition I.12.15]{BD}, a real or complex Banach $^*$-algebra is a real or complex Banach algebra $\mathscr{A}$ together with an algebra involution $^*$ satisfying $\|a^*\| = \|a\|$ for all $a\in \mathscr{A}$.\smallskip

A \emph{$^*$-homomorphism} (respectively, a \emph{$^*$-isomorphism}) between real or complex Banach $^*$-algebras $\mathscr{A}$ and $\mathscr{B}$ is a homomorphism (respectively, isomorphism) $\Phi : \mathscr{A}\to \mathscr{B}$ satisfying $\Phi (a^*) = \Phi(a)^*$ for all $a\in \mathscr{A}$. Conjugate-linear $^*$-homomorphisms and conjugate-linear $^*$-isomorphisms are similarly defined.\smallskip

If $\mathscr{A}$ is a complex Banach $^*$-algebra and $\tau: \mathscr{A}\to \mathscr{A}$ is an involution and a conjugate-linear $^*$-homomorphism, then the real form $\mathscr{A}^{\tau}$ is a real Banach $^*$-subalgebra of $\mathscr{A}$. In the other direction, by assuming that ${A}$ is a real $^*$-algebra, ${A}_c$ can be endowed with the involution $(a_1+ia_2)^*=a_1-ia_2$. Furthermore, if ${A}$ is a real Banach $^*$-algebra satisfying $\|a\|=\|a^*\|$ for all $a\in {A}$, then the standard extension of the involution $^*$ to $A_c$ also is an isometry with respect to the norm defined by $$\||a+i b |\|_{T,a} = \max\{\|L_{a+ i b}\|_{\mathcal{B}(A_c, \|\cdot\|_T)}, \|R_{a+ i b}\|_{\mathcal{B}(A_c, \|\cdot\|_T)}\}.$$ Namely, for the Taylor norm, we have $$\begin{aligned}\|(a + i b)^* \|_T &= \|a^* - i b^* \|_T:=\sup_{t\in[0,2\pi]}\|\cos(t) a^* + \sin(t) b^*\|\\
&=\sup_{t\in[0,2\pi]}\|(\cos(t) a + \sin(t) b)^*\| = \sup_{t\in[0,2\pi]} \|\cos(t) a + \sin(t) b\| = \|a +i b \|_T,
\end{aligned} $$ and it follows from this that $$\begin{aligned} \|L_{a^*- i b^*}\|_{\mathcal{B}(A_c, \|\cdot\|_T)} & = \sup_{\|(x+i y)\|_T\leq 1}\|(a^* - i b^*) (x+i y)\|_T \\
&= \sup_{\|(x+i y)\|_T\leq 1}\|(x+i y)^* (a + i b)\|_T = \|R_{a+ i b}\|_{\mathcal{B}(A_c, \|\cdot\|_T)},
\end{aligned}$$ and similarly $\|R_{a^*- i b^*}\|_{\mathcal{B}(A_c, \|\cdot\|_T)} = \|L_{a+ i b}\|_{\mathcal{B}(A_c, \|\cdot\|_T)}$, which implies that $$\||(a+i b)^* |\|_{T,a} = \||a^* -i b^* |\|_{T,a} = \||a+i b |\|_{T,a}.$$ Clearly, the Taylor norm on the complexification can be replaced by any of the norms on the complexification defined in Subsection \ref{subsec: standard complexification Banach algebras}. \smallskip

We arrive now to one of the starring models in mathematics, operator algebras, and a source of models for physics. A C$^*$-algebra is a complex Banach $^*$-algebra $\mathscr{A}$ satisfying the celebrated \emph{Gelfand--Naimark} axiom: \begin{equation}\label{eq GN axiom} \|a^* a \| = \| a\|^2\quad \hbox{ for all } a\in \mathscr{A}.
\end{equation}

Given a locally compact Hausdorff space $\Omega$ and a C$^*$-algebra $\mathscr{A}$, we write $C_b(\Omega, \mathscr{A})$ for the Banach space of all bounded continuous functions from $\Omega$ to $\mathscr{A}$ and by $C_0(\Omega,\mathscr{A})$ the closed subspace of all functions $f\in C_b(\Omega,\mathscr{A})$ such that $\|f\|$ vanishes at infinity equipped with the supremum norm $\|f\| = \sup_{t\in \Omega} \|f(t)\|$. When equipped with the pointwise sum, product, and involution, these spaces are C$^*$-algebras. The space $\mathcal{B}(\mathscr{H}),$ of all bounded linear operators on a complex Hilbert space $\mathscr{H},$ with the operator sum, product, and norm and with the adjoint operation as an involution, is a C$^*$-algebra, which is noncommutative when dim$(\mathscr{H}) > 1$. In the particular case in which $\mathscr{H}$ is $n$-dimensional for some natural $n$, the C$^*$-algebra $\mathcal{B}(\mathscr{H})$ naturally identifies with the algebra $\mathcal{M}_n(\mathbb{C})$ of (complex) $n \times n$ matrices. The subalgebra $K(\mathscr{H})$ of all compact linear operators on $\mathscr{H}$ is a C$^*$-subalgebra of $\mathcal{B}(\mathscr{H})$. A subalgebra of a C$^*$-algebra $\mathscr{A}$ is a subspace that is also closed for products. A subset $S$ of $\mathscr{A}$ is called self-adjoint if $a^*\in S$ for all $a\in S$. Actually, every norm closed self-adjoint subalgebra of some $\mathcal{B}(\mathscr{H})$ is a C$^*$-algebra because the Gelfand--Naimark axiom \eqref{eq GN axiom} is automatically inherited in this case.\smallskip

The celebrated Gelfand--Naimark theorem establishes that every C$^*$-algebra is isometrically $*$-isomorphic to a norm closed self-adjoint subalgebra of some $\mathcal{B}(\mathscr{H})$ (see \cite{GelNai43}, \cite[Theorem V.38.10]{BD}, and \cite[Theorem I.9.18]{Tak}). In the original result stated by Gelfand and Naimark \cite{GelNai43}, the definition of C$^*$-algebra included an extra axiom assuring that $1 + x^* x$ is invertible for all $x\in \mathscr{A}$. This extra axiom was shown to be superfluous by Fukamiya \cite{Fukamiya}, Kelley and Vaught \cite{KelleyVaught}, and Kaplansky \cite{Kap51}.\smallskip

There is a natural way to define a real analogue of a known algebraic-analytic structure by considering real forms under conjugations preserving some required algebraic-analytic structure. For example, if $\tau: \mathscr{A}\to \mathscr{A}$ is a conjugation on a C$^*$-algebra preserving the product, then the real form $\mathscr{A}^{\tau}$ is a norm closed self-adjoint real subalgebra of $\mathscr{A}$. According to this procedure, a real C$^*$-algebra $A$ is a real form of a C$^*$-algebra $\mathscr{A}$ under a conjugation $\tau$ which is also a $^*$-automorphism on $\mathscr{A}$, equivalently, a real Banach $^*$-algebra whose complexification $A_c$ admits a structure of C$^*$-algebra with a norm that extends the norm of $A$ and the involution of $A_c$ is the standard extension of the involution on $A$ (see \cite[Definition 5.1.1]{LI} or \cite{Goodearl, Inglestam64, IsRo, Palmer70}). This is one of the equivalent definitions of \emph{real C$^*$-algebras}; however, its handicap resides in the need of working with a superstructure of a C$^*$-algebra. In order to have an intrinsic definition not requiring an external structure, we recall the following result borrowed from the book of Li \cite{LI}.

\begin{theorem}\label{t characterizations of real Cstar}{\rm\cite[Corollary 5.2.11 and Proposition 7.3.4]{LI}} Let $A$ be a real Banach $^*$-algebra. Then the following statements are equivalent:\begin{enumerate}[$(1)$]\item $A$ is a real C$^*$-algebra;
\item $A$ is isometrically $^*$-isomorphic to a norm closed self-adjoint subalgebra of $\mathcal{B}(H)$ for some real Hilbert space $H$;
\item $A$ is hermitian {\rm(}i.e., ${\rm sp}(a)\subseteq \mathbb{R}$ for all $a\in A_{sa}${\rm)} and $\|a^* a\| = \|a\|^2$ for all $a\in A$;
\item $A$ is symmetric {\rm(}i.e., $a^*a\geq 0$ for all $a\in A${\rm)} and $\|a^* a\| = \|a\|^2$ for all $a\in A$;
\item $1 + a^* a$ is invertible in $A$ {\rm(}if $A$ is nonunital, then we consider its unitization{\rm)} and $\|a^* a\| = \|a\|^2$ for all $a\in A$;
\item The inequality $$\|a^*\| \ \|a\| \leq \| a^* a + b^* b\|$$ holds for all $a,b\in A$.
\end{enumerate}
\end{theorem}

The equivalence $(1)\Leftrightarrow (3)$ is due to Ingelstam \cite{Inglestam64}. It should be noted that the equivalent definition provided by statement $(5)$ is the notion employed by Goodearl \cite{Goodearl}, Chu et al. \cite{ChuDangRuVen}, and Isidro and Rodr{\'i}guez-Palacios \cite{IsRo}.\smallskip

Clearly, every C$^*$-algebra is a real C$^*$-algebra when it is regarded as a real Banach $^*$-algebra.\smallskip

Let us consider the C$^*$-algebra $\mathcal{M}_n(\mathbb{C})$ and a conjugation $\tau: \mathcal{M}_n(\mathbb{C})\to \mathcal{M}_n(\mathbb{C})$, $\tau((a_{ij}))= (\overline{a_{ij}})$. Clearly, $\tau$ is a conjugate-linear $^*$-automorphism on $\mathcal{M}_n(\mathbb{C})$ and the real form $\mathcal{M}_n(\mathbb{R}) = \mathcal{M}_n(\mathbb{C})^{\tau}$ is a real C$^*$-algebra whose algebra involution is just the transposition on $\mathcal{M}_n(\mathbb{R})$, the product is the matrix product, and the C$^*$-norm is the operator norm.\smallskip

Another interesting example of a real noncommutative unital C$^*$-algebra is the algebra of quaternions $\mathbb{H}$ described in page \pageref{def quaternions}. In this case, we consider the involution $^*$ on $\mathbb{H}$ defined by $(\alpha+\beta i+\gamma j+\delta k)^*=\alpha-\beta i-\gamma j-\delta k$. It is not hard to check that for each $a = \alpha+\beta i+\gamma j+\delta k\in \mathbb{H}$, we have
$$\|a^* a \| = \alpha^2 +\beta^2 +\gamma^2 +\delta^2= \|a\|^2,$$ and since $\mathbb{H}$ is a division algebra, we can deduce that it is a real C$^*$-algebra. \medskip

It should be added here that the mapping $$\alpha+\beta i+\gamma j+\delta k \mapsto
\left( \begin{matrix}
\alpha & -\beta &-\gamma&-\delta \\
\beta & \alpha &-\delta &\gamma\\
\gamma&\delta&\alpha&-\beta\\
\delta&-\gamma&\beta&\alpha
\end{matrix}\right)$$ is an isometric $^*$-monomorphism from $\mathbb{H}$ into $\mathcal{M}_4(\mathbb{R})= \mathcal{B}(\ell_2^4(\mathbb{R}))$ (see \cite[Example (2)]{LI}).\smallskip

In the real setting, the extra axiom ``$1 + a^* a$ is invertible in $A$'' does not follow from the other assumptions. For example, if we equip $\mathbb{C}$ (with its usual product and module) with the involution $\lambda^\star = \lambda$, the axiom $\| \lambda^\star \lambda \| = |\lambda^2| = |\lambda|^2$ holds for all $\lambda\in \mathbb{C}$. However, $1 + i^\star i = 0$ is not invertible in $\mathbb{C}$.\smallskip

We can now get back to the Gleason--Kahane--\.{Z}elazko and Kowalski--S{\l}odkowski theorems as a source of inspiration for new results on preservers. For example, by relaxing the hypothesis of linearity in the Gleason--Kahane--\.{Z}elazko, and replacing it by the preservation of products, Tour{\'e}, Schulz and Brits \cite[Problem 1.5]{TouSchuBrits2017} consider the converse of this result in the following preserver problem: Let $A$ be a complex unital Banach algebra, and suppose that $\phi: A\to \mathbb{C}$ is a continuous and multiplicative mapping satisfying $\phi(x)\in {\rm sp}(x)$ for all $x\in A$. Is $\phi$ automatically linear? \smallskip

Under the above conditions, assuming additionally that ${\rm sp}(x)$ is totally disconnected for each $x\in A$, then a multiplicative mapping $\phi: A\to \mathbb{C}$ with $\phi(x)\in {\rm sp}(x)$ for each $x\in A$, is linear if and only if it is continuous on $A$ (see \cite[Corollary 2.3]{TouSchuBrits2017}). Assuming that $A$ is a unital C$^*$-algebra and $\phi$ satisfies the commented assumptions, then there exists a (unique) character $\psi_{\phi}$ on $A$ satisfying $$\phi(e^{\lambda x}) = e^{\lambda \psi_{\phi} (x)},$$ for all $x\in A, \lambda\in \mathbb{C}$ \cite[Theorem 3.2 and Corollary 3.3]{TouSchuBrits2017}. The mapping $\phi$ and the induced character $\psi_{\phi}$ share many linear properties, and if $A$ is a von Neumann algebra or a commutative C$^*$-algebra, then it turns out that $\phi$ itself is linear, and that it coincides with its induced character (see \cite[Theorem 3.13]{TouSchuBrits2017}, and \cite[Theorem 2.5]{TouBrits2020}). The culminating point is the next result due to Brits, Mabrouk and Tour{\'e}:

\begin{theorem}\rm{\cite[Theorem 2.1]{BritMabrTou2021}} Let $A$ be a C$^*$-algebra, and let $\phi :A\to \mathbb{C}$ be a continuous multiplicative mapping such that $\phi(x) \in {\rm sp}(x)$ for all $x\in A$. Then $\phi$ is a character of $A$.
\end{theorem}

Another related problem, more in the line of preservers, reads as follows: Let $A$ be a complex and unital Banach algebra. Suppose $\phi:A\to \mathbb{C}$ is a mapping satisfying the following properties:\begin{enumerate}[$(P1)$]
\item $\phi(x) \phi(y) \in {\rm sp}(xy),$ for all $x, y\in A$;
\item $\phi$ is unital, i.e., $\phi(\textbf{1}) =1$;
\item $\phi$ is continuous on $A$.
\end{enumerate} Is $\phi$ a character?

\begin{theorem}\rm{\cite{TouBrits2020}} Let $A$ be a complex and unital Banach algebra, and let $\phi :A \to \mathbb{C}$ be a map satisfying the properties $(P1)$-$(P3)$ above. Then the following statements hold:
\begin{enumerate}[$(a)$]\item If $\sigma (x)$ is totally disconnected for each $x \in A$, then $\phi$ is a character of $A$;
\item If $A$ is a unital C$^*$-algebra, the mapping $$\psi_{\phi} (x):= \phi\left(\frac{x+x^*}{2}\right) + i \phi \left(\frac{x-x^*}{2 i}\right)$$ defines a character on $A$. Moreover, $\phi$ and $\psi_{\phi}$ coincide on the principal component of the invertible group of $A$.
\end{enumerate}
\end{theorem}
	
For additional results and generalizations in this line the reader can consult \cite{BurFerGarPe2015,BurFerGarPe2015tripl,Molnar19}. The preserver problems on real structures are on their very early stages.


\subsection{The Russo--Dye theorem}\ \smallskip

An element $u$ in a unital real or complex C$^*$-algebra is called \emph{unitary} if $u u^* = \11 = u^* u$. The set of all unitary elements in a unital real or complex C$^*$-algebra $A$ is actually a subgroup of $A$. The Russo--Dye theorem is a key result in the field of functional analysis, which is nowadays contained in most books and basic references. The theorem, whose statement appears below, was originally proved by Russo and Dye \cite{RuDye} in 1966. A surprisingly elementary proof was given by Gardner \cite{Gard84}.

\begin{theorem}\label{t RussoDye theorem}{\rm(Russo--Dye theorem, \cite{RuDye})} For each $($complex$)$ unital C$^*$-algebra $\mathscr{A}$, the closure of the convex hull of the unitary elements in $\mathscr{A}$ is the closed unit ball.
\end{theorem}

Actually, the open unit ball of each unital C$^*$-algebra $\mathscr{A}$ is inside the convex hull of the unitary elements in $\mathscr{A}$ (see \cite{Gard84}). Kadison and Pedersen \cite{KadPed85} sharpened this conclusion by showing that if $a\in \mathscr{A}$ with $\|a\| < 1-\frac{2}{n}$, then $a$ is the arithmetic mean of $n$ unitaries in $\mathscr{A}$. Haagerup \cite{Haa1990} proved that the same is true when $\|a\|\leq 1-\frac{2}{n}$, and a simplified version of the latter statement was given by Haagerup, Kadison, and Pedersen \cite{HaaKadPed2007}.\smallskip

The Russo--Dye theorem does not hold for unital real C$^*$-algebra. For example, $A=C([0,1],\mathbb{R})$ clearly is a unital real C$^*$-algebra with the identity as involution. The set of unitaries in this real C$^*$-algebra $A$ is so small that reduces to $\pm \11$. Thus the convex hull of the unitaries in $A$ is far from covering the whole closed unit ball.\smallskip

There are subtle differences between the real and complex structures. For example, it follows from the local theory of C$^*$-algebras that every hermitian or self-adjoint element in the closed unit ball of a unital C$^*$-algebra $\mathscr{A}$ can be written as the average of a unitary $u$ and its transposed $u^*$. Consequently, then every $x\in \mathscr{A}$ is a linear combination of four unitary elements (see \cite[Proposition I.4.9]{Tak}). This conclusion is not true in the real setting; however, each skew-hermitian element in the open unit ball of a unital real C$^*$-algebra $A$ is the average of two unitaries in $A$ (see \cite[Lemma 3.1.3]{LI}).\smallskip

To explore the Russo--Dye theorem in the setting of unital real C$^*$-algebras, we refresh some well-known results on local theory and continuous functional calculus. Let $A$ be a real C$^*$-algebra. Let us recall that $\sigma (a) =\sigma_{A} (a) = \overline{\sigma_{A} (a)}$ for each $a\in A$ (see \eqref{eq the spectrum in the real setting is self adjoint}). If $a \in A$ is a normal element (i.e., $a^*a = aa^*$), then the real C$^*$-subalgebra $C^*(a)$ of $A$ generated by $a$ and $\11$ is isometrically $*$-isomorphic to $$C({\rm sp}(a), \overline{\ \cdot \ }) =\{ f\in C({\rm sp}(a)) : f(\overline{\lambda}) = \overline{f(\lambda)},\hbox{ for all } \lambda\in {\rm sp}(a)\},$$ and under this identification, the element $a$ corresponds to the identity mapping on ${\rm sp}(a)$ (see \cite[Proposition 5.1.6$(2)$]{LI}).\smallskip

In particular, for each self-adjoint element $h\in A$, the real C$^*$-subalgebra $C^*(h)$ is identified with $C({\rm sp}(h), \mathbb{R})$. Therefore, for each continuous function $f:{\rm sp}(a)\to \mathbb{R}$, there exists a unique element $f(h)\in C^*(h)$ that is identified with the function $f$ under this $^*$-isomorphism. The element $f(h)$ is called the \emph{continuous functional calculus of $f$ at the element $h$}.\smallskip

A close statement to the Russo--Dye theorem in the real setting was given by Li \cite{Li75,LI} between 1975 and 1979. The concrete result can be stated in the following terms.

\begin{theorem}\label{t RD real I}\cite[Theorem 7.2.4 and Proposition 7.2.7]{LI} Let $A$ be a unital real C$^*$-algebra and let $\mathcal{B}_{A}$ denote the closed unit ball of $A$. Then the convex hull of the subset $$\{\cos(b) e^a \ | \ a,b \in A, a^* = -a, b^* = b\}$$
is dense in $\mathcal{B}_{A}$. Furthermore, $$ \hbox{int}( \mathcal{B}_{A}) \subseteq co \{\cos(b) e^a \ | \ a,b \in A, a^* = -a, b^* = b\}\subseteq \mathcal{B}_{A},$$ where int$(\mathcal{B}_A)$ denotes the open unit ball of $A$.
\end{theorem}

After considering the counterexample to the Russo--Dye theorem for real C$^*$-algebras, it seems natural to ask whether the original statement in this theorem holds under stronger hypotheses on the C$^*$-algebra, for example, when we have a real von Neumann algebra. First, we recall the definition of real and complex von Neumann algebras.\smallskip

Let $\mathscr{H}$ be a real or complex Hilbert space. Following the standard notation, for each subset $M$ of $\mathcal{B}(\mathscr{H})$, we write $M^{\prime}$ for the set
of all bounded operators on $\mathscr{H}$ commuting with every operator in $M$. The set $M^{\prime}$ is a Banach algebra of operators containing the identity operator $\11$. If $M$ is self-adjoint (i.e., $x^*\in M$ for all $x\in M$), then $M^{\prime}$ is a real or complex C$^*$-algebra acting on $\mathscr{H}$, which is closed with respect to all the standard locally convex topologies on $\mathcal{B}(\mathscr{H})$ (see \cite[\S II.2]{Tak} and \cite[\S 4.2]{LI}). A \emph{von Neumann algebra} (respectively, a \emph{real von Neumann algebra} is a $^*$-subalgebra $\mathscr{M}$ of $\mathcal{B}(\mathscr{H})$ for some complex (respectively, real) Hilbert space $\mathscr{H}$ whose bicommutant coincides with $\mathscr{M}$ itself, that is, $\mathscr{M}^{\prime\prime} = \mathscr{M}$ (see \cite[Definition II.3.2]{Tak} and \cite[Definition 4.3.1]{LI}). The famous \emph{von Neumann's bicommutant theorem} shows that for each complex (respectively, real) Hilbert space $\mathscr{H}$ and each complex (respectively, real) C$^*$-subalgebra $\mathscr{M}$ of $\mathcal{B}(\mathscr{H})$ containing the identity, the following conditions are equivalent:\begin{enumerate}[$(a)$]\item $\mathscr{M} = \mathscr{M}^{\prime\prime};$
\item $\mathscr{M}$ is weakly closed;
\item $\mathscr{M}$ is strongly closed.
\end{enumerate}  (see \cite[Theorem 2.2.2]{Ped} or \cite[Theorem II.3.9]{Tak} in the complex setting and \cite[Proposition 4.3.2]{LI} in the case of real von Neumann algebras, see also \cite[Theorems 8.1.30 and 8.1.31 and Proposition 8.1.106]{Cabrera-Rodriguez-vol2}). A celebrated theorem due to Sakai asserts that a C$^*$-algebra $\mathscr{A}$ is isometrically $*$-isomorphic to a von Neumann algebra (these algebras are frequently termed W$^*$-algebras) if and only if it is a dual Banach space, and in such a case, it admits a unique isometric predual and its product is separately weak$^*$ continuous (see \cite[Theorem III.3.5 and Corollary II.3.9]{Tak} or \cite[Theorem 3.9.8]{Ped}), which is also equivalent to say that $\mathscr{A}$ is monotone closed and admits sufficiently many normal positive linear functionals (see \cite[Theorem III.3.16]{Tak}).\smallskip

Similarly, a real W$^*$-algebra $M$ is defined as a real C$^*$-algebra whose complexification is a W$^*$-algebra, in such a case, there exists a weak$^*$-continuous conjugate-linear $^*$-homomorphism of period-2 on $M_c$ such that $M = (M_c)^{\tau}$, and defining $$\tau_*: (M_c)_*\to (M_c)_*,$$ $$\tau_* (\varphi) (x) := \overline{\varphi(\tau(x))} \ \ (x\in M_c, \varphi\in (M_c)_*),$$ we get another conjugation on $M_c$ satisfying that $M$ is the dual space of the real form $M_* := ((M_c)_*)^{\tau_*}$, and the product of $M$ is separately weak$^*$ continuous (see \cite[Proposition 6.1.3]{LI}). Actually, a real C$^*$-algebra $M$ is a real W$^*$-algebra if and only if it is a dual Banach space, and in such a case, its product is separately weak$^*$-continuous (see \cite[Theorem 6.1.7]{LI} and \cite[Theorem 1.11]{IsRo}). For these reasons, we shall not distinguish between real von Neumann algebras and real W$^*$-algebras.\smallskip

{Thanks to the Gelfand-Naimark axiom, we can now establish a strengthened version of Theorem \ref{t Ingelstam  characterization of complex type} and Proposition \ref{p original real norm is a complex norm}.
	
\begin{theorem}\label{t Inglestam real Cstar algebra} Let $A$ be a real C$^*$-algebra. Then $A$ is a complex C$^*$-algebra regarded as a real C$^*$-algebra if and only if there exists a linear operator $J$ on $A$ satisfying: \begin{enumerate}[$(a)$] \item $J$ is an $A$-module homomorphism, that is, $$J(ab) =J(a) b = a J(b), \hbox{ for all } a,b\in A;$$
			\item $-J^2$ is the identity map on $A$.
\end{enumerate} 	
\noindent Consequently, if $A$ is unital, then $A$ is a complex C$^*$-algebra regarded as a real C$^*$-algebra if and only if there exists an element $\boldsymbol{\iota}$ in the center of $A$, satisfying $\boldsymbol{\iota}^2 = -\mathbf{1}$.
\end{theorem}
	
\begin{proof} The necessary condition is clear by taking $J(a) = i a$ ($a\in A$). For the sufficient implication we shall simply prove that $\|\cdot\|$ is a complex norm for the product by complex scalars defined by $(\alpha + i \beta ) a = \alpha + \beta J(a)$.\smallskip
	
We begin by observing that the mapping $J$ is continuous. Indeed, every real C$^*$-algebra admits an approximate unit \cite[Proposition 5.2.4]{LI}, and thus the continuity of $J$ follows from \cite[Proposition 2.1]{Inglestam66}.\smallskip

We shall next prove that \begin{equation}\label{eq J(a*) = -J(a)*} J(a)^* = - J(a^*), \hbox{ for all } a \in A.
\end{equation}

To get the desired conclusion we first observe that for a central element $z$ in a unital C$^*$-algebra $A$ the condition \begin{equation}\label{eq zsquare = -1 implies z anti} \hbox{$z^2 = -\mathbf{1}$ implies $z^* = -z$.}
\end{equation} Namely, let us write $z = h + k$ with $h^* = h$ and $k^*= -k$. Clearly $h$ and $k$ are central elements too. The identity $z^2 = -\mathbf{1}$ is equivalent to $ h^2 + k^2 + 2 h k = -\mathbf{1},$ and gives $h k =0$ because the latter is anti-symmetric. Now the equality $h^2 + k^2 = -\mathbf{1}$ implies 
$$h h^* h = h^3 = h^3 + k^2 h = -h,$$ and thus $0\leq h^4 = -h^2\leq 0$, which proves that $h= 0$, as desired.\smallskip

Since $J$ is continuous, the bitransposed mapping $J^{**} : A^{**}\to A^{**}$ is weak$^*$ continuous. Combining this continuity with the separate weak$^*$ continuity of the product of the real von Neumann algebra $A^{**}$, Goldstine's theorem and the hypotheses on $J$, it is not hard to check that $J^{**}$ satisfies $-(J^{**})^2 = Id_{A}$ and $$J^{**} (a b) = J^{**} (a) b = a J^{**} (b), \hbox{ for all } a,b\in A^{**}.$$ We shall prove \eqref{eq J(a*) = -J(a)*} for $J^{**}$. The element $J^{**} (\mathbf{1})$ (which is central by the properties of $J^{**}$) satisfies $$J^{**} (\mathbf{1})^2 =  J^{**} (\mathbf{1}) J^{**} (\mathbf{1}) =\mathbf{1} (J^{**})^2 (\mathbf{1}) = - \mathbf{1}.$$ It follows from \eqref{eq zsquare = -1 implies z anti} that $J^{**}(\mathbf{1})^* = - J^{**} (\mathbf{1})$. Now, the properties of $J^{**}$ lead to $$J^{**} (a) = J^{**} (\mathbf{1} a) = J^{**} (\mathbf{1}) a, \hbox{ for all } a\in A^{**}.$$ Therefore $$J^{**} (a)^{*} = a^* J^{**} (\mathbf{1})^* = - J^{**} (\mathbf{1}) a^* = - J^{**} (a^*), \hbox{ for all } a\in A^{**},$$ which proves \eqref{eq J(a*) = -J(a)*} for $J^{**}$.\smallskip
		
Finally, by applying \eqref{eq J(a*) = -J(a)*} we get
$$\begin{aligned} \| (\alpha + i \beta ) a\|^2 & = \left\| ((\alpha + i \beta ) a) ((\alpha + i \beta ) a)^*\right\| = \left\| (\alpha a + \beta J(a)) (\alpha a + \beta J(a))^*\right\| \\
&= \left\| (\alpha a + \beta J(a)) (\alpha a^* -  \beta J(a^*))\right\| \\
&= \left\| \alpha^2 a a^* - \beta^2 J(a)J(a^*) + \alpha \beta J(a) a^* -  \alpha \beta a J(a^*) \right\|\\
&= \left\| \alpha^2 a a^* - \beta^2 a J^2(a^*)  \right\| = \left\| \alpha^2 a a^* + \beta^2 a a^*  \right\|\\
&=  (\alpha^2 + \beta^2) \|aa^*\| =  |\alpha + i \beta|^2 \ \|a\|^2.
\end{aligned}$$	
\end{proof}

The second conclusion in the above theorem appears in \cite[Exercise (13A)]{Goodearl}.\smallskip
}

The question of whether the original statement in the Russo--Dye theorem is valid for real von Neumann algebras was open for several years. For example, this problem was explicitly posed by Becerra Guerrero et al. \cite[p. 98 and Problems (P1) and (P3)]{BeBuKaRod08} in the particular case of $\mathcal{B}(H)$ for a real Hilbert space $H$ in 2008. The question was addressed in 2012 by Navarro-Pascual and Navarro-Pascual (see \cite[Theorem 5 and Corollary 6]{NavNav2012}), who gave an explicit and positive answer to this question by proving that actually the Russo--Dye theorem holds true for any real von Neumann algebra. However, as observed by Mori and Ozawa \cite[Proof of Corollary 3]{MoriOza2020}, the desired solution can be easily obtained from the results established by Li \cite{Li75,LI}, which have been reviewed in Theorem \ref{t RD real I}. Namely, let $b$ be a self-adjoint element in a real von Neumann algebra $M$. By \cite[Proposition 5.1.6$(2)$ and Theorem 6.3.3]{LI}, the real von Neumann subalgebra $M_b$ of $M$ generated by $b$ and the unit element $\11$ is isometrically $^*$-isomorphic to $C(\Omega, \mathbb{R})$, where $\Omega$ is a hyperstonean compact Hausdorff space and $b$ is a positive generator (this can also be obtained from \cite[Lemma 4.1.11]{HOS}). We recall that for a compact Hausdorff space $K$, the closed unit ball of $C(K,\mathbb{R})$ coincides with the closed convex hull of its extreme points (i.e., the unitary elements in $C(K,\mathbb{R})$) if and only if $K$ is totally disconnected (see, for example, \cite[p. 191]{Da} or \cite{Phelps65}). We recall that a Hausdorff space is said to be \emph{extremally disconnected} if the closure of every open subset is open. A compact extremely disconnected space is called a \emph{stonean} space. It is known that every extremally disconnected space is totally disconnected and that every hyperstonean compact Hausdorff space is a stonean compact Hausdorff space satisfying additional properties (cf. \cite[Definition III.1.14]{Tak}). Therefore, the elements in the closed unit ball of $M_b$ can be approximated in the norm by convex combinations of unitaries in $M_b$ and hence in $M$. Since, by the definition of the continuous functional calculus, $\cos(b)$ lies in the real C$^*$-subalgebra of $M$ generated by $b$, it follows that $\cos(b)\in M_b,$ and thus it can be approximated in norm by convex combinations of unitaries in $M$. Having in mind that the set of unitaries in $M$ is a subgroup, elements of the form $\cos(b) e^a,$ with $a,b\in M$, $a=-a^*$, and $b=b^*$, can be approximated in norm by convex combinations of unitaries in $M$. Theorem \ref{t RD real I} implies that the convex hull of the unitaries in $M$ is norm dense in its closed unit ball.\smallskip

Let us finally observe that the results by Kadison and Pedersen \cite{KadPed85}, Haagerup \cite{Haa1990}, and Haagerup, Kadison, and Pedersen \cite{HaaKadPed2007} on means and convex combinations of unitaries in unital C$^*$-algebras have not been fully explored in the setting of unital real C$^*$-algebras nor real von Neumann algebras.

\subsection{Kadison--Schwarz inequality}\label{subsec: KadisonSchwarz inequ}\ \smallskip

It is well known that an element $a$ in a real or complex C$^*$-algebra $\mathscr{A}$ is called \emph{positive},
denoted by $a \geq 0$, if $a^* = a$ and $\sigma (a) \subseteq \mathbb{R}^+_0$ (see \cite[Definition I.6.2]{Tak} or \cite[Lemma 1.3.1]{Ped} and \cite[Definition 5.2.1 and Proposition 5.2.2]{LI} for the basic properties in the complex and real settings, respectively). The closed cone of positive elements in $\mathscr{A}$ is denoted by $\mathscr{A}^+$. Suppose now that $A$ is a real C$^*$-algebra. Let $A_c$ denote the complexification of $A$ and let $\tau : A_c\to A_c$ be the involutive conjugate-linear $^*$-automorphism satisfying $A = (A_c)^{\tau}$. It is known that $$A^+ = A\cap A_c^+ = (A_c)^{\tau} = \{ b^* b : b\in A\}$$ (see \cite[Proposition 5.2.2]{LI}). These are the usual properties of the cone of positive elements in a complex C$^*$-algebra (see \cite[\S 1.4 and Theorem 1.4.4]{S}).\smallskip

Given a (complex) C$^*$-algebra $\mathscr{A},$ a linear mapping $\varphi : \mathscr{A}\to \mathbb{C}$ is called \emph{positive} if $\varphi( \mathscr{A}^+) \subseteq \mathbb{R}_0^+.$ Each positive linear functional $\varphi$ satisfies the \emph{Cauchy--Schwarz inequality} (see \cite[Theorem 3.1.3]{Ped} and \cite{MOS1, MOS2}):
\begin{equation}\label{eq CS inequality complex Cstar} |\varphi (b^*a)|^2 \leq \varphi (a^* a) \varphi(b^* b)\quad \hbox{ for all } a, b \in \mathscr{A}.
\end{equation} It is known that every positive linear functional $\varphi$ on a complex C$^*$-algebra is continuous and satisfies $\|\varphi\| = \varphi (\11)$ if $\mathscr{A}$ is unital (see \cite[Lemma I.9.9]{Tak}) and that $\|\varphi\| = \lim_{i} \varphi(u_i)$ for some/any approximate unit $\{u_i\}_i$ in $\mathscr{A}$ (see \cite[Proposition 3.1.4]{Ped}). Actually a bounded linear functional $\varphi\in \mathscr{A}^*$ is positive if and only if the condition $\|\varphi\| = \lim_{i} \varphi(u_i)$ holds for some/any approximate unit $\{u_i\}_i$ in $\mathscr{A}$ if and only if $\varphi$ attains its norm at a positive element in $\mathscr{A}$ (see \cite[Proposition 1.5.2]{S}). \smallskip

The notion of positive linear functional changes in the real setting. For example, when $\mathbb{C}$ is regarded as a real C$^*$-algebra, the linear functional $\varphi (a+ i b)= \frac{1}{\sqrt{2}} (a + b)$ maps positive elements to positive elements; however, $\varphi (\11) = \frac{1}{\sqrt{2}} < \|\varphi\| =1$. This functional $\varphi$ does not satisfy the Cauchy--Schwarz inequality in \eqref{eq CS inequality complex Cstar}. Indeed, for $\lambda,$ $\mu\in \mathbb{C}$, we have \begin{equation}\label{eq CS ineq is not true in the real case functional} \frac12 \left( \Re\hbox{e} (\lambda \overline{\mu}) + \Im\hbox{m} (\lambda \overline{\mu}) \right)^2 = |\varphi (\lambda \overline{\mu})|^2 \nleqq \varphi (|\lambda|^2) \varphi(|\mu|^2) = \frac12 |\lambda|^2 \ |\mu|^2 = \frac12 |\lambda \ \overline{\mu}|^2 .
\end{equation}

Let $A$ be a real C$^*$-algebra. A linear mapping $\varphi : A\to \mathbb{R}$ is called \emph{positive} if $\varphi( {A}^+) \subseteq \mathbb{R}_0^+$ and $\varphi|_{A_{skew}} \equiv 0$ (see \cite[Definition 5.2.5]{LI}). Every positive linear functional on a real C$^*$-algebra satisfies the real version of the Cauchy--Schwarz inequality in \eqref{eq CS inequality complex Cstar} (see \cite[Proposition 5.2.6$(1)$]{LI}). As in the complex setting, a bounded linear functional $\varphi$ in the dual of a real C$^*$-algebra is positive if and only if it attains its norm at a positive element (see \cite[Proposition 5.2.6$(3)$]{LI}).\smallskip

Let $A$ be a real C$^*$-algebra with complexification $A_c$, and let $\tau: A_c\to A_c$ be a conjugate-linear $^*$-automorphism such that $A$ identifies with the real form $A_c^{\tau}$. Similar arguments to those employed in the previous subsections allow us to conclude that the mapping $\tau^{\sharp} : A_c^*\to A_c^*$, $\tau^{\sharp} (\varphi) (a) = \overline{\varphi(\tau(a))}$, is a conjugation on $A^*_c$ and the corresponding real form $(A_c^*)^{\tau^{\sharp}} = \{ \varphi\in A_c^* :\tau^{\sharp} (\varphi) = \varphi\}$ identifies with $A^*$ via the following linear isometric surjection:
$$ (A_c^*)^{\tau^{\sharp}} \to A^*, \quad \varphi\mapsto \Re\hbox{e} \varphi|_{A}= \varphi|_{A}.$$ If we write $\mathscr{A}^*_{+}$ for the set of all positive linear functionals on a real or complex C$^*$-algebra $\mathscr{A}$, then, in the case where $A$ is a real C$^*$-algebra, we can actually deduce (see \cite[Proposition 5.2.6$(2)$ and $(4)$]{LI}) that
$$ ((A_c^*)_+)^{\tau^{\sharp}}\equiv A^*_{+}, \quad \varphi\mapsto \Re\hbox{e} \varphi|_{A}= \varphi|_{A}.$$ 

Note that with $b=\textbf{1}$ and $a=a^*$ in \eqref{eq CS inequality complex Cstar}, we get
\[\varphi(a)^2=|\varphi(a)|^2 \leq \varphi(\textbf{1}^*\textbf{1}) \varphi(a^*a)=\varphi(a^2).\]
An interesting questions asks what happens if we assume that $\varphi$ takes its values in a general C$^*$-algebra instead of ${\mathbb C}$.\smallskip

In 1952, Kadison \cite{Kad52} proved a generalized Schwarz inequality for positive linear mappings between C$^*$-algebras ---now called a generalized \emph{Kadison--Schwarz inequality} for C$^*$-algebras. It asserts that if $\Phi: \mathscr{A} \to \mathscr{B}$ is a unital positive linear map and $a\in \mathscr{A}$ is self-adjoint, then $$\Phi(a)^2 \leq \Phi(a^2).$$

Concrete version of the Kadison--Schwarz inequality states that for each positive bounded linear mapping $\Phi: \mathscr{A} \to \mathscr{B}$ between two C$^*$-algebras, the inequality
\begin{equation}\label{eq GKS inequality} \Phi (a) \circ \Phi (a)^{*} \leq \| \Phi \| \ \Phi (a \circ
a^{*}),
\end{equation}holds for all $a\in \mathscr{A}$, where $\circ$ denotes the natural Jordan product given by
$x\circ y = \frac{1}{2} (x y + y x)$ (see \cite[Theorem 1]{Kad52}).\smallskip

Another basic result in the theory of C$^*$-algebras asserts that for each self-adjoint element $a$ in a C$^*$-algebra $\mathscr{A}$, there exists a norm-one positive linear functional $\varphi\in \mathscr{A}$ such that $\|a\| = |\varphi(a)|$ \cite[Proposition 1.5.4]{S}. It is further known that an element $a\in \mathscr{A}$ is positive (respectively, self-adjoint) if $\varphi (a)\geq 0$ (respectively, $\varphi (a)\in \mathbb{R}$) for all positive linear functionals $\varphi\in A^*$ (see \cite[Theorem 4.3.4]{KR1}). By combining the Russo--Dye theorem with the Kadison--Schwarz inequality and the characterization of positive linear functionals, it can be seen that a linear mapping $\Phi$ between unital C$^*$-algebras $\mathscr{A}$ and $\mathscr{B}$ with $\Phi(\11) = \11$ is positive if and only if $\Phi$ is continuous with $\|\Phi\| = 1$ (see \cite[Corollary 3.2.6]{BraRo}).\smallskip

Let $A$ and $B$ be two real C$^*$-algebras. A linear mapping $\Phi : A\to B$ sending positive elements to positive elements (i.e., $\Phi (A^+)\subseteq B^+$) need not satisfy the generalized Kadison--Schwarz inequality \eqref{eq GKS inequality} nor its consequences. We have already seen that a linear functional $\varphi: \mathbb{C}_{_\mathbb{R}}\to \mathbb{R}$ mapping positive elements to positive elements that does not satisfy the Cauchy--Schwarz inequality (see \eqref{eq CS ineq is not true in the real case functional}). Additional counterexamples to the consequences of the Kadison--Schwarz inequality can be given as follows.

\begin{example}\label{example first March} The real linear mapping $\Phi_1: \mathbb{C}\to \mathbb{C}$, $\Phi_1(a + i b) = a + i 3 b$ is clearly unital ($\Phi_1 (\11) = \11$), maps positive elements to positive elements, and $\|\Phi_1\| = 3$.\smallskip

Consider next the linear mapping $\Phi_2 : \mathcal{M}_{2}
(\mathbb{R}) \to \mathcal{M}_{2} (\mathbb{R})$ defined by
$$\Phi_2 \begin{pmatrix}
a_{11} & a_{12} \\
a_{21} & a_{22}
\end{pmatrix} = \begin{pmatrix}
a_{11} & a_{12} \\
- a_{12} & a_{11}
\end{pmatrix}.$$ Clearly, $\Phi_2$ is unital, but it does not map positive elements to positive elements because $\Phi_2 \left(
\begin{array}{cc}
2 & 1 \\
1 & 2 \\
\end{array}
\right) = \left(
\begin{array}{cc}
2 & 1 \\
-1 & 2 \\
\end{array}
\right)$, which is skew symmetric. The element $p = \left(
\begin{array}{cc}
1 & 0 \\
0 & 0 \\
\end{array}
\right)$ is a projection (i.e., a symmetric idempotent) in $\mathcal{M}_2 (\mathbb{R})$, and hence the mapping $$\left(
\begin{array}{cc}
a_{11} & a_{12} \\
a_{21} & a_{22} \\
\end{array}
\right) \mapsto p \left(
\begin{array}{cc}
a_{11} & a_{12} \\
a_{21} & a_{22} \\
\end{array}
\right) = \left(
\begin{array}{cc}
a_{11} & a_{12} \\
0 & 0 \\
\end{array}
\right)$$ is linear and contractive. It is well known that $\left\| \left(
\begin{array}{cc}
a_{11} & a_{12} \\
0 & 0 \\
\end{array}
\right) \right\| = \sqrt{|a_{11}|^2 + |a_{12}|^2 }$, even in $\mathcal{M}_2(\mathbb{C})$. On the other hand, the matrices $\11 = \left(
\begin{array}{cc}
1 & 0 \\
0 & 1 \\
\end{array}
\right)$ and
$\widehat{i} = \left(
\begin{array}{cc}
0 & 1 \\
-1 & 0 \\
\end{array}
\right)$ have a very special behavior in $\mathcal{M}_2(\mathbb{R})$ because they generate an isometric copy of $\mathbb{C}$ ---a conclusion, which is no longer true for $\mathcal{M}_2(\mathbb{C})$---. We therefore infer that $$\left\| \Phi_2 \begin{pmatrix}
a_{11} & a_{12} \\
a_{21} & a_{22}
\end{pmatrix} \right\| = \left\| \begin{pmatrix}
a_{11} & a_{12} \\
- a_{12} & a_{11}
\end{pmatrix} \right\| = \left\| a_{11} \ \11 + a_{12} \ \widehat{i} \ \right\| = \sqrt{|a_{11}|^2 + |a_{12}|^2 },$$ witnessing that $\Phi_2$ is a contractive mapping, and hence $\|\Phi_2\| = 1$.

\end{example}

The previous counterexamples can also be employed to deduce that the natural extension of a bounded linear mapping between two real C$^*$-algebras to the corresponding complexifications need not be, in general, norm preserving. Actually, the extension of a linear mapping preserving positive elements between real C$^*$-algebras to the corresponding complexifications need not send positive elements to positive elements.

\begin{example}\label{example second March} For the mappings ${\Phi}_1$ and ${\Phi}_2$ in Example \ref{example first March}, the mapping $\widehat{\Phi}_1: \left(\mathbb{C}_{_\mathbb{R}}\right)_{c}=\mathbb{C}_{_\mathbb{R}} \oplus i \mathbb{C}_{_\mathbb{R}} \to \left(\mathbb{C}_{_\mathbb{R}}\right)_{c}$, defined by $$\widehat{\Phi}_1(\lambda + i \mu) = \Phi_1 (\lambda) + i \Phi_1(\mu)$$ is clearly bounded complex linear and unital. In this setting, $\widehat{\Phi}_1$ maps positive elements to positive elements if and only if $\|\widehat{\Phi}_1\|= 1$. To simplify the notation, let us write $\mathscr{X}$ for the complex C$^*$-algebra $\left(\mathbb{C}_{_\mathbb{R}}\right)_{c}$. It is not hard to see that $\mathscr{X}_{sa} = \{ \lambda + i \mu : \lambda\in \mathbb{R}, \ \mu\in i \mathbb{R} \},$ that $\mathscr{X}^{+} = \{ \lambda^2 - \mu^2 + 2 i \lambda \mu : \lambda\in \mathbb{R}, \ \mu\in i \mathbb{R} \}$, and that the element $x = (2^2 + 1^2) + i (4 i)$ is positive in $\mathscr{X}$, while $\widehat{\Phi}_1 (5 + i (4 i)) = 5 + i (12 i)\notin \mathscr{X}^+$. Since $1< \|\Phi_1\|\leq \|\widehat{\Phi}_1\|$, we also deduce that $\widehat{\Phi}_1$ is noncontractive.\smallskip

The counterexample given by the mapping $\Phi_2 : \mathcal{M}_{2}(\mathbb{R}) \to \mathcal{M}_{2} (\mathbb{R})$ above admits a nice algebraic-geometric reinterpretation. Let $\mathcal{C}$ denote the real C$^*$-subalgebra of $\mathcal{M}_2(\mathbb{R})$ generated by $\11$ and let $s = \left(
\begin{array}{cc}
0 & 1 \\
1 & 0 \\
\end{array}
\right). $
Since $s^* =s$ and $s^2 = 1$, it is easy to see that $\mathcal{C} = \mathbb{R} \11 \oplus \mathbb{R} s$ is a commutative unital real C$^*$-algebra whose involution is the identity. Furthermore, since for $a,b\in \mathbb{R}$, the eigenvalues of the matrix $ a \11 + b s= \left(
\begin{array}{cc}
a & b \\
b & a \\
\end{array}
\right)$
are $\{a\pm |b|\}$, it can be easily seen that $\|a \11 + b s\| = |a| + |b|,$ which gives a concrete expression of the C$^*$-norm on $\mathcal{C}$ ---we observe that $\mathcal{C}$ is an example of a two-dimensional real spin factor studied by Kaup \cite[\S 4]{Ka97}, and we will find them again in subsequent subsections---. The restriction $\Psi= \Phi_2|_{\mathcal{C}} : \mathcal{C} \to \mathcal{M}_{2} (\mathbb{R})$ is a real linear bijection from $\mathcal{C}$ onto $\mathbb{C}_{\mathbb{R}},$ when the latter is regarded as the real C$^*$-subalgebra of $\mathcal{M}_2 (\mathbb{R})$ generated by $\11$ and $\widehat{i}$. Both real C$^*$-algebras $\mathcal{C}$ and $\mathbb{C}_{\mathbb{R}}$ are commutative. It follows from what we have seen before that $\Psi$ is unital with $\|\Psi\| =1$. The canonical complex linear extension $\widehat{\Psi} : \mathcal{C}_c = \mathcal{C} + i \mathcal{C} \to \left(\mathbb{C}_{_\mathbb{R}}\right)_{c}=\mathbb{C}_{_\mathbb{R}} \oplus i \mathbb{C}_{_\mathbb{R}} =\mathscr{X}$ is unital, but it cannot be contractive nor positive since $\Psi$ does not map positive elements to positive elements. Furthermore, $$\|\11 + i s\|_{\mathcal{C}_c}^2 = \| (\11 + i s)^* (\11 + i s) \|_{\mathcal{C}_c}^2= \| 2 \11 \|_{\mathcal{C}_c}^2 = 2,$$ while $\widehat{\Psi} (\11 + i s) = \11 + i (i)= y \in \mathscr{X}_{sa}$ with $y^2 = 2 y$, and thus $\| y\|_{\mathscr{X}}^2 = \|y^2\|_{\mathscr{X}} = \|2 y \|_{\mathscr{X}} = 2 \|y\|_{\mathscr{X}}$. Therefore, $\|\widehat{\Psi} (\11 + i s) \|_{\mathscr{X}} = \|y\|_{\mathscr{X}} = 2 \nleq \|\11 + i s\|_{\mathcal{C}_c} = \sqrt{2}$.\smallskip

Let us finally note that $\Psi^{-1} : \mathbb{C}_{\mathbb{R}} \to \mathcal{C}$ maps positive elements to positive elements.
\end{example}

Due to the previous counterexamples, the theory of completely bounded and completely positive linear maps gains prominence in the real setting. Let us simply recall the basic notions. For each real or complex C$^*$-algebra $\mathscr{A}$ and each natural number $n$, there exists a unique real or complex C$^*$-norm on the space $\mathcal{M}_n(\mathscr{A}),$ of all $n\times n$-matrices $a = (a_{ij})$ with entries $a_{ij}$ in $\mathscr{A}$, with respect to the natural linear space structure, matrix multiplication, and algebra involution given by $a^* = (a_{ij})^* = (a_{ji}^*)$, making $\mathcal{M}_n(\mathscr{A})$ a real or complex C$^*$-algebra (see \cite[\S IV.3]{Tak}, \cite[\S 2]{ChuDangRuVen}, and \cite[Proposition 5.1.10]{LI} for references in the complex and real case, respectively). This can be done via the standard procedure for operator spaces (see \cite{EffRuanBook, Paulsen, PisBook} and \cite[IV.3]{Tak}). If $A$ is a real C$^*$-algebra represented as a real self-adjoint closed subalgebra of some $\mathcal{B}(H)$ for an appropriate real Hilbert space $H$ (see Theorem \ref{t characterizations of real Cstar}), for each natural $n$, we can consider the real Hilbert space $H^{(n)} = \ell_n^2(H),$ and each $a= (a_{ij})\in \mathcal{M}_n(\mathscr{A})$ can be regarded as a bounded linear operator on $H^{(n)}$ defined by $$ (a_{ij}) \left(
\begin{array}{c}
\xi_1 \\
\vdots \\
\xi_n \\
\end{array}
\right) = \left(
\begin{array}{c}
\sum\limits_{j=1}^n a_{1j} (\xi_j) \\
\vdots \\
\sum\limits_{j=1}^n a_{nj} (\xi_j) \\
\end{array}
\right).$$
This assignment defines a $^*$-isomorphism from $\mathcal{M}_n(\mathcal{B}(H))$ onto $\mathcal{B}(H^{(n)})$. Since $A$ is represented as a norm closed self-adjoint subalgebra of some $\mathcal{B}(H)$, it turns out that $\mathcal{M}_n(A)$ can be represented as a real C$^*$-algebra, and this construction does not depend on the representation of $A$ inside $\mathcal{B}(H)$ because the norm is unique on a real C$^*$-algebra (see \cite[Proposition 5.1.9]{LI} and \cite[Corollary 1.2.5]{S}).\smallskip

Suppose that $\mathscr{A}$ and $\mathscr{B}$ are two real or complex C$^*$-algebras. For each bounded linear mapping $\Phi: \mathscr{A}\to \mathscr{B}$ and each natural $n$, we can consider a linear mapping $\Phi_n: \mathcal{M}_n (\mathscr{A}) \to \mathcal{M}_n(\mathscr{B})$ defined by $\Phi_n((a_{ij})) := (\Phi(a_{ij}))$. The mapping $\Phi$ is called \emph{$n$-positive} if $\Phi_n : \mathcal{M}_n(\mathscr{A}) \to \mathcal{M}_n(\mathscr{B})$ is positive. If $\Phi_n$ is $n$-positive for all $n$, then $\Phi$ is said to be \emph{completely positive} (see \cite[Definition IV.3.3]{Tak}). There is a vast literature on completely positive and bounded operators between C$^*$-algebras. Here, we shall limit ourselves to comparing some basic properties in the real and complex settings.\smallskip

Some real C$^*$-algebras already hide a complete structure of complex C$^*$-algebra inside. Recall from Section \ref{subsect: complex structure} that a (real) Banach space $X$ has a complex structure if there exists a bounded linear operator $\sigma:X\to X$ satisfying $\sigma^2=-Id$. One can further define a complex norm on $X$ given by \eqref{norm by sigma}. There exist infinite-dimensional Banach spaces admitting no complex structure, and more surprisingly, as shown by Koszmider, Mart{\'i}n, and Mer{\'i} \cite[Corollaries 2.4 and 3.6]{KoszMarMe}, there exist examples of \emph{extremely noncomplex} Banach spaces, that is, Banach spaces that not only do not admit a bounded linear operator $\sigma$ with $\sigma^2=-Id$, but every bounded linear operator $T$ on such a space satisfies $\|Id+T^2\|= 1+ \|T^2\|$. The results in the just quoted reference show that there are several different compact (Hausdorff) spaces $K$ such that the corresponding real C$^*$-algebra $C(K,\mathbb{R})$ is extremely noncomplex.\smallskip

The existence of a complex structure on a real Banach space $X$ determines the presence of multiplicative real linear functionals on the Banach algebra $B(X)$. As it is masterfully explained by \.{Z}elazko \cite{Ze94} and Mankiewicz \cite{Mank71}, \emph{the existence of a  nontrivial linear multiplicative functional on the Banach algebra of all continuous endomorphisms of a Banach space $\mathscr{X}$ implies that $\mathscr{X}$ is not isomorphic to any finite Cartesian power of any Banach space} (see \cite[Remark 6.4]{Mank71}). It is well known that if $\mathscr{X}$ is a complex Banach space, then there does not exist any real linear multiplicative functional $\varphi: \mathcal{B}(\mathscr{X}) \to \mathbb{R}$; since otherwise {it would contradict the presence of a complex structure on $\mathscr{X}_{r}$.} However, the case of real Banach spaces is a bit different. {There are several folk classic arguments showing that, for $n\geq 2$, $\mathcal{B}(\mathbb{R}^n)$ --aka $M_n(\mathbb{R})$-- does not admit a non-zero multiplicative linear functional. Namely, each non-zero multiplicative functional $\phi : M_n(\mathbb{R})\to \mathbb{R}$ satisfies $\phi (a b)= \phi (b a )$ for all $a,b\in M_n(\mathbb{R})$. This property characterizes the normalized trace, $tr(.),$ on $M_n(\mathbb{R})$ up to a scalar multiple. Therefore, $\phi = tr$ because $\phi(I_n)=1$, contradicting that $tr$ is not multiplicative. Alternatively, for any such functional $\phi$, its kernel would be a proper ideal of $M_n(\mathbb{R}),$ which leads to a contradiction.\smallskip
	
It is further known that for each infinite dimensional complex Hilbert space $H$, $\mathcal{B}(H)$ does not admit a non-zero multiplicative real linear functional. Indeed, if $\phi : B(H)\to \mathbb{C}$ is a non-zero multiplicative real linear functional, we can find two orthogonal infinite projections $p$ and $q$ and a partial isometry $e$ such that $p+ q =Id,$ $ee^* =p$ and $e^* e=q$. These facts together imply that $\phi(p) \phi(q) = \phi(p q)=0$, $\phi (p) = \phi (e e^*) = \phi (e^* e) = \phi (q)$ and $1=\phi (Id) = \phi (p) +\phi (q),$ which is impossible.}\smallskip

Mityagin and Edelstein found an example of a real Banach space $X$ such that $\mathcal{B}(X)$ admits a non-trivial real linear multiplicative functional, they concretely showed that this is the case when $X$ is the James space or the space $C(\Gamma_{\omega_1})$ of all continuous scalar valued functions on the set of ordinals not exceeding the first uncountable ordinal with its usual order topology, equipped with the supremum norm (see \cite{MitiEldel70} or \cite{Mank71, MankTomcz2003}). {However, to the best of our knowledge the first to prove that the James space does not admit a complex structure was J. Dieudonn\'e \cite{DIE}. Apart from the James space, the famous Gowers-Maurey example of a Banach space not having a basic sequence, which in its turn is also heriditarily indecomposable ($HI$), provides yet another example of a Banach space lacking a complex structure}. P. Mankiewicz proved in \cite[Theorem 1.1]{Mank71} (see also \cite[\S 9]{MankTomcz2003}) the existence of a separable superreflexive real Banach space $Y$ with the following properties: \begin{enumerate}[$(1)$]
\item $Y$ has a finite-dimensional decomposition;
\item $\mathcal{B}(Y)$ admits a continuous homomorphism onto the Banach algebra $C(\beta \mathbb{N})$ of all continuous scalar-valued functions on the compactification $\beta \mathbb{N}$ of $\mathbb{N}$ equipped with the supremum norm;
\item For each $t\in \mathbb{R}$ there are a projection $P_t\in \mathcal{B}(Y)$ and a linear multiplicative functional $\phi_t$ on $\mathcal{B}(Y)$ such that for every $t_1,t_2\in \mathbb{R}$, $\phi_{t_1}(P_{t_2})$ is equal to $1$ for $t_1=t_2$ and equal to $0$ otherwise.
\end{enumerate}

Consequently, the space $Y$ constructed by Mankiewicz is not isomorphic to any finite Cartesian power of any Banach space. {The reason being that, for every Banach space $Z$ which is the Cartesian product of $n$ copies of another Banach space ($n \geq 2$) there exists a unital homomorphic embedding of $\mathcal{B}(\mathbb{R}^n)$ into $\mathcal{B}(Z)$, $\mathcal{B}(Y)$ admits many non-zero multiplicative real linear functionals, and $\mathcal{B}(\mathbb{R}^n)$ lacks of non-zero multiplicative functionals.}\smallskip

Let us focus on the real C$^*$-algebra $\mathcal{M}_2(A)$, where $A$ is unital real C$^*$-algebra. Given $\alpha,\beta \in \mathbb{R}$, the matrix
$$ w= w_{\alpha,\beta} = \alpha \left(
\begin{array}{cc}
\11 & 0 \\
0 & \11 \\
\end{array}
\right) + \beta \left(
\begin{array}{cc}
0 & \11 \\
-\11 & 0 \\
\end{array}
\right)\in \mathcal{M}_2(A)$$ satisfies that $w^* w = w w^* = (\alpha^2 + \beta^2) \left(
\begin{array}{cc}
\11 & 0 \\
0 & \11 \\
\end{array}
\right)$. Therefore, for $\alpha^2 + \beta^2 \neq 0$, the matrix $u=\frac{1}{\sqrt{\alpha^2+\beta^2}} w$ is a unitary element in $\mathcal{M}_2(A)$. Since the left (respectively, right) multiplication operator by a unitary element in a real C$^*$-algebra is an isometry, the mapping $$ L_u : \mathcal{M}_2(A)\to \mathcal{M}_2(A), \quad x\mapsto L_u (x) = u x$$ is a surjective linear isometry. Taking $\alpha =0$, $\beta =1$, and $u_0 = \left(
\begin{array}{cc}
0 & \11 \\
-\11 & 0 \\
\end{array}
\right)\in \mathcal{M}_2(A)$, the mapping $\sigma = L_{u_0}$ is an isometry on $\mathcal{M}_2(A)$ with $\sigma^2 = L_{u_0}^2 = L_{u_{0}^2} = -Id$. Therefore, $\mathcal{M}_2(A)$ admits a complex structure. The product by complex scalars given by this structure is defined as follows:
$$ (\alpha + i \beta) x = \alpha x + \beta \sigma(x) = \alpha \left(
\begin{array}{cc}
\11 & 0 \\
0 & \11 \\
\end{array}
\right) x + \beta u_0 x = w_{\alpha,\beta} \ x .$$ Now, by the Gelfand--Naimark axiom, we have
$$\begin{aligned}\| (\alpha + i \beta) x \|^2 &= \| w_{\alpha,\beta} \ x \|^2 = \| x^* w_{\alpha,\beta}^* w_{\alpha,\beta} x\| = \left\|(\alpha^2+\beta^2) x^* \left(
\begin{array}{cc}
\11 & 0 \\
0 & \11 \\
\end{array}
\right) x
\right\| \\
& = (\alpha^2+\beta^2) \|x^* x\| = (\alpha^2+\beta^2) \|x\|^2 = |\alpha + i \beta|^2 \|x\|^2\quad (\alpha + i \beta\in \mathbb{C}),
\end{aligned}$$ witnessing that the norm on $\mathcal{M}_2(A)$ is actually a complex norm (cf. Theorem \ref{t Inglestam real Cstar algebra}). Consequently, $\mathcal{M}_2(A)$ is a complex C$^*$-algebra for the corresponding complex structure that we just defined and the original C$^*$-norm. If $A$ is not unital, then we can consider its unitization.\smallskip

We observe next that $\mathcal{M}_2(A)$ contains the algebraic complexification of $A$ as a C$^*$-subalgebra. Namely, let $$\mathcal{A}_c = \left\{ \left(
\begin{array}{cc}
a & b \\
-b & a \\
\end{array}
\right)\in \mathcal{M}_2(A) : a,b\in A \right\}.$$ Clearly, $\mathcal{A}_c$ is a real closed subspace of $\mathcal{M}_2(A)$.
Elements $a + i b\in A_c = A \oplus i A$ are identified with elements $a \11 + b u_0 \equiv \left(
\begin{array}{cc}
a & b \\
-b & a \\
\end{array}
\right)\in \mathcal{M}_2(A)$, and we note that $$(\alpha + i \beta) \left(
\begin{array}{cc}
a & b \\
-b & a \\
\end{array}
\right) = w_{\alpha,\beta} \left(
\begin{array}{cc}
a & b \\
-b & a \\
\end{array}
\right) = \left(
\begin{array}{cc}
\alpha a - \beta b & \alpha b + \beta a \\
-\alpha b - \beta a & \alpha a - \beta b \\
\end{array}
\right).$$ Therefore, $\mathcal{A}_c$ is a norm closed complex subspace of $\mathcal{M}_2(A)$. Similar arguments to those given above show that $\mathcal{A}_c$ is a complex C$^*$-subalgebra of $\mathcal{M}_2(A)$. It follows that $\mathcal{A}_c$ is isometrically $*$-isomorphic to the complexification of $A$ by the uniqueness of the C$^*$-norm (cf. \cite[Corollary 1.2.5]{S}). This procedure can be compared with the construction in \cite[\S 2]{ChuDangRuVen}.\smallskip

Let $T: A\to B$ be a linear mapping between two real C$^*$-algebras. We say that $T$ is \emph{complexifiably positive} if the canonical complex linear extension $\widehat{T}= T_c : A_c \to B_c$ defined by $T_c (a+ i b) = T(a) + i T(b)$ is positive. Clearly, $T$ is complexifiably positive if it is 2-positive. However, the reciprocal statement is not always true. For example, for $A= \mathcal{M}_2(\mathbb{R})$, the transposition $T: \mathcal{M}_2(\mathbb{R})\to \mathcal{M}_2(\mathbb{R})$ defined by $T((\alpha_{ij})) = (\alpha_{ji})$, is positive but not 2-positive (see \cite[p. 5]{Paulsen}). The complexification of $T$ is precisely the transposition on $\mathcal{M}_2(\mathbb{C})= \mathcal{M}_2(\mathbb{R})_c$, which is positive, and thus $T$ is complexifiably positive.\smallskip

One of the fundamental results on completely positive maps, essentially due to Stinespring \cite{Stinespring}, assures that if $\mathscr{A}$ and $\mathscr{B}$ are two C$^*$-algebras and one of them is commutative, then every positive operator $T: \mathscr{A} \to \mathscr{B}$ is completely positive (see \cite[Theorems 3.9 and 3.11]{Paulsen} or \cite[Corollary IV.3.5 and Proposition IV.3.9]{Tak} as well as \cite{MOS3}). This conclusion does not hold in the real setting (see Example \ref{example second March}).\smallskip

Suppose that $T: A\to B$ is a linear operator between two real C$^*$-algebras, and let us assume that one of them is commutative. Since the complexification of a commutative real C$^*$-algebra is a commutative C$^*$-algebra, we can deduce from the above arguments that the following statements are equivalent: \begin{enumerate}[$(a)$]\item $T$ is complexifiably positive, that is, the natural complex linear extension $T_c : A_c\to B_c$ is positive;
\item $T$ is 2-positive;
\item $T$ is completely positive.
\end{enumerate}

There are many open questions to explore about (completely) positive maps in the setting of real C$^*$-algebras, which are not treated here for the sake of brevity.\smallskip

Most of the procedures described in the preceding paragraphs hold in the wider setting of operator spaces. A (complex) \emph{operator space} is a Banach space $\mathscr{X}$ equipped with an isometric embedding $\mathscr{X} \hookrightarrow B(\mathscr{H})$ into the C$^*$-algebra of all bounded linear operators on some complex Hilbert space $\mathscr{H}$. As commented above, the embedding $\mathscr{X}\hookrightarrow B(\mathscr{H})$ induces a norm on each space $\mathcal{M}_n(\mathscr{X})$ of $n\times n$ matrices with entries in $\mathscr{X}$, obtained by regarding any element of $\mathcal{M}_n(\mathscr{X})$ as an operator acting on the Hilbert space $\mathscr{H}^{(n)}$. The resulting sequence of matrix norms is called the operator space structure of $\mathscr{X}$. Then morphisms between operator spaces are the completely bounded maps, that is, the linear mappings $T:\mathscr{X}\to \mathscr{Y}$ which induce uniformly bounded mappings between the matrix spaces $\mathcal{M}_n(\mathscr{X})$ and $\mathcal{M}_n(\mathscr{Y})$ (cf. \cite{EffRuanBook, Paulsen, PisBook}).\smallskip
	
After fifteen years of successful developing of the theory of complex operator spaces, Ruan, one of the founders of operator space theory, introduced real operator spaces in \cite{RUA03}. A \emph{real operator space} on a real Hilbert space $H$ is a norm closed subspace $V$ of $B(H)$ together with the canonical matrix norm inherited from $B(H)$. According to this definition, every real C$^*$-algebra is a real operator space with a canonical matrix norm (actually, a real C$^*$-algebra matrix norm). Ruan described in this paper representations of $\mathbb{C}$ and of the real quaternion ring $\mathbb{H}$ as real operator spaces, as well as similar procedures to complexify a real operator space as the one discussed above for real C$^*$-algebras. Ruan also proved interesting examples which have no counterpart in the complex case. \smallskip

In a continuation paper Ruan investigated the complexification of a real operator space (see \cite{RUA}). Suppose $V$ is a real operator space, and let $V_c=V + i V$ be the algebraic complexification of $V$. In this case, the norm on the complexification must enjoy additional properties linked to the operator space structure. An operator space structure on $V_c$, given by a sequence of matrix norms $\{\|\cdot\|_n\}$, is called \emph{reasonable} if the mapping $x\mapsto x + i 0$ is a complete isometry from $V$ into $V_c$ and $$ \|x+iy\|_n=\|x- iy\|_n \hbox{}$$ for any $n\geq 1$ and any $x,y\in \mathcal{M}_n (V)$. The main result, established by Ruan \cite{RUA}, proves that $V_c$ admits a unique reasonable operator space structure $\{\|\cdot\|_n\}$. Furthermore, for any $x,y\in \mathcal{M}_n(V)$, $\|x+iy\|_n$ is equal to the norm of the matrix $\left(\begin{array}{cc}
	x & y \\
	-y & x \\
\end{array} \right)$ in $\mathcal{M}_{2n}(V)$. Therefore, up to a complete isometry, there is a unique reasonable complex operator space structure on the complexification of a real operator space. This result is employed to characterize complex operator spaces which can be expressed as the complexification of some real operator space. \smallskip

\subsection{Surjective linear isometries}\label{subsec: surjective linear isometries 1}\ \smallskip

Suppose that $K_1$ and $K_2$ are two compact Hausdorff spaces. Most of basic references and basic courses in functional analysis cover the result known as \emph{Banach--Stone theorem}, which asserts (see \cite{Ba,Sto}) that for each surjective linear isometry $T: C(K_1)\to C(K_2)$, there exist a homeomorphism $\sigma: K_2\to K_1$ and a unimodular (unitary) continuous function $u\in C(K_2)$ such that $$T(f) (s) = u(s) f(\sigma(s)),\quad \hbox{ for all } f\in C(K_1) .$$ The spaces involved in this result are commutative unital C$^*$-algebras, and the conclusion implies that, although not every surjective linear isometry between $C(K)$ spaces preserves the product nor the involution, each one of them preserves products of the form $\{f,g,h\} = f \overline{g} h = f g^* h$, that is, $$T \{f,g,h\} = \{T(f), T(g), T(h)\} .$$ The mapping $T$ is precisely given by a composition operator, $C(K_1)\to C(K_2)$, $f\mapsto f\circ \sigma$, multiplied by a unitary element in $C(K_2)$.\smallskip

In the noncommutative setting, we find one of the most influencing results in the theory of C$^*$-algebras, which was established by Kadison in his study on isometries of operator algebras (cf. \cite{Kad51}). Given two unital C$^*$-algebras $\mathscr{A}$ and $\mathscr{B}$, for each surjective linear isometry $T: \mathscr{A} \to \mathscr{B}$, there exist a unitary element $u$ in $\mathscr{B}$ and a Jordan $^*$-isomorphism $\Phi: \mathscr{A} \to \mathscr{B}$ (i.e., a linear bijection preserving Jordan products $\Phi (a\circ b ) = \Phi (a) \circ \Phi(b)$, where $a\circ b := \frac12 (a b + ba)$ and the involution $\Phi (a^*) = \Phi(a)^*$) such that \begin{equation}\label{eq form a surjective linear isometry Kadison} T (x) = u \Phi (x)\quad \hbox{ for all $x\in \mathscr{A}.$ }
\end{equation} Jordan $^*$-isomorphisms were called C$^*$-isomorphisms by Kadison (see \cite[Theorem 7]{Kad51}).\smallskip

A subsequence result by Paterson and Sinclair \cite{PaSi} indicates that, at the unique cost of considering the unitary $u$ in the
multiplier algebra of the C$^*$-algebra in the codomain, the conclusion in Kadison's theorem remains true for surjective linear isometries between non-necessarily unital C$^*$-algebras $\mathscr{A}$ and $\mathscr{B}$. That is, if $T: \mathscr{A} \to \mathscr{B}$ is a surjective linear isometry, then there exist a unitary element $u$ in $M(\mathscr{B}) = \{ b\in \mathscr{B}^{**} : ba , ab \in \mathscr{B} \hbox{ for all } a\in \mathscr{B}\}$ and a Jordan $^*$-isomorphism $\Phi: \mathscr{A} \to \mathscr{B}$ such that the identity in \eqref{eq form a surjective linear isometry Kadison} holds for all $x\in \mathscr{A}$.\smallskip

A surjective linear isometry $T$ between C$^*$-algebras $\mathscr{A}$ and $\mathscr{B}$ need not preserve, in general, neither associative nor Jordan products. However, it is easy to check from \eqref{eq form a surjective linear isometry Kadison} that any such surjective linear isometry $T$ preserves the triple products of the form $\{a,b,c\} = \frac12 ( a b^* c + c b^* a)$ ($a,b,c\in \mathscr{A}$), that is, $$T \{a,b,c\} = \{T(a),T(b),T(c)\}\quad \hbox{ for all } a,b,c\in \mathscr{A}.$$

Those linear maps preserving the above triple products are called \emph{triple homomorphisms}.\smallskip

The problem of studying those surjective linear isometries between real C$^*$-algebras was addressed by Chu et al. \cite{ChuDangRuVen}, where, in a real tour de force, they obtained the following conclusion.

\begin{theorem}\label{t ChuDangRussoVentura sl isometries real Cstar algebras}\cite[Theorem 6.4]{ChuDangRuVen} Let $A$ and $B$ be real C$^*$-algebras. Suppose that $T: A\to B$ is a surjective linear isometry. Then $T$ preserves triple products of the form $\{a,b,c\} = \frac12 ( a b^* c + c b^* a)$, that is, $$T \{a,b,c\} = \{T(a),T(b),T(c)\}\quad \hbox{ for all } a,b,c\in A.$$ We can actually conclude that $T$ is a triple isomorphism.
\end{theorem}

The original result obtained by Chu et al. does not include a description of the form given by Kadison, Paterson, and Sinclair in \eqref{eq form a surjective linear isometry Kadison}. There is a method to derive this concrete expression. First, we recall that given a bounded linear operator $T$ between real C$^*$-algebras $A$ and $B$ (or between real Banach spaces), finding a norm preserving complex linear extension to the corresponding complexifications is not an easy task, which is actually impossible in some cases (see Example \ref{example second March}). Let us present some case in which this norm preserving extension is possible. The self-adjoint part $\mathscr{A}_{sa}$ of a C$^*$-algebra $\mathscr{A}$ is a closed real subspace of $\mathscr{A}$, which is not, in general, a subalgebra of $\mathscr{A}$. However, if we replace the associative product by the Jordan product $a\circ b = \frac12 (a b + ba)$, which is commutative but non-necessarily associative, $\mathscr{A}_{sa}$ is a norm closed real Jordan subalgebra of $\mathscr{A}$. Kadison \cite[Theorem 2]{Kad52} proved that every surjective (real) linear isometry $T: \mathscr{A}_{sa}\to \mathscr{B}_{sa},$ where $\mathscr{B}$ is another C$^*$-algebra, admits an extension to a surjective complex linear isometry from $\mathscr{A}$ onto $\mathscr{B}$.\smallskip

Let us see how to apply Theorem \ref{t ChuDangRussoVentura sl isometries real Cstar algebras} for our purposes. Let $T: A\to B$ be a surjective isometry between two real C$^*$-algebras. Let $A_c$ and $B_c$ denote the corresponding complexifications, and let $\tau_1$ and $\tau_2$ be conjugate-linear $^*$-automorphisms of order-2 on $A_c$ and $B_c$, respectively, such that $A = (A_c)^{\tau_1}$ and $B = (B_c)^{\tau_2}$.
Since, by Theorem \ref{t ChuDangRussoVentura sl isometries real Cstar algebras}, $T$ preserves triple products of the form $\{a,b,c\} = \frac12 (a b^* c + c b^* a)$, it can be easily checked that $T_c : A_c\to B_c$ is a surjective complex linear mapping preserving triple products. Therefore, $$\hbox{$T_c\{x,x,x\} = \{T_c(x), T_c(x), T_c(x)\}$ \quad for all $x\in A_c$.}$$ Let us observe that the Gelfand--Naimark axiom is equivalent to $\|x\|^3 = \|\{x,x,x\}\|$ for all $x\in A_c$. Thus the inequalities $$\|T_c(x)\|^3 =\| \{T_c(x), T_c(x), T_c(x)\}\| = \| T_c\{x,x,x\}\| \leq \|T_c\| \|x\|^3$$ hold for all $x\in A_c$, which implies that $T_c$ is nonexpansive. We similarly get $\|T_c^{-1}\|\leq 1,$ and thus $T_c$ is an isometry. Therefore, there exist a unitary $u$ in the multiplier algebra of $B_c$ and a Jordan $^*$-isomorphism $\widehat\Phi : A_c\to B_c$ such that $T_c (x) = u \widehat\Phi (x)$ for all $x\in A_c$. By considering $A^{**}$ and $B^{**}$ as a real forms of $A_c^{**}$ and $B_c^{**},$ via conjugate-linear $^*$-automorphisms $\widehat{\tau}_1$ and $\widehat{\tau}_2$ extending $\tau_1$ and $\tau_2$, respectively (see \cite[Theorem 1.6 and its proof]{ChuDangRuVen}), it is easy to check that $u\in B^{**}$ actually lies in the multiplier algebra of $B$, and since the identity $$u \widehat\Phi (a)= T(a) = \tau_2 T(a) = \widehat{\tau}_2 (u) \widehat{\tau}_2 \widehat\Phi (a) = u \widehat{\tau}_2 \widehat\Phi (a)$$ holds for all $a\in A$, the mapping $\Phi= \widehat\Phi|_{A} : A\to B$ is a (real linear) Jordan $^*$-isomorphism and $T (a) = u \Phi (a)$ for all $a\in A$.\smallskip

It seems from the just surveyed results that, in what concerns surjective real linear isometries, there is no substantial difference between real and complex C$^*$-algebras. We see in the next section that this parallelism will vanish when considering more general structures like real JB$^*$-triples.

\subsection{Jordan structures and contractive projections}\ \smallskip

We have already caught a glimpse of the Jordan structure underlying a C$^*$-algebra $\mathscr{A}$ with the triple product defined by $\{a,b,c\} = \frac12 ( a b^* c + c b^* a)$ ($a,b,c\in \mathscr{A}$). The main motivation to introduce (complex) JB$^*$-triples resides in the results of holomorphic theory on arbitrary complex Banach spaces and the seeking of a generalization of the celebrated Riemann mapping theorem to classify bounded symmetric domains in complex Banach spaces of dimension bigger than or equal to 2 (see, for example, the introduction and the main result in \cite{Ka83}). Since this point of view is well referenced in the literature, we shall introduce ourselves to the notion of JB$^*$-triples from another perspective and advance our incursion into the topic of contractive projections. Let $p$ be a rank-one projection in $\mathcal{B}(\mathscr{H}),$ where $\mathscr{H}$ is an infinite-dimensional complex Hilbert space and consider the mapping $P: \mathcal{B}(\mathscr{H})\to \mathcal{B}(\mathscr{H})$ defined by $P(a) = p a$ that is a linear contractive projection whose image is $\mathscr{H}$. It is well known from results due to Gal{\'e}, Ransford, and White \cite{GalRansWhit92} (see also the article \cite{Mathieu89}), that a C$^*$-algebra is reflexive if and only if it is finite-dimensional. Therefore, the image of the projection $P$ is not a C$^*$-algebra. In other words, C$^*$-algebras are not stable under contractive projections.\smallskip

In the commutative setting, Friedman and Russo \cite[Theorem 2]{FriRu82c0} proved that the range of a norm-one projection $P$ on a commutative C$^*$-algebra $\mathscr{A}$ has a ternary product structure for the triple product defined by $$\{a,b,c\}_{_P} :=P (a b^* c) \quad(a,b,c\in \mathscr{A}).$$ This provides a link with the notion of ternary ring of operators studied by Zettl \cite{Zettl}. In the same article, they also described and characterized all such projections in terms of extreme points of the unit ball of the image of the dual, and they gave necessary and sufficient conditions for the range to be isometric to a C$^*$-algebra. Several years earlier, Arazy and Friedman \cite{ArFri78} gave an encyclopedic work, a complete description of all contractive projections on the C$^*$-algebra $K(H)$ of all compact operators on a complex Hilbert space $H$ and on its dual space of all trace class operators on $H$. \smallskip

Before presenting additional results, we introduce some notions and definitions. We recall that a JC-algebra is a norm closed real Jordan subalgebra of the self-adjoint part of some $\mathcal{B}(\mathscr{H}),$ where $\mathscr{H}$ is a complex Hilbert space (see \cite{Topp65,HOS}). Concerning contractive projections, Effros and St{\o}rmer \cite{EffStor79} observed that for each positive unital projection $P$ on a unital C$^*$-algebra $\mathscr{A}$, the image of the hermitian part of $\mathscr{A}$ under $P$ is itself a Jordan algebra when provided with the new Jordan multiplication given by $x\circ_{_P} y:= P(x\circ y)$.\smallskip

A \emph{J$^*$-algebra}, in the sense introduced by Harris \cite{Harr}, is a norm closed complex linear subspace of $\mathcal{B}(\mathscr{H},\mathscr{K})$, the Banach space of all bounded linear operators from a complex Hilbert space $\mathscr{H}$ to a complex Hilbert space $\mathscr{K}$, which is closed under the product $a\mapsto aa^*a$. A J$^*$-algebra is a concrete example of a JB$^*$-triple, in the sense we will see in the next paragraph, and is also known under the name of JC$^*$-triple. Clearly, the class of J$^*$-algebras contains all C$^*$-algebras, all complex Hilbert spaces, and the spaces $B(H,K)$. The next step in our story takes us to another work by Friedman and Russo. In \cite{FriRu85contractive}, these authors proved that the class of J$^*$-algebras is stable under the action of norm-one projections. More concretely, if $P$ is a contractive projection on a J$^*$-algebra $M$, then $P(M)$ is a Jordan triple system with triple product $\{a,b,c\}_{_P} = \frac12 P( a b^* c+ c b^* a)$ ($a,b,c\in P(M)$); and $(P(M),\{\cdot,\cdot,\cdot \}_{_P})$ admits a faithful representation as a J$^*$-algebra.\smallskip

In 1984, Kaup \cite{Ka84} gave an elegant and sharp example of how holomorphy can be applied in functional analysis by proving that the class of JB$^*$-triples is also stable under contractive projections. The result was also independently established by Stach{\'o} \cite{Stacho82}. We have naturally met the elements in the exceptional class of complex Banach spaces called JB$^*$-triples, which were originated in holomorphic theory, and whose definition, from the point of view of functional analysis, can be stated with the algebraic-analytic axioms presented below.\smallskip

A JB$^*$-triple is a complex Banach space $\mathcal{E}$ admitting a continuous triple product $\{ \cdot,\cdot,\cdot\} : \mathcal{E}\times
\mathcal{E}\times \mathcal{E} \to \mathcal{E},$ which is conjugate-linear in the central variable and symmetric and bilinear in the outer variables and satisfies the following conditions:
\begin{enumerate}[{\rm (a)}] \item The triple product satisfies the Jordan identity \begin{equation}\label{eq Jordan identity} L(a,b) L(x,y) = L(x,y) L(a,b) + L(L(a,b)x,y) - L(x,L(b,a)y),
\end{equation} for all $a,b,x,y\in \mathcal{E}$, where $L(a,b)$ is the linear operator on $\mathcal{E}$ defined by $L(a,b) x = \{ a,b,x\};$
\item For each $a\in \mathcal{E}$, the mapping $L(a,a): \mathcal{E}\to \mathcal{E}$ is a hermitian operator with nonnegative spectrum;
\item $\|\{a,a,a\}\| = \|a\|^3$ for all $a\in \mathcal{E}$.
\end{enumerate} We recall that a bounded linear operator $T$ on a complex Banach space $\mathscr{X}$ is said to be \emph{hermitian} if $\| \exp (i \alpha T) \|_{_{\mathcal{B}(\mathscr{X})}} = 1$ for all real $\alpha$, that is, $\exp (i \alpha T)$ is a surjective linear isometry for all real $\alpha$ (see \cite[\S 10 and Corollary 10.13]{BD} or page \pageref{numerical range} for the connections with the numerical range). This is the definition found by Kaup \cite{Ka83} in the study of bounded symmetric domains in arbitrary complex Banach spaces.\smallskip

We have already commented that all J$^*$-algebras ---and in particular, all C$^*$-algebras--- are examples of JB$^*$-triples with the triple product defined by \begin{equation}\label{eq triple product J*algebras} \{a,b,c\}=\frac12 (a b^* c + c b^* a). \end{equation} Let us observe that for this triple product, axiom $(c)$ in the definition of JB$^*$-triple writes in the form $\|a a^* a \| = \|a\|^3,$ which is equivalent to the Gelfand--Naimark axiom. \smallskip

Several Jordan structures have been introduced to provide a mathematical model for the algebra of observables in quantum mechanics, which is the case of Jordan algebras introduced by Jordan, von Neumann, and Wigner \cite{Jordan33,JorvNWig34}. Friedman \cite{Fri94} presented several examples in theoretical physics, where JB$^*$-triples theory plays an essential role. For example, the M\"{o}bius--Potapov--Harris transformations (see \cite{Harr}) of the automorphism group of a bounded symmetric domain occur as transformations of signals in an ideal transmission line and as velocity transformations between two inertial systems in special relativity. The velocity transformation is similar to a conformal map, and the operators occurring in these transformations have a natural physical meaning. The just quoted author struggles to present the theory of Jordan algebraic structures (especially, JB$^*$-triples) from the point of view of mathematical physics (special relativity, spinors, and foundational quantum mechanics), in a clear exposition suitable both for experts and nonexperts in the monograph \cite{Fri05}. Besides the classical applications of Jordan theory in well-established areas of physics, like special relativity including fermions and quantum mechanics, Jordan algebras are also employed in string theory, quantum gravity, and $M$-theory; the interested reader may consult the book of Iord\u{a}nescu \cite{Iord03}.\smallskip

A real or complex Jordan algebra is a non-necessarily associative algebra $B$ over $\mathbb{R}$ or $\mathbb{C}$ whose multiplication, denoted by $\circ$, is commutative and satisfies the \emph{Jordan identity}: \begin{equation}\label{eq Jordan identity in Jordan algebras} ( x \circ y ) \circ x^2 = x\circ ( y\circ x^2 )\quad\hbox{ for all } x,y\in B.
\end{equation} For each element $a$ in a Jordan algebra $B$, the symbol $U_a$ will stand for the linear mapping on $B$ defined by $$U_a (b) := 2(a\circ b)\circ a - a^2\circ b\quad (b\in B).$$ A real or complex Jordan Banach algebra $B$ is a real or complex Jordan algebra together with a complete norm satisfying $\|a\circ b\|\leq \|a\| \cdot \|b\|$ for all $a,b\in B$. A \emph{JB-algebra} is a real Jordan Banach algebra $J$ satisfying the following axioms: \begin{enumerate}[$(i)$]\item $\|a^2\| =\|a\|^2$;
\item $\|a^2\|\leq \|a^2 +b^2\|$ for all $a, b\in J$.
\end{enumerate}

A complex Jordan Banach algebra $B$ admitting an involution $^*$ satisfying \begin{equation}\label{eq GN axiom Jordan}
\| U_{a} (a^*) \| = \|a\|^3,
\end{equation} for all $a\in B$ is called a \emph{JB$^*$-algebra} (see \cite{youngson1978vidav}, \cite[Definition 3.3.1]{Cabrera-Rodriguez-vol1}). As in the case of C$^*$-algebras, the involution in a JB$^*$-algebra is automatically a conjugate-linear isometry (see \cite[Lemma 4]{youngson1978vidav} and also \cite[Proposition 3.3.13]{Cabrera-Rodriguez-vol1}).\smallskip

A non-necessarily associative algebra $A$, with product denoted by juxtaposition, is called \emph{flexible} if it satisfies the ``flexibility'' condition $(ab)a = a(ba)$, for all $a,b\in A$ (cf. \cite[Definition 2.3.54]{Cabrera-Rodriguez-vol1}). The algebra $A$ is said to be a \emph{non-commutative Jordan algebra} (cf. \cite[Definition 2.4.9]{Cabrera-Rodriguez-vol1}) if it is flexible and a \emph{Jordan-admissible algebra} (i.e., $A$ is a Jordan algebra when equipped with the natural Jordan product $a\circ b = \frac12 (a b + b a)$).\smallskip

In coherence with the notation in the associative setting of C$^*$-algebras, the self-adjoint part of a JB$^*$-algebra $B$ will be denoted by $B_{sa}.$ It is known that (real) \emph{JB-algebras} are precisely the self-adjoint parts of JB$^*$-algebras (see \cite{Wri77}). Any JB$^*$-algebra also admits a structure of a JB$^*$-triple when equipped with the triple product defined by \begin{equation}\label{eq triple product JB*-algebras} \{x,y,z\} = (x\circ y^*) \circ z + (z\circ y^*)\circ x - (x\circ z)\circ y^*,
\end{equation} and in particular, $U_a (b)=\{a,b^*,a\}$ (see \cite[Theorem 3.3]{BraKaUp78}). The reader interested in knowing additional details may consult the monographs \cite{HOS,Cabrera-Rodriguez-vol1}.\smallskip

As in \cite[Definition 3.3.1]{Cabrera-Rodriguez-vol1} a \emph{non-commutative JB$^*$-algebra} is a complete normed	non-commutative Jordan complex $^*$-algebra (say $\mathscr{A}$) satisfying the axiom in \eqref{eq GN axiom Jordan}. JB$^*$-algebras are precisely those non-commutative JB$^*$-algebras which are commutative. The involution of every non-commutative JB$^*$-algebra is an isometry (see \cite[Proposition 3.3.13]{Cabrera-Rodriguez-vol1}). Non-commutative JB$^*$-algebras include all alternative C$^*$-algebras. The recent monographs \cite{Cabrera-Rodriguez-vol1, Cabrera-Rodriguez-vol2} contain a thorough study on the theory of non-commutative JB$^*$-algebras, JB$^*$-triples, and their real counterparts. For example, in \cite[Corollary 3.4.7]{Cabrera-Rodriguez-vol1} we can find a Russo–Dye–Palmer-type theorem for unital noncommutative JB$^*$-algebras. Each non-commutative JB$^*$-algebra becomes a JB$^*$-triple under its own norm and the natural triple product \cite[Theorem 4.1.45]{Cabrera-Rodriguez-vol1}.\smallskip

A JBW$^*$-triple is a JB$^*$-triple that is also a dual Banach space. A triple version of the celebrated Sakai's theorem established by Barton and Timoney \cite{BarTi}, asserts that each JBW$^*$-triple admits a unique (isometric) predual and that its triple product is separately weak$^*$ continuous.\smallskip

Now, since the notion of JB$^*$-triple has been presented, we can state the previously advanced result on contractive projections.

\begin{theorem}\label{t contractive projection}{\rm(Contractive projection principle, \cite{Stacho82,Ka84})} Let $P: \mathcal{E}\to \mathcal{E}$ be a contractive projection on a JB$^*$-triple. Then $P(\mathcal{E})$ is a JB$^*$-triple with respect to the triple product $$\hbox{$\{x,y,z\}_{_P} :=P\{x,y,z\}$ \quad {\rm(}$x,y,z\in P(\mathcal{E})${\rm)}.}$$
\end{theorem}

In the previous theorem, the image of $P$ need not be a JB$^*$-subtriple of $E$. However, if $P :M \to M$ is a weak$^*$-continuous contractive
projection on a JBW$^*$-triple, then there exists a JBW$^*$-subtriple $C$ of $M$ such that $C$ is linearly isometrically isomorphic to $P(M),$ and such that $C$ is the image of a weak$^*$-continuous projection on $M$ (see \cite[\S 5]{EdRu96structural} or \cite[Theorem 2]{FriRu87bicontractive}). \smallskip

A projection $P$ on a Banach space $X$ is called bicontractive if $\|P\|\leq 1$ and $\|Id-P\|\leq 1$. For each linear isometry of order-2, $T,$ on $X$, the mapping $P=\frac12 (Id+T)$ is a bicontractive projection. Friedman and Russo \cite{FriRu87bicontractive} established that in the setting of JB$^*$-triples the reciprocal statement is also true.

\begin{theorem}\cite[Proposition 3.1 and Theorem 4]{FriRu87bicontractive} Let $P: \mathcal{E}\to \mathcal{E}$ be a bicontractive projection on a JB$^*$-triple. Then $P(\mathcal{E})$ is a JB$^*$-subtriple of $\mathcal{E}$. Furthermore, there exists a surjective linear isometry of order-2 $T: \mathcal{E}\to \mathcal{E}$ satisfying $P=\frac12 (Id+T)$. The same conclusion holds for duals of JB$^*$-triples and preduals of JBW$^*$-triples.
\end{theorem}

We will see in the next section that the contractive projection principle does not hold for real JB$^*$-triples. Bicontractive projections on real C$^*$-algebras and real JB$^*$-triples have not been fully studied.

\subsection{Back to surjective linear isometries }\ \smallskip

JB$^*$-triples constitute a suitable setting to study real forms. Few classes of complex Banach spaces offer a better algebraic-analytic structure to describe surjective linear isometries. We observed in Subsection \ref{subsec: surjective linear isometries 1} that each surjective linear isometry between real or complex C$^*$-algebras is a triple isomorphism for the natural triple product associated with C$^*$-algebras given in \eqref{eq triple product J*algebras}. An outstanding generalization of the commented results crystallized in a Banach--Stone type theorem for JB$^*$-triples obtained by Kaup (see \cite[Proposition 5.5]{Ka83}).

\begin{theorem}\label{t Kaup-Banach-Stone thm}{\rm(Kaup--Kadison--Banach--Stone theorem, \cite[Proposition 5.5]{Ka83})} Let $T: \mathcal{E}\to \mathcal{F}$ be a linear bijection between JB$^*$-triples. Then $T$ is an isometry if and only if $T$ is a triple isomorphism.
\end{theorem}

Alternative proofs of this result were given by Dang, Friedman, and Russo \cite{DaFriRu} and by Fern{\'a}ndez-Polo, Mart{\'i}nez, and the third author of this paper \cite{FerMarPeGeometric04} (see also \cite[Corollary 3.4]{BePe04}).\smallskip

It should be noted here that the ``only if'' implication in Theorem \ref{t Kaup-Banach-Stone thm} does not hold when the mapping $T$ merely is a real linear bijection (see \cite[Remark 2.7]{Dang92}). However, every surjective real linear isometry $T: \mathcal{E}\to \mathcal{F}$ between complex JB$^*$-triples preserves cubes of elements (i.e., $T\{x,x,x\} = \{T(x), T(x), T(x)\}$ for all $x\in \mathcal{E}$), and if we further assume that $\mathcal{E}^{**}$ does not have a nontrivial Cartan factor of rank-one as a summand , then $T$ is a triple isomorphism (see \cite[Proposition 1.1 and Theorem 3.1]{Dang92} or \cite[Proposition 3.8]{IsKaRo95}) ---actually, $\mathcal{E}$ is the direct sum of two orthogonal JB$^*$-subtriples $\mathcal{E}_{1}$ and $\mathcal{E}_2$ such that $T|_{\mathcal{E}_1}$ is a (complex) linear and $T|_{\mathcal{E}_2}$ is a conjugate-linear homomorphism---.\smallskip

If $\mathcal{E}$ is a JB$^*$-triple, then the complex conjugate $\overline{\mathcal{E}}$ of $\mathcal{E}$, constructed in Subsection \ref{complex conjugation1}, is also a JB$^*$-triple. Thus each conjugation $\tau$ (i.e., a conjugate-linear isometry of period-2) on $\mathcal{E}$ must preserve triple products (see also \cite[Corollary 1.2]{Dang92}), and the real form $\mathcal{E}^{\tau}=\{ x\in \mathcal{E} : \tau (x) =x\}$ is a norm closed real subtriple of $\mathcal{E}$. Contrary to the case of real C$^*$-algebras, no additional assumptions on $\tau$ are required. This is an equivalent re-statement of the definition of real JB$^*$-triple. As defined in \cite{IsKaRo95}, a real Banach space, $E$, together with a trilinear map $\{\cdot,\cdot,\cdot\}:E\times E\times E\rightarrow E$ is called a \emph{real JB$^*$-triple} if there exist a JB$^*$-triple, $\mathcal{E}$, and a real linear isometry, $\lambda : E\to \mathcal{E}$ preserving triple products, that is, $$\lambda \{x,y,z\}=\{\lambda (x),\lambda (y),\lambda (z)\}$$ for all
$x,y,z$ in $E$. As commented above, this is equivalent to say that $E$ is a real form of a complex JB$^*$-triple under a conjugation (see \cite[Proposition 2.2]{IsKaRo95}). A real JBW$^*$-triple is a real JB$^*$-triple, which is also a dual Banach space. The original definition of real JBW$^*$-triples in \cite[Definition 4.1 and Theorem 4.4]{IsKaRo95} requires an extra axiom assuming that the triple product is separately
w$^*$-continuous. This extra axiom was shown to be superfluous in \cite{MarPe}.\smallskip

Clearly, every real C$^*$-algebra is a real JB$^*$-triples; real and complex Hilbert spaces, J$^*$-algebras, JB$^*$-algebras, and JB$^*$-triples are also real JB$^*$-triples. Further examples include the self-adjoint parts of C$^*$-algebras and all JB-algebras. \smallskip

Unfortunately, an equivalent definition of real JB$^*$-triples in terms of a set of algebraic-analytic axioms, like the one we have for real C$^*$-algebras in Theorem \ref{t characterizations of real Cstar}$(3)$--$(6)$, is not known. This is actually one of the current open problems in the theory of JB$^*$-triples. The best positive partial answers were contained in \cite{DaRu94,Pe03} in the cases of commutative triples and real JB$^*$-triples admitting a unitary element, respectively. Let us revisit these concrete results.\smallskip

A real \emph{Jordan Banach triple} is a real Banach space $A$ together with a continuous trilinear product $$ A \times A \times A \rightarrow A, \quad (x,y,z) \mapsto \{x,y,z\}, $$ which is symmetric in the outer variables and satisfies the Jordan identity seen in \eqref{eq Jordan identity}. A similar notion works in the complex setting.\smallskip

A real or complex Jordan Banach triple $A$ is called \emph{commutative} or \emph{abelian} if the identity $$ \{ \{ x,y,z\},u,v\} = \{ x,y,\{ z,u,v\}\} = \{ x,\{ y,z,u\},v\}$$ holds for all $x,y,z,u,v \in A$. An element $u\in A$ is said to be \emph{unitary} if the mapping $L(u,u)$ coincides with the identity map on $A$. In this case, $A$ is a unital Jordan $^*$-algebra with product $x\circ_u y:=\{x,u,y\}$ and the involution $x^{*_u} :=\{u,x,u\}$ ($x,y,\in A$). \smallskip

A first attempt to find an axiomatic definition of real JB$^*$-triples was conducted by Dang and Russo \cite{DaRu94}. These authors proposed the following definition.

\begin{definition}\label{def of real J*B-triples}\cite[Definition 1.3]{DaRu94} A \emph{J$^*$B-triple} is a real Banach space $E$ equipped
with a structure of real Jordan Banach triple satisfying the following axioms:
\begin{enumerate}
\item[{\small (J$^*$B1)}] $\left\| \{x,x,x\} \right\| =\left\| x\right\| ^3$ for all $x$ in ${E};$
\item[{\small (J$^*$B2)}] $\left\| \{x,y,z\} \right\| \leq \left\| x\right\| \left\| y\right\| \left\| z\right\|$
for all $x,y,z$ in ${E};$
\item[{\small (J$^*$B3)}] $\sigma_{\mathcal{B}(E)} (L(x,x)) \subseteq [0, + \infty)$ for
all $x \in E$;
\item[{\small (J$^*$B4)}] $\sigma_{\mathcal{B}(E)} (L(x,y)-L(y,x)) \subseteq i \mathbb{R}$ for all $x,y \in E$.
\end{enumerate}
Here, the symbol $\sigma_{\mathcal{B}(E)} (T)$ stands for the spectrum of $T\in \mathcal{B}(E)$ when the latter is regarded as a unital real Banach algebra (see Subsection \ref{subsec: standard complexification Banach algebras}).
\end{definition}

Each closed subtriple of a J$^*$B-triple is a J$^*$B-triple (see \cite[Remark 1.5]{DaRu94}). The class of J$^*$B-triples encompasses all real and complex C$^*$-algebras and all real and complex JB$^*$-triples. Furthermore, complex JB$^*$-triples are precisely those complex Jordan Banach triples whose underlying real Banach space is a J$^*$B-triple (see \cite[Proposition 1.4]{DaRu94}). The class of (real) J$^*$B-triples is very huge.\smallskip

In the setting of commutative J$^*$B-triples, Dang and Russo proved that their definition coincides with those mathematical objects called real JB$^*$-triples by Isidro, Kaup, and Rodr{\'i}guez-Palacios \cite{IsKaRo95} one year later.

\begin{theorem}\cite[Theorem 3.11]{DaRu94} Let $E$ be a commutative Jordan Banach triple. Then the following statements are equivalent:\begin{enumerate}[$(a)$]\item $E$ is a J$^*$B-triple $E$;
\item The complexification of $E$ is a complex JB$^*$-triple in some
norm extending the norm on $E$, that is, $E$ is a real JB$^*$-triple.
\end{enumerate}
\end{theorem}

The proof of the previous result is based on a local ``Gelfand'' theory for commutative Jordan Banach triple systems. The questions of whether the complexification of every J$^*$B-triple is a complex JB$^*$-triple in some norm extending the original norm, and if the second dual of a J$^*$B-triple is a J$^*$B-triple with a separately weak$^*$-continuous triple product, remain open (see \cite[Problems 1 and 2 in p. 137]{DaRu94}). One can find positive partial answers to these questions in noncommutative structures.\smallskip

We recall first some definitions. Let $\mathscr{B}$ be a JB$^*$-algebra. Clearly, the involution on $\mathscr{B}$ defines a conjugate-linear isometric Jordan $^*$-automorphism of period-2 on $\mathscr{B}$, and the real form $\mathscr{B}_{sa} = \{ a\in \mathscr{B} : a^* = a\}$ is precisely a (real) JB-algebra. If we replace $^*$ by a conjugate-linear isometric Jordan $^*$-automorphism of period-2 on $\mathscr{B}$, then the corresponding real form is called a \emph{real JB$^*$-algebra}. For these concrete models, Alvermann \cite{Alv86} found the following axiomatic definition: A J$^*$B-algebra, in the sense of Alvermann, is a real Jordan algebra $A$ with unit and an involution $^*$ equipped with a complete algebra
norm satisfying the following axioms: \begin{enumerate}[$\checkmark$] \item $\|U_{x} (x^*) \| = \|x\|^3$;
\item $\|x^*\circ x \| \leq \|x^*\circ x + y^*\circ y \|$ for all $x,y\in A$.
\end{enumerate} Alvermann \cite[Theorem 4.4]{Alv86} proved that the norm of each J$^*$B-algebra $A$ can be extended to its complexification $A_c = A + i A$ making the latter a JB$^*$-algebra. Consequently, every J$^*$B-algebra is a real form of a JB$^*$-algebra under a conjugate-linear isometric Jordan $^*$-automorphism of period-2.\smallskip

A (real or complex) \emph{numerical range space} is a (real or complex) Banach space $\mathscr{X}$ with a fixed norm-one element $u\in \mathscr{X}$. The set of \emph{states of $\mathscr{X}$ relative to $u$}, $D(\mathscr{X},u)$, is defined as the nonempty (by virtue of the Hahn--Banach theorem), convex, and weak$^*$-compact subset of $\mathscr{X}^{*}$ defined as $$D(\mathscr{X},u) := \{ \phi \in \mathscr{X}^{*} : \|\phi\| =1, \ \phi (u)=1 \}.$$ For $x\in \mathscr{X}$, the \emph{numerical range}\label{numerical range} of $x$ relative to $u$, $V(\mathscr{X},u,x)$, is defined as the set $V(\mathscr{X},u,x):= \{ \phi (x) : \phi \in D(\mathscr{X},u) \}$. The \emph{numerical radius} of $x$ relative to $u$, $v(\mathscr{X},u,x)$, is given by $$v(\mathscr{X},u,x) := \max \{ | \lambda | : \lambda \in V(\mathscr{X},u,x)\}.$$
It is well known that a bounded linear operator $T$ on a complex Banach space $\mathscr{X}$ is hermitian if and only if $V( \mathcal{B}(\mathscr{X}), Id, T) \subseteq \mathbb{R}$ (see \cite[Corollary 10.13]{BD}).
The \emph{numerical index} of the numerical range space $(\mathscr{X},u)$ is defined as \begin{align*}n(\mathscr{X},u)= n(\mathscr{X}) &:= \inf\{ v(x) : x\in \mathscr{X}, \ \|x\|=1 \} \\&= \max\{ \alpha \geq 0 : \alpha \|x\| \leq v(x) \hbox{ for all } x\in \mathscr{X}\}.\end{align*} The element $u$ is called a \emph{geometrically unitary} element of $X$ if and only if $n(X, u) > 0$. See \cite[\S 2.1]{Cabrera-Rodriguez-vol1} for a complete survey on numerical ranges. \smallskip

Let us revisit the connections with some previous results. As we have already seen in \eqref{eq Bohoneblust Karlin Hilbert}, for each complex Hilbert space $\mathscr{H}$, the inequality $$\frac{1}{2}\|T\| \leq w(T) \leq \|T\| \hbox{ holds for all } T\in \mathcal{B}(\mathscr{H}),$$ where $w(T)$ stands for the spatial numerical radius of $T$ (see \cite[\S 9, Theorems 3 and 4]{BDnr} and page \pageref{page numerical radius} for the connections and coincidence of the spatial numerical radius of an operator $T\in \mathcal{B} (\mathscr{X})$ and its numerical range in $(\mathcal{B}(\mathscr{X}),Id_{\mathscr{X}})$). The celebrated \emph{Bohnenblust--Karlin theorem} \cite{BK1} proves that if $\mathscr{A}$ is a norm-unital (associative) Banach algebra with unit $\textbf{1}$, then the numerical radius is a norm on $\mathscr{A}$, which is equivalent to the original norm of this Banach algebra. Furthermore, $n(\mathscr{A},\textbf{1}) \geq \frac1e,$ and thus $$v(a) \leq \|a\|\leq e \ v(a)$$ for all $a\in \mathscr{A}$ (see \cite[Theorem 2.6.4]{Palmer}). Subsequent results show that the hypothesis concerning the associativity of $\mathscr{A}$ in the Bohnenblust--Karlin theorem can be actually replaced by a weaker condition. Namely, suppose that $\mathscr{B}$ is a norm-unital (non-necessarily associative) normed complex algebra. Then $n(\mathscr{B},\textbf{1})\geq \frac1e$, and thus $$v(a) \leq \|a\|\leq e \ v(a)$$ for all $a\in \mathscr{B}$ (see \cite[Proposition 2.1.11]{Cabrera-Rodriguez-vol1}). For real Banach algebras, this conclusion is not true, in general. A version of the Bohnenblust--Karlin theorem for unital real Banach algebras was explored by Ingelstam \cite{Inglestam62}.\label{label results by Ingelstam}\smallskip

We recall that a real algebra $A$ is of \emph{complex type} if it is the realization of a complex algebra $\mathscr{A}$, that is, $A =\mathscr{A}_r$. We say that $A$ is of \emph{real type} if it is not of complex type. An element $a$ in $A$ is called \emph{right} (\emph{left}) \emph{quasi-regular} if there exists a $b$ such that $a + b - ab = 0$ ($a + b - ba = 0$). A real algebra $A$ is of \emph{strongly real type} if the element $-x^2$ is quasi-regular for every $x\in A$.\smallskip

Ingelstam \cite[Theorem 2]{Inglestam62} proved that the unit element is a vertex point of a unital real Banach algebra $A$ if and only if $\exp (\alpha x)$ is unbounded as a function of $\alpha$ (real) for each $x\neq 0$. The author also showed that if $x\neq 0$ belongs to the radical of a real Banach algebra, then $\exp(\alpha x)$ is unbounded, and that each real Banach algebra of strongly real type with identity has the vertex property (see \cite[Theorems 3 and 4]{Inglestam62}).\smallskip

We return now to the setting of Jordan Banach triple systems. It is shown in \cite{Pe03} that, by adding an additional axiom to the definition of J$^*$B-triples, we can actually conclude that the previously commented question posed by Dang and Russo admits a positive answer in the case of J$^*$B-triples admitting a unitary element.

\begin{theorem}\cite[Theorem 2.6]{Pe03} Let $E$ be a J$^*$B-triple admitting a unitary element $u$. Then the following
assertions are equivalent: \begin{enumerate}[$(a)$]
\item $E$ is a numerically positive real J$^*$B-triple, that is, $E$ satisfies the following additional axiom: $V(\mathcal{B}(E),Id,L(x,x)) \subseteq [0,+\infty)$ for all $x\in E$;
\item $E$ is a J$^*$B-algebra or a unital real JB$^*$-algebra with product $x\circ_u y := \{x,u,y\}$ and involution $x^{*_u}:= \{ u,x,u\}$;
\item $E$ is a real JB$^*$-triple, that is, the complexification of $E$ is a complex JB$^*$-triple in some norm extending the original norm on $E$.
\end{enumerate}
\end{theorem}

The question of whether every numerically positive J$^*$B-triple admitting no unitary elements, is a real JB$^*$-triple remains open.\smallskip

We can now resume our narrative about contractive projections on real structures. For commutative real C$^*$-algebras, it was shown by Chu et al. that the image of a contractive projection is a real JB$^*$-triple.

\begin{theorem}\cite[Proposition 7.4]{ChuDangRuVen} Let $P$ be a contractive projection on a commutative real C$^*$-algebra $A$. Then $P(A)$ is a real JB$^*$-triple for the triple product defined by $\{a,b,c\}_{_P} = P\{a,b,c\}$ for all $a,b,c\in P(A)$.
\end{theorem}

The reader should not get the impression that all previously known results for surjective linear isometries and contractive projections have been confirmed for real C$^*$-algebras. As stated by Chu et al. \cite[Problem 7.5]{ChuDangRuVen}, the following remains a challenging and important open problem in the study of real JB$^*$-triples: Is the range of a contractive projection on a real C$^*$-algebra isometric to a linear subspace of some real C$^*$-algebra, closed for the natural triple product associated with each J$^*$-algebra?\smallskip

It had been conjectured that, as in the complex setting, the image of a real JB$^*$-triple under a contractive linear projection is a real JB$^*$-triple with respect to the projected product. However, in 2002, Stach{\'o} \cite[Proposition 2.1]{Stacho02} found a counterexample of a contractive real linear projection on a four-real-dimensional JB$^*$-triple whose image is not a real JB$^*$-triple for the projected triple product because the projected triple product violates the Jordan identity. Let us observe that the counterexample found by Stach{\'o} is a rank-one JB$^*$-triple. We do not know whether the result holds for real JB$^*$-triples not admitting rank-one real or complex Cartan factors as summands in their bidual spaces.\smallskip

We have already seen how rank-one Cartan factors and JB$^*$-triples produce subtle problems to determine an algebraic characterization of surjective linear isometries (see \cite[Remark 2.7]{Dang92}) and contractive projections. In what concerns surjective linear isometries between real JB$^*$-triples, this seems to be the unique obstacle to getting triple isomorphisms. We conclude this article with the most general answer known in this line until this moment.

\begin{theorem}\label{t BS real}\cite[Theorem 3.2 and Corollary 3.4]{FerMarPe04} Let $T : E \to F$ be a surjective linear isometry between two real JB$^*$-triples. Suppose that $E^{**}$ does not contain (real or complex) rank-one Cartan factors as direct summands in its atomic part. Then $T$ is a triple isomorphism. Consequently, every surjective linear isometry between two J$^*$B-algebras is a real triple isomorphism.
\end{theorem}

Let us observe that Theorem \ref{t BS real} implies that under the corresponding hypotheses, each surjective linear isometry $T : E \to F$ admits an extension to a surjective complex linear isometry between the complexifications.\smallskip

Another interesting real structure worth to be considered by itself is the class of real non-commutative JB$^*$-algebras. As well as real C$^*$-algebras and real JB$^*$-algebras are defined as closed real $^*$-subalgebras of (complex) C$^*$- and JB$^*$-algebras, respectively, a real non-commutative JB$^*$-algebra is a closed real $^*$-subalgebra of a (complex)  non-commutative JB$^*$-algebra (cf. \cite[Definition 4.2.45]{Cabrera-Rodriguez-vol1}). Every real non-commutative JB$^*$-algebra becomes a real
JB$^*$-triple under its own norm and the same triple product employed in the complex case \cite[Example 4.2.51]{Cabrera-Rodriguez-vol1}. \smallskip

A norm-one element $x$ in a real or complex Banach space $\mathscr{X}$ is called a \emph{vertex} of the closed unit ball of $\mathscr{X}$ (respectively, a \emph{geometric unitary} of $\mathscr{X}$) if the set $D(\mathscr{X},x),$ of all states of $\mathscr{X}$ relative to $x,$ separates the points of $\mathscr{X}$ (respectively, spans $\mathscr{X}^*$).\smallskip

Many results have been derived from the celebrated paper of Kadison \cite{Kad51}[42], on surjective linear isometries of C$^*$-algebras; one of them is an implicit Banach space characterization of unitary elements in unital C$^*$-algebras. It is well explained by Rodr{\'i}guez-Palacios \cite{Rod2010} that the mentioned characterization can be deduced from results of Kadison as well as  Bohnenblust and Karlin \cite{BK1}, and an explicit statement was included by Akemann and Weaver 	 \cite{AkWea2002}.

\begin{theorem}\label{t AKWeRod unitaries Cstar}\cite[Theorem 2]{AkWea2002} \cite[Theorem 2.1]{Rod2010} Let $\mathscr{A}$ be a unital C$^*$-algebra, and let $u$ be a norm-one element of $\mathscr{A}$. Then the following conditions are equivalent:\begin{enumerate}[$(1)$] \item $u$ is unitary; \item $u$ is a geometric unitary of $\mathscr{A}$;
\item $u$ is a vertex of the closed unit ball of $\mathscr{A}$.
\end{enumerate}
\end{theorem}

An element $u$ in a real or complex JB$^*$-triple $\mathscr{E}$ is called \emph{unitary tripotent} or \emph{unitary} if $L(u,u)$ is the idenity mapping on $\mathscr{E}$, that is, $\{u,u,x\}= x$ for all $x\in \mathscr{E}$. This definition produces no contradiction when unital JB$^*$-algebras are regarded as JB$^*$-triples because unitary elements in a unital JB$^*$-algebra $\mathscr{A}$ are precisely the unitary tripotents in $\mathscr{A}$ when the latter is regarded as a JB$^*$-triple (cf. \cite[Proposition 4.3]{BraKaUp78}).\smallskip

As shown in \cite[Theorem 3.1]{Rod2010} and \cite[Theorem 4.2.24]{Cabrera-Rodriguez-vol1}, the conclusion in Theorem \ref{t AKWeRod unitaries Cstar} remains true when the C$^*$-algebra $\mathscr{A}$ is replaced by a JB$^*$-triple. However, in the real setting the conclusions are rather different.\smallskip

The case of JB-algebras was treated by Leung, Ng, and Wong \cite{LeNgWon2009}. An element $s$ in a unital JB-algebra is called a \emph{symmetry} if $s^2=\textbf{1}$.

\begin{theorem}\label{t unitaries JB-algebras Leung Wong} \cite[Theorem 2.6]{LeNgWon2009}, \cite[Proposition 3.1.15]{Cabrera-Rodriguez-vol1} Suppose $x$ is a norm-one element in a JB-algebra $N$, then the following statements are equivalent:\begin{enumerate}[$(a)$]\item $x$ is a geometric unitary in $N$;
		\item $x$ is a vertex of the closed unit ball of $N$;
		\item $x$ is an isolated point of the set Symm$(N)$ of all symmetries in $N$ (endowed with the norm topology);
		\item $x$ is a central unitary in $N$;
		\item The multiplication operator $M_x: z\mapsto x\circ z$ satisfies $M_x^2 = \hbox{id}_{N}$,
	\end{enumerate}
\end{theorem}

In the case of real JB$^*$-triples, it is shown in \cite{FerMarPeGeometric04} that the existence of a geometrically unitary element in a real JB$^*$-triple $E$ is equivalent to the fact that $E$ is triple-isomorphic to a unital JB-algebra.

\begin{theorem}\label{t unitaries real JBstar triples}\cite[Proposition 2.8]{FerMarPeGeometric04}, \cite[Theorem 4.2.53]{Cabrera-Rodriguez-vol1}
Let $E$ be a real JB$^*$-triple, and let $u$ be a norm-one element in $E$. Then the following
conditions are equivalent:\begin{enumerate}[$(1)$]\item  $u$ is a geometrically unitary element of $E$;
\item $u$ is a vertex of the closed unit ball of $E$;
\item The Banach space of $E$, endowed with the product $x \circ y := \{x,u,y\}$, becomes a JB-algebra with unit $u$.
\end{enumerate}	
\end{theorem}

A unitary element in a general real JB$^*$-triple need not be, in general, a vertex nor a geometric unitary. The previous theorem should be compare with the conclusions of Ingelstam's version of the Bohnenblust--Karlin theorem for unital real Banach algebras (see page \pageref{label results by Ingelstam}).\smallskip

Certain classical properties of C$^*$- and JB$^*$-algebras have been shown to be true for real non-commutative JB$^*$-algebras. For example, every Jordan $^*$-homomorphism between real non-commutative JB$^*$-algebras is automatically contractive. Furthermore, every Jordan $^*$-monomorphism between real non-commutative JB$^*$-algebras is an isometry \cite[Proposition 5.1.47]{Cabrera-Rodriguez-vol2}. Each
closed ideal of a real non-commutative JB$^*$-algebra is $^*$-invariant or self-adjoint (cf.\cite[Proposition 5.1.48]{Cabrera-Rodriguez-vol2}). Actually, closed ideals
of a real non-commutative JB$^*$-algebra are $M$-ideals \cite[Proposition 5.1.53]{Cabrera-Rodriguez-vol2}. A version of Sakai's theorem for real non-commutative JB$^*$-algebra, in the line of \cite[Theorem 6.1.7]{LI}, \cite[Theorem 1.11]{IsRo} and \cite{MarPe}, is established in \cite[Proposition 5.7.62]{Cabrera-Rodriguez-vol2}.\smallskip

Let $\mathscr{X}$ be a real or complex Banach space. The \emph{spatial numerical range} of such an operator $T\in \mathcal{B}(\mathscr{X})$ is the subset $V (T )\subseteq \mathbb{K}$ defined by
$$V (T ):=\{\varphi (T(x)) : x\in \mathscr{X}, \ \varphi\in \mathscr{X}^*,\  1=\|x\|=\|\varphi\|= \varphi(x)\}.$$ The closed convex hull of the spatial numerical range coincides with the numerical range of $T$ in the $(\mathcal{B}(\mathscr{X}), Id_{\mathscr{X}})$, and thus the numerical radius of $T$ in $(\mathcal{B}(\mathscr{X}), Id_{\mathscr{X}})$ coincides with the supremum of the modulus of those elements in its spatial numerical range\label{page numerical radius} (cf. \cite[Theorems 9.3 and 9.4]{BDnr} and \cite[\S 2.1 ]{Cabrera-Rodriguez-vol1}). The \emph{numerical index} of the space $\mathscr{X}$, $n(\mathscr{X})$, is defined as $$n(\mathscr{X)} = \inf\{ v(T) : T \in \mathcal{B}(\mathscr{X}), \ \|T\| =1\}.$$

Several isomorphic properties of real infinite-dimensional Banach spaces with numerical index $1$ were obtained by L{\'o}pez, Mart{\'i}n
and Pay{\'a} in \cite{LoMarPa}. The next result is a consequence of the just quoted study.

\begin{theorem}\label{thm LoMarPay}\cite{LoMarPa}
Every reflexive real Banach space with numerical index 1 must be finite-dimensional.
\end{theorem}

The result by L{\'o}pez, Mart{\'i}n and Pay{\'a} is deeper and shows that an infinite-dimensional real Banach space with numerical index 1 satisfying the Radon-Nikod{\'y}m property contains $\ell_1$ \cite[Theorem 3]{LoMarPa}. In particular a reflexive or quasi-reflexive real Banach space cannot be re-normed to have numerical index 1, unless it is finite-dimensional.\smallskip

For complex spaces, the existence of reflexive infinite-dimensional Banach spaces with numerical index 1 is a long standing open problem. The problem is related to the validity of the technical Proposition 2 in \cite{LoMarPa} in the complex setting.

\smallskip\smallskip

\textbf{Acknowledgements} We would like to thank V. Kadets, M. Mart{\'i}n, A. Rodr{\'i}guez-Palacios, B. Russo and  A. Sofi for several constructive comments and hints during the redaction of this paper. We also appreciate the valuable constructive comments from the reviewers.\smallskip\smallskip

The first author was supported by a grant from Ferdowsi University of Mashhad (No. 2/55373). The third author was partially supported by the Spanish Ministry of Science, Innovation and Universities (MICINN) and European Regional Development Fund project no. PGC2018-093332-B-I00, Junta de Andaluc{\'i}a grants number A-FQM-242-UGR18, PY20$\underline{\ }$00255 and FQM375, and  by the IMAG--Mar{\'i}a de Maeztu grant CEX2020-001105-M/AEI/10.13039/501100011033. The second and fourth authors were supported by the Spanish Ministry of Science, Innovation and Universities (MICINN) no. PGC2018-097286-B-I00. \smallskip

\bibliographystyle{amsplain}

\begin{thebibliography}{999}

\bibitem{ABR} Y.A. Abramovich, C.D. Aliprantis, G. Sirotkin, V.G. Troitsky, Some open problems and conjectures associated with the invariant subspace problem, \emph{Positivity} \textbf{9} (2005), no. 3, 273--286.

\bibitem{Aron} M.D. Acosta, R.M. Aron, D. Garc\'{i}a, M. Maestre, The Bishop-Phelps-Bollob\'{a}s theorem for operators, \textit{J. Funct. Anal.} \textbf{254} (2008), no. 11, 2780--2799.

\bibitem{AkWea2002} C.A. Akemann, N. Weaver, Geometric characterizations of some classes of operators in C$^*$-algebras and von Neumann algebras, \textit{Proc. Amer. Math. Soc.} \textbf{130}, no. 10  (2002), 3033--3037.

\bibitem{AO1953} A. Alexiewicz, W. Orlicz, Analytic operations in real Banach spaces, \textit{Studia Math.} \textbf{14} (1953), 57--78.

\bibitem{Alv86} K. Alvermann, Real normed Jordan algebras with involution, \textit{Arch. Math.} \textbf{47} (1986), 135--150.

\bibitem{AR2006} V. Anagnostopoulos, S. R\'ev\'esz, Polarization constants for products of linear functionals over ${\mathbb R}^2$ and ${\mathbb C}^2$ and Chebyshev constants of the unit sphere, \textit{Publ. Math. Debrecen} \textbf{68} (2006), no. 1-2, 63--75.

\bibitem{MOS1} L. Aramba\v{s}i\'{c}, D. Baki\'{c}, M.S. Moslehian, A treatment of the Cauchy-Schwarz inequality in $C^*$-modules, \textit{J. Math. Anal. Appl.} \textbf{381} (2011), no. 2, 546--556.

\bibitem{Araujo2020} G. Ara\'{u}jo, P.H. Enflo, G.A. Mu\~{n}oz-Fern\'{a}ndez, D.L. Rodr\'{\i}guez-Vidanes, J.B. Seoane-Sep\'{u}lveda, Quantitative and qualitative estimates on the norm of products of polynomials, \textit{Israel J. Math.} \textbf{236} (2020), no. 2, 727--745.

\bibitem{ArFri78} J. Arazy, Y. Friedman, Contractive projections in $C_1$ and $C_0$, \textit{Mem. Amer. Math. Soc.} \textbf{13} (1978), no. 200, {\rm iv}+165 pp.


\bibitem{AdR1998} J. Arias-de-Reyna, Gaussian variables, polynomials and permanents, \textit{Linear Algebra Appl.} \textbf{285} (1998), 107--114.

\bibitem{01} R.M. Aron, L. Bernal-Gonz\'alez, D.M. Pellegrino, J.B. Seoane-Sep\'ulveda, Lineability: the search for linearity in mathematics. Monographs and Research Notes in Mathematics. CRC Press, Boca Raton, FL, 2016. xix+308 pp. ISBN: 978-1-4822-9909-0.

\bibitem{02} R.M. Aron, V.I. Gurariy, J.B. Seoane-Sep\'ulveda, Lineability and spaceability of sets of functions on $\mathbb{R}$, \textit{Proc. Amer. Math. Soc.} \textbf{133} (2005), no. 3, 795--803.

\bibitem{aronboydryanzalduendo2003} R. M. Aron, C. Boyd, R. A. Ryan, I. Zalduendo, Zeros of polynomials on Banach spaces: the real
story, \textit{Positivity} \textbf{7} (2003), no. 4, 285--295.

\bibitem{arongonzalozagorodnyuk2000} R. M. Aron, R. Gonzalo, A. Zagorodnyuk, Zeros of real polynomials, \textit{Linea Multilinear Algebra} \textbf{48} (2000), no. 2, 107--115.

\bibitem{aronhajek2007}
R. M. Aron, P. H\'ajek, Odd degree polynomials on real Banach spaces, \textit{Positivity} \textbf{11} (2007), no. 1,
143--153.
	
\bibitem{aronhajek2006}
R. M. Aron, P. H\'ajek, Zero sets of polynomials in several variables, \textit{Arch. Math.} (Basel) \textbf{86} (2006), no. 6, 561--568.

\bibitem{AronRueda} R.M. Aron, M.P. Rueda, A problem concerning zero-subspaces of homogeneous polynomial, Dedicated to Professor Vyacheslav Pavlovich Zahariuta, \textit{Linear Topol. Spaces Complex Anal.} \textbf{3} (1997), 20--23.

\bibitem{Aupe91} B. Aupetit, \textit{A primer on spectral theory}, Universitext. Springer-Verlag, New York, 1991.

\bibitem{A2016}  Avil\'es, Antonio; S\'anchez, F\'elix Cabello; Castillo, Jes\'us M. F.; Gonz\'alez, Manuel; Moreno, Yolanda Separably injective Banach spaces. Lecture Notes in Mathematics, 2132. Springer, [Cham], 2016. xxii+217 pp.

\bibitem{avilestodorcevic2009} A. Avil\'es, S. Todorcevic, Zero subspaces of polynomials on $l_1(\Gamma)$, \textit{J. Math. Anal. Appl.} \textbf{350} (2009), no. 2, 427--435.

\bibitem{Ball1991} K.M. Ball, The plank problem for symmetric bodies, \textit{Invent. Math.} \textbf{104} (1991), no. 3, 535--543.

\bibitem{BAA} K.M. Ball, The complex plank problem, \textit{Bull. London Math. Soc.} \textbf{33} (2001), no. 4, 433--442.

\bibitem{Ba} S. Banach, \textit{Th{\' e}orie des op{\' e}rations lin{\' e}aires}, Warsaw, 1932.

\bibitem{Banach} S. Banach, \"Uber homogene Polynome in $(L^2)$, \textit{Studia Math.} \textbf{7} (1938), 36--44.

\bibitem{Banakh} T. Banakh, Every 2-dimensional Banach space has the Mazur-Ulam property. arXiv:2103.09268

\bibitem{BanakhCabello} T. Banakh, J. Cabello S{\'a}nchez, Every non-smooth 2-dimensional Banach space has the Mazur-Ulam property,\textit{ Linear Algebra Appl.} \textbf{625} (2021), 1--19.

\bibitem{Bang} T. Bang, A solution of the ``plank problem'', \textit{Proc. Amer. Math. Soc.} \textbf{2} (1951), 990--993.

\bibitem{BarTi} T. Barton, R.M. Timoney, Weak$^*$-continuity of Jordan triple products and its applications, \textit{Math. Scand.} \textbf{59} (1986), 177--191.

\bibitem{bps} F. Bayart, D. Pellegrino, J.B. Seoane-Sep\'ulveda, The Bohr radius of the $n$-dimensional polydisk is equivalent to $\log n/n$, \textit{Adv. Math.} \textbf{264} (2014), 726--746.

\bibitem{BeBuKaRod08} J. Becerra-Guerrero, M.J. Burgos, A. Kaidi, A. Rodr{\'i}guez Palacios, Banach spaces whose algebras of operators have a large group of unitary elements, \textit{Math. Proc. Cambridge Philos. Soc.} \textbf{144} (2008), no. 1, 97--108.

\bibitem{BeCuFerPe2021} J. Becerra-Guerrero, M. Cueto-Avellaneda, F.J. Fern{\'a}ndez-Polo, A.M. Peralta, On the extension of isometries between the unit spheres of a JBW$^*$-triple and a Banach space, \textit{J. Inst. Math. Jussieu} \textbf{20} (2021), no. 1, 277--303.

\bibitem{BePe04} J. Becerra-Guerrero, A.M. Peralta, Subdifferentiability of the norm and the Banach-Stone theorem for real and complex JB$^*$-triples, \textit{Manuscripta Math.} \textbf{114} (2004), no. 4, 503--516.

\bibitem{Sarant1993} C. Ben\'itez, Y. Sarantopoulos, Characterization of real inner product spaces by means of symmetric bilinear forms, \textit{J. Math. Anal. Appl.} \textbf{180} (1993), no. 1, 207--220.

\bibitem{BST1998} C. Ben\'itez, Y. Sarantopoulos, A. Tonge, Lower bounds for norms of products of polynomials \textit{Math. Proc. Cambridge Philos. Soc.} \textbf{128} (1998), no. 3, 395--408.

\bibitem{RACSAM2021} L. Bernal-Gonz\'alez, H. J. Cabana, D. Garc\'ia, M. Maestre, G. A. Mu\~{n}oz-Fern\'andez, J.B. Seoane-Sep\'ulveda, A new approach towards estimating the $n$-dimensional Bohr radius, \textit{Rev. R. Acad. Cienc. Exactas Fís. Nat. Ser. A Mat. RACSAM} \textbf{115} (2021), no. 2, Paper No. 44, 10 pp.

\bibitem{03} L. Bernal-Gonz\'alez, D.M. Pellegrino, J.B. Seoane-Sep\'ulveda, Linear subsets of nonlinear sets in topological vector spaces, \textit{Bull. Amer. Math. Soc. (N.S.)} 51 (2014), no. 1, 71--130.

\bibitem{BP} E. Bishop, R.R. Phelps, A proof that every Banach space is subreflexive, \textit{Bull. Amer. Math. Soc.} \textbf{67} (1961), 97--98.

\bibitem{Boas} H.P. Boas, The football player and the infinite series, \textit{Notices Amer. Math. Soc.} \textbf{44} (1997), no. 11, 1430--1435.

\bibitem{boas} H.P. Boas, D. Khavinson, Bohr's power series theorem in several variables, \textit{Proc. Amer. Math. Soc.} \textbf{125} (1997), no. 10, 2975--2979.

\bibitem{Bochnak1970} J. Bochnak, Analytic functions in Banach spaces, \textit{Studia Math.} \textbf{35} (1970), 273--292.

\bibitem{Bochnak1971bis} J. Bochnak, J. Siciak, Polynomials and multilinear mappings in topological vector spaces, \textit{Studia Math.} \textbf{39} (1971) 59--76.

\bibitem{Bochnak1971} J. Bochnak, J. Siciak, Analytic functions in topological vector spaces, \textit{Studia Math.} \textbf{39} (1971), 77--112.

\bibitem{bh} H.F. Bohnenblust, E. Hille, On the absolute convergence of Dirichlet series, \textit{Ann. of Math.} (2) \textbf{32} (1931), no. 3, 600--622.

\bibitem{BK1} H.F. Bohnenblust, S. Karlin, Geometrical properties of the unit sphere of Banach algebras, \textit{Ann. of Math.} (2) \textbf{62} (1955), 217--229.

\bibitem{BDnr} F.F. Bonsall, J. Duncan, \textit{Numerical ranges of operators on normed spaces and of elements of normed algebras}, London Mathematical Society Lecture Note Series, 2 Cambridge University Press, London-New York, 1971.

\bibitem{BD} F.F. Bonsall, J. Duncan, \textit{Complete normed algebras}, Ergebnisse der Mathematik und ihrer Grenzgebiete, Band 80. Springer-Verlag, New York-Heidelberg, 1973.

\bibitem{BOU} J. Bourgain, Real isomorphic complex Banach spaces need not be complex isomorphic, \textit{Proc. Amer. Math. Soc.} \textbf{96} (1986), no. 2, 221--226.

\bibitem{BR1998} C. Boyd, R. Ryan, The norm of the product of polynomials in infinite dimensions, \textit{Proc. Edinb. Math. Soc.} (2) \textbf{49} (2006), no. 1, 17--28.

\bibitem{Boyd2018} C. Boyd, R.A. Ryan, N. Snigireva, Radius of analyticity of analytic functions on Banach spaces, \textit{ J. Math.Anal.Appl.} \textbf{463} (2018) 40--49.

\bibitem{BraRo} O. Bratteli, D.W. Robinson, \textit{Operator algebras and quantum statistical mechanics I}, Springer Verlag,
New York, 1979.

\bibitem{BraKaUp78} R. Braun, W. Kaup, H. Upmeier, A holomorphic characterisation of Jordan-C$^*$-algebras, \textit{Math. Z.} \textbf{161} (1978), 277--290.

\bibitem{BritMabrTou2021} R. Brits, M. Mabrouk, C. Tour{\'e}, A multiplicative Gleason--Kahane--\.{Z}elazko theorem for C$^*$-algebras, \textit{J. Math. Anal. Appl.} \textbf{500} (2021), no. 1, 125089, 4 pp.


\bibitem{BurFerGarPe2015} M.J. Burgos, F.J. Fern{\'a}ndez-Polo, J.J. Garc{\'e}s, A.M. Peralta, A Kowalski-S{\l}odkowski theorem for 2-local $^*$-homomorphisms on von Neumann algebras, \textit{Rev. R. Acad. Cienc. Exactas Fís. Nat. Ser. A Mat. RACSAM} \textbf{109} (2015), no. 2, 551--568.

\bibitem{BurFerGarPe2015tripl} M.J. Burgos, F.J. Fern{\'a}ndez-Polo, J.J. Garc{\'e}s, A.M. Peralta, 2-local triple homomorphisms on von Neumann algebras and JBW$^*$-triples, \textit{J. Math. Anal. Appl.} \textbf{426} (2015), no. 1, 43--63.


\bibitem{CabMol2002} F. Cabello S{\'a}nchez, L. Moln{\'a}r, Reflexivity of the isometry group of some classical spaces, \textit{Rev. Mat. Iberoamericana} \textbf{18} (2002), no. 2, 409--430.

\bibitem{CabSan2004} F. Cabello S{\'a}nchez, The group of automorphisms of $L_{\infty}$ is algebraically reflexive, \textit{Studia Math.} \textbf{161} (2004), no. 1, 19--32.

\bibitem{Cabrera-Rodriguez-vol1} M. Cabrera~Garc\'{\i}a, A. Rodr\'{\i}guez~Palacios, \textit{Non-associative normed algebras. {V}ol. 1}, The Vidav-Palmer and Gelfand-Naimark theorems, vol.~154 of {\em Encyclopedia of Mathematics and its Applications}, Cambridge University Press, Cambridge, 2014.

\bibitem{Cabrera-Rodriguez-vol2} M. Cabrera~Garc\'{\i}a, A. Rodr\'{\i}guez~Palacios, \textit{Non-associative normed algebras. Vol. 2.}
Representation theory and the Zel'manov approach. Encyclopedia of Mathematics and its Applications, 167. Cambridge University Press, Cambridge, 2018.

\bibitem{CPR2013} D. Carando, D. Pinasco, J.T. Rodr\'iguez, Lower bounds for norms of products of polynomials on $L_ p$ spaces, \textit{Studia Math.} \textbf{214} (2013), 157--166.

\bibitem{CPR2017} D. Carando, D. Pinasco, J.T. Rodr\'iguez, On the linear polarization constants of finite dimensional spaces, \textit{Math. Nachr.} \textbf{290} (2017), no. 16, 2547--2559.

\bibitem{CPR2017_2} D. Carando, D. Pinasco, J.T. Rodr\'iguez, Non-linear plank problems and polynomial inequalities, \textit{Rev. Mat. Complut.} \textbf{30} (2017), no. 3, 507--523.

\bibitem{CP2011} I. Chalendar, J.R. Partington, Modern approaches to the invariant-subspace problem, \textit{Cambridge Tracts in Mathematics} \textbf{188}, Cambridge University Press. Cambridge, 2011.

\bibitem{Cho} M. Cho, An elementary proof of Gleason--Kahane--\.{Z}elazko's theorem for complex Banach algebra with a Hermitian involution,
\textit{Sci. Rep. Niigata Univ. Ser. A} No. \textbf{11} (1974), 1--4.

\bibitem{ChodaNakamura71} H. Choda, M. Nakamura, Elementary proofs of Gleason--Kahane--\.{Z}elazko's theorem for B$^*$-algebras, \textit{Mem. Osaka Kyoiku Univ.} III Natur. Sci. Appl. Sci. \textbf{20} (1971), 111--112.

\bibitem{Miguel} Y.S. Choi, S.K. Kim, H.J. Lee, M. Mart\'{i}n, The Bishop-Phelps--Bollob\'{a}s theorem for operators on $\mathcal{L}_1(\mu)$, \textit{J. Funct. Anal.} \textbf{267} (2014), no. 1, 214--242.

\bibitem{ChuDangRuVen} C.H. Chu, T. Dang, B. Russo, B. Ventura, Surjective isometries of real C$^*$-algebras, \textit{J. London Math. Soc.}, \textbf{47} (1993), 2, 97--118.

\bibitem{04} K.C. Ciesielski, J.B. Seoane-Sep\'ulveda, Differentiability versus continuity: restriction and extension theorems and monstrous examples, \textit{Bull. Amer. Math. Soc. (N.S.)} \textbf{56} (2019), no. 2, 211--260.

\bibitem{Vandercorput} J.G. Van der Corput, G. Schaake, Berichtigung zu: Ungleichungen f\"ur Polynome und trigonometrische Polynome. (German) \textit{Compositio Math.} \textbf{3} (1936), 128--128.

\bibitem{CUE} C.W. Cuellar, A Banach space with a countable infinite number of complex structures, \textit{J. Funct. Anal.} \textbf{267} (2014), no. 5, 1462--1487.

\bibitem{Cuellar} W. Cuellar Carrera, Complex structures on Banach spaces with a subsymmetric basis, \textit{J. Math. Anal. Appl.} \textbf{440} (2016), no. 2, 624--635.

\bibitem{Da} H.G. Dales, F.K.Jr. Dashiell, A.T.-M. Lau, D. Strauss, \textit{Banach spaces of continuous functions as dual spaces}, CMS Books in Mathematics/Ouvrages de Math{\'e}matiques de la SMC. Springer, Cham, 2016.

\bibitem{Dang92} T. Dang, Real isometries between JB$^*$-triples, \textit{Proc. Amer. Math. Soc.} \textbf{114} (1992), no. 4, 971--980.

\bibitem{DaFriRu} T. Dang, Y. Friedman, B. Russo, Affine geometric proofs of the Banach-Stone theorems of Kadison and Kaup, \textit{Proceedings of the Seventh Great Plains Operator Theory Seminar (Lawrence, KS, 1987). Rocky Mountain J. Math.} \textbf{20} (1990), no. 2, 409--428.

\bibitem{DaRu94} T. Dang, B. Russo, Real Banach Jordan triples, \textit{Proc. Amer. Math. Soc.} \textbf{122} (1994), 135--145.

\bibitem{DF1993} A. Defant, K. Floret, \textit{Tensor norms and operator ideals}, North-Holland Mathematics Studies, 176. North-Holland Publishing Co., Amsterdam, 1993.

\bibitem{annals2011} A. Defant, L. Frerick, J. Ortega-Cerd\`a, M. Ouna\"ies, K. Seip, The Bohnenblust-Hille inequality for homogeneous polynomials is hypercontractive, \textit{Ann. of Math.} (2) \textbf{174} (2011), no. 1, 485--497.

\bibitem{ManDom} A. Defant, D. Garc\'ia, M. Maestre, P. Sevilla-Peris, \textit{Dirichlet series and holomorphic functions in high dimensions}, New Mathematical Monographs, 37. Cambridge University Press, Cambridge, 2019.

\bibitem{DIE} J. Dieudonn\'{e}, Complex structures on real Banach spaces, \textit{Proc. Amer. Math. Soc.} \textbf{3} (1952), 162--164.

\bibitem{Dimant2021} V. Dimant, D. Galicer, J.T. Rodr\'iguez, The polarization constant of finite dimensional complex spaces is one, \textit{Math. Proc. Camb. Phil. Soc}, (2021), 1--19.

\bibitem{Dineen1999} S. Dineen, \textit{Complex analysis on infinite-dimensional spaces}, Springer Monographs in Mathematics. Springer-Verlag London, Ltd., London, 1999.

\bibitem{DD} D. Diniz, G.A. Mu\~{n}oz-Fern{\'a}ndez, D. Pellegrino, J.B. Seoane-Sep\'ulveda, The asymptotic growth of the constants in the Bohnenblust-Hille inequality is optimal, \textit{J. Funct. Anal.} \textbf{263} (2012), 415--428.

\bibitem{diniz2} D. Diniz, G.A. Mu\~{n}oz-Fern\'andez, D. Pellegrino, J.B. Seoane-Sep\'ulveda, Lower bounds for the constants in the Bohnenblust-Hille inequality: the case of real scalars, \textit{Proc. Amer. Math. Soc.} \textbf{142} (2014), 575--580.

\bibitem{DowHuMup96} P.N. Dowling, Z. Hu, D. Mupasiri, Complex convexity in Lebesgue-Bochner function spaces, \textit{Trans. Amer. Math. Soc.} \textbf{348} (1996), no. 1, 127--139.

\bibitem{EdRu96structural} C.M. Edwards, G.T. R\"{u}ttimann, Structural projections on JBW$^*$-triples, \textit{J. London Math. Soc.} \textbf{53} (1996), 354--368.

\bibitem{EffRuanBook} E. Effros, Z.-J. Ruan, \textit{Operator spaces}, London Math. Soc. Monographs, New Series, 23, Oxford University Press, New York, 2000.

\bibitem{EffStor79} E.G. Effros, E. St{\o}rmer, Positive projections and Jordan structure in operator algebras, \textit{Math. Scand.} \textbf{45} (1979), no. 1, 127--138.

\bibitem{Enflo1987} P.H. Enflo, On the invariant subspace problem for Banach spaces, \textit{Acta Math.} \textbf{158} (1987), no. 3-4, 213--313.


\bibitem{ERDOS} P. Erd\"os, Some remarks on polynomials, \textit{Bull. Amer. Math. Soc.}, \textbf{53} (1947), 1169--1176.

\bibitem{FER} V. Ferenczi, Uniqueness of complex structure and real hereditarily indecomposable Banach spaces, \textit{Adv. Math.} \textbf{213} (2007), no. 1, 462--488.

\bibitem{FG} V. Ferenczi, E.M. Galego, {Countable groups of isometries on Banach spaces}, \textit{Trans. Amer. Math. Soc.} \textbf{362} (2010), no. 8, 4385--4431.

\bibitem{FerMarPeGeometric04} F.J. Fern{\'a}ndez-Polo, J. Mart{\'i}nez, A.M. Peralta, Geometric characterization of tripotents in real and complex JB$^*$-triples, \textit{J. Math. Anal. Appl.} \textbf{295} (2004), no. 2, 435--443.

\bibitem{FerMarPe04} F.J. Fern{\'a}ndez-Polo, J. Mart{\' i}nez, A.M. Peralta, Surjective isometries between real JB$^*$-triples, \textit{Math. Proc. Cambridge Phil. Soc.}, \textbf{137} (2004), 709--723.

\bibitem{fernandez2006} M. Fern\'andez-Unzueta, Zeroes of polynomials on $l_\infty$, \textit{J. Math. Anal. Appl.} \textbf{324} (2006), no. 2, 1115--1124.

\bibitem{ferrer2007} J. Ferrer, Zeroes of real polynomials on C(K) spaces, \textit{J. Math. Anal. Appl.} \textbf{336} (2007), no. 2, 788--796.

\bibitem{ferrerprims2007} J. Ferrer, On the zero-set of real polynomials in non-separable Banach spaces, \textit{Publ. Res. Inst. Math. Sci.} \textbf{43} (2007), no. 3, 685--697.

\bibitem{ferrer2009} J. Ferrer, A note on zeroes of real polynomials on $C(K)$ spaces, \textit{Proc. Amer. Math. Soc.} \textbf{137} (2009), no. 2, 573--577.

\bibitem{ferrer2019} J. Ferrer, D. Garc\'ia, M. Maestre, J.B. Seoane-Sep\'ulveda, On the zero-set of 2-homogeneous polynomials in Banach spaces, \textit{Linear Multilinear Algebra} \textbf{67} (2019), no. 10, 1958--1970.

\bibitem{FM} J. Ferrera, G.A. Mu\~{n}oz, A characterization of real Hilbert spaces using the Bochnak complexification norm, \textit{Arch. Math. (Basel)} \textbf{80} (2003), no. 4, 384--392.

\bibitem{Frenkel}
P. E. Frenkel, Hafnians and products of real linear functionals, \textit{Math. Res. Lett.} \textbf{15} (2008), no. 2, 351--358.


\bibitem{Fri94} Y. Friedman, Bounded symmetric domains and the JB$^*$-triple structure in physics, \textit{Jordan algebras (Oberwolfach, 1992)}, 61--82, de Gruyter, Berlin, 1994.

\bibitem{Fri05} Y. Friedman, \textit{Physical applications of homogeneous balls.} With the assistance of Tzvi Scarr. Progress in Mathematical Physics, 40. Birkh\"{a}user Boston, Inc., Boston, MA, 2005.

\bibitem{FriRu82c0} Y. Friedman, B. Russo, Contractive projections on $C_0(K)$, \textit{Trans. Amer. Math. Soc.} \textbf{273} (1982), 57--73.

\bibitem{FriRu85contractive} Y. Friedman, B. Russo, Solution of the contractive projection problem, \textit{J. Funct. Anal.} \textbf{60} (1985), no. 1, 56--79.

\bibitem{FriRu87bicontractive} Y. Friedman, B. Russo, Conditional expectation and bicontractive projections on Jordan C$^*$-algebras and their generalizations, \textit{Math. Z.} \textbf{194} (1987), no. 2, 227--236.

\bibitem{Fukamiya} M. Fukamiya, On a theorem of Gelfand and Neumark and the B$^*$-algebra, \textit{Kumamoto J. Sci.} \textbf{1} (1952), 17--22.

\bibitem{GalRansWhit92} J.E. Gal{\'e}, T.J. Ransford, M.C. White, Weakly compact homomorphisms, \textit{Trans. Amer. Math. Soc.} \textbf{331} (1992), no. 2, 815--824.

\bibitem{GMPS2012} J.L. G{\'a}mez-Merino; G.A. Mu\~noz-Fern\'andez, D. Pellegrino, J.B. Seoane-Sep{\'u}lveda, Bounded and unbounded polynomials and multilinear forms: characterizing continuity, \textit{Linear Algebra Appl.} \textbf{436} (2012), no. 1, 237--242.

\bibitem{GV1999} J.C. Garc\'ia-V\'azquez, R. Villa, Lower bounds for multilinear forms defined on Hilbert spaces, \textit{Mathematika} \textbf{46} (1999) 315--322.

\bibitem{Gard84} L.T. Gardner, An elementary proof of the Russo-Dye theorem, \textit{Proc. Amer. Math. Soc.} \textbf{90} (1984), no. 1, 171.

\bibitem{GelNai43} I. Gelfand, M. Naimark, On the imbedding of normed rings into the ring of operators in Hilbert space, \textit{Mat. Sbornik} \textbf{12} (1943), 197--213.

\bibitem{Gle} A.M. Gleason, A characterization of maximal ideals, \textit{J. Analyse Math.} \textbf{19} (1967), 171--172.

\bibitem{Glob75} J. Globevnik, On complex strict and uniform convexity, \textit{Proc. Amer. Math. Soc.} \textbf{47} (1975), 175--178.

\bibitem{Goodearl} K.R. Goodearl, \textit{Notes on real and complex C$^*$-algebras}, Shiva Mathematics Series, 5. Shiva Publishing Ltd., Nantwich, 1982.

\bibitem{GK} G. Godefroy and N.J. Kalton, Lipschitz-free Banach spaces. Dedicated to Professor Aleksander Pe\l czy\'{n}ski on the occasion of his 70th birthday, \textit{Studia Math.} 159 (2003), no. 1, 121--141.

\bibitem{Haa1990} U. Haagerup, On convex combinations of unitary operators in C$^*$-algebras, \textit{Mappings of operator algebras} (Philadelphia, PA, 1988), 1--13, Progr. Math., 84, Birkh\"{a}user Boston, Boston, MA, 1990.

\bibitem{HaaKadPed2007} U. Haagerup, R.V. Kadison, G.K. Pedersen, Means of unitary operators, revisited, \textit{Math. Scand.} \textbf{100} (2007), no. 2, 193--197.

\bibitem{Hajek2014} P. H\'ajek, M. Johanis, \textit{Smooth Analysis in Banach Spaces}, De Gruyter Ser. Nonlinear Anal. Appl. vol.19, DeGruyter, Berlin, 2014.

\bibitem{HAL} P.R. Halmos, \textit{A Hilbert space problem book}, Second edition. Graduate Texts in Mathematics, 19. Encyclopedia of Mathematics and its Applications, 17. Springer-Verlag, New York-Berlin, 1982.

\bibitem{HOS} H. Hanche-Olsen, E. St{\o}rmer, \textit{Jordan operator algebras}, Pitman, London, 1984.

\bibitem{Harr} L.A. Harris, Bounded symmetric homogeneous domains in infinite dimensional spaces, \textit{Proceedings on Infinite Dimensional Holomorphy} (Internat. Conf., Univ. Kentucky, Lexington, Ky., 1973), pp. 13--40. Lecture Notes in Math., Vol. \textbf{364}, Springer, Berlin, 1974.

\bibitem{Harris1975} L.A. Harris, Bounds on the derivatives of holomorphic functions of vectors, \textit{Analyse fonctionnelle et applications} (Comptes Rendus Colloq. Analyse, Inst. Mat., Univ. Federal Rio de Janeiro, Rio de Janeiro, 1972), pp. 145--163. Actualités Aci. Indust., No. 1367, Hermann, Paris, 1975.

\bibitem{Harris1997} L.A. Harris, A Bernstein-Markov theorem for normed spaces, \textit{J. Math. Anal. Appl.} \textbf{208} (1997), no. 2, 476--486.

\bibitem{Harris2010} L.A. Harris, A proof of Markov's theorem for polynomials on Banach spaces, \textit{J. Math. Anal. Appl.} \textbf{368} (2010), no. 1, 374--381

\bibitem{HaMiOkTak07} O. Hatori, T. Miura, H. Oka, H. Takagi, 2-local isometries and 2-local automorphisms on uniform algebras, \textit{Int. Math. Forum} \textbf{50} (2007), 2491--2502.

\bibitem{HOO} N. D. Hooker, Lomonosov's hyperinvariant subspace theorem for real spaces, \textit{Math. Proc. Cambridge Philos. Soc.} \textbf{89} (1981), no. 1, 129--133.

\bibitem{H1973} O. Hustad, A note on complex spaces. \textit{Isr. J. Math.} \textbf{16} (1973), 117--119.

\bibitem{IliKuzLiPoon} D. Ili\v{s}evic, B. Kuzma, Ch.-K. Li, E. Poon, Complexifications of real Banach spaces and their isometries, \textit{Linear Algebra Appl.} \textbf{589} (2020), 222--241.

\bibitem{Inglestam62} L. Ingelstam, A vertex property for Banach algebras with identity, \textit{Math. Scand.} \textbf{11} (1962), 22--32.

\bibitem{Inglestam64} L. Ingelstam, Real Banach algebras, \textit{Ark. Math.} \textbf{5} (1964), 239--270.

\bibitem{Inglestam66} L. Ingelstam, Symmetry in real Banach algebras, \emph{Math. Scand.} \textbf{18} (1966), 53-68.

\bibitem{Iord03} R. Iord\u{a}nescu, \textit{Jordan structures in geometry and physics. With an appendix on Jordan structures in analysis}, Editura Academiei Rom\^{a}ne, Bucharest, 2003.

\bibitem{IsKaRo95} J.M. Isidro, W. Kaup, A. Rodr{\'\i}guez-Palacios, On real forms of JB$^*$-triples, \textit{Manuscripta Math.} \textbf{86} (1995), 311--335.

\bibitem{IsRo} J.M. Isidro, A. Rodr{\'i}guez-Palacios, On the definition of real W$^*$-algebras, \textit{Proc. Amer. Math. Soc.}, \textbf{124} (1996), 3407--3410.

\bibitem{JimVMorVill2010} A. Jim{\'e}nez-Vargas, A. Morales Campoy, M. Villegas-Vallecillos, Algebraic reflexivity of the isometry group of some spaces of Lipschitz functions, \textit{J. Math. Anal. Appl.} \textbf{366} (2010), no. 1, 195--201.

\bibitem{JimVVill20} A. Jim{\'e}nez-Vargas, M. Villegas-Vallecillos, 2-iso-reflexivity of pointed Lipschitz spaces, \textit{J. Math. Anal. Appl.} \textbf{491} (2020), no. 2, 124359.

\bibitem{Jordan33} P. Jordan, \"{U}ber Verallgemeinerungsm\"{o}glichkeiten des Formalismus der Quantenmechanik, \textit{Nachr. Akad. Wiss. G\"{o}ttingen. Math. Phys. Kl. I}, \textbf{41} (1933), 209--217.

\bibitem{JorvNWig34} P. Jordan, J. von Neumann, E. Wigner, On an algebraic generalization of the quantum mechanical formalism, \textit{Ann. of Math. (2)}, \textbf{35} (1934), no. 1, 29--64.

\bibitem{KadKell2000} V.M. Kadets, A.Yu. Kellerman, On complex strictly convex complexification of Banach spaces, \textit{Mat. Fiz. Anal. Geom.} \textbf{7} (2000), no. 3, 299--307.

\bibitem{Kad51} R.V. Kadison, Isometries of operator algebras, \textit{Ann. Math.} \textbf{54} (1951), 325--338.

\bibitem{Kad52} R.V. Kadison, A generalized Schwarz inequality and algebraic invariants for operators algebras, \textit{Ann. of Math.} \textbf{56} (1952), 494--503.

\bibitem{KadPed85} R.V. Kadison, G.K. Pedersen, Means and convex combinations of unitary operators, \textit{Math. Scand.} \textbf{57} (1985), no. 2, 249--266.

\bibitem{KR1} R.V. Kadison, J.R. Ringrose, \textit{Fundamentals of the theory of operator algebras. {V}ol. {I}, Elementary theory}, vol.~100 of \textit{Pure and Applied Mathematics}, Academic Press, Inc. [Harcourt Brace Jovanovich, Publishers], New York, 1983.

\bibitem{KaZe} J.P. Kahane, W. \.{Z}elazko, A characterization of maximal ideals in commutative Banach algebras, \textit{Studia Math.} \textbf{29} (1968), 339--343.

\bibitem{KalPe2020} O.F.K. Kalenda, A.M. Peralta, Extension of isometries from the unit sphere of a rank-$2$ Cartan factor, \textit{Anal. Math. Phys.} \textbf{11} (2021), no. 1, Paper No. 15, 25 pp.

\bibitem{KAL} N.J. Kalton, An elementary example of a Banach space not isomorphic to its complex conjugate, \textit{Canad. Math. Bull.} \textbf{38} (1995), no. 2, 218--222.

\bibitem{KW} N.J. Kalton, G.V. Wood, Orthonormal systems in Banach spaces and their applications, \textit{Math. Proc. Cambridge Philos. Soc}, \textbf{79} (1976), no. 3, 493--510.

\bibitem{Kap49} I. Kaplansky, Normed algebras, \emph{Duke. Math. J.} \textbf{16} (1949), 399-418.

\bibitem{Kap51} I. Kaplansky, A theorem on rings of operators, \textit{Pacific J. Math.} \textbf{1} (1951), 227--232.

\bibitem{Ka83} W. Kaup, A Riemann mapping theorem for bounded symmentric domains in complex Banach spaces, \textit{Math. Z.} \textbf{183} (1983), 503-529.

\bibitem{Ka84} W. Kaup, Contractive projections on Jordan C$^*$-algebras and generalizations, \textit{Math. Scand.} \textbf{54} (1984), no. 1, 95--100.

\bibitem{Ka97} W. Kaup, On real Cartan factors, \textit{Manuscripta Math.} \textbf{92} (1997), 191--222.

\bibitem{K1952} J.L. Kelley, Banach spaces with the extension property. \textit{Trans. Am. Math. Soc.} \textbf{72} (1952), 323--326.

\bibitem{Kellogg} O.D. Kellogg, On bounded polynomials in several variables, \textit{Math. Z.} \textbf{27} (1928), no. 1, 55--64.

\bibitem{KelleyVaught} J.L. Kelley, R.L., Vaught, The positive cone in Banach algebras, \textit{Trans. Amer. Math. Soc.} \textbf{74} (1953), 44--55.

\bibitem{MOS2} M. Kian, M.S. Moslehian, R. Nakamoto, Asymmetric Choi-Davis inequalities, \textit{Linear Multilinear Algebra}, 2021, doi: 10.1080/03081087.2020.1836115.

\bibitem{K1994} J. L. King, Three problems in search of a measure, \textit{Amer. Math. Monthly} \textbf{101} (1994), no. 7, 609--628

\bibitem{Kirwan} P. Kirwan, Complexification of multilinear mappings and polynomials, \textit{Math. Nachr.} \textbf{231} (2001), 39--68.

\bibitem{Sarant1999} P. Kirwan, Y. Sarantopoulos, A. Tonge, Extremal homogeneous polynomials on real normed spaces, \textit{J. Approx. Theory} \textbf{97} (1999), no. 2, 201--213.

\bibitem{KoszMarMe} P. Koszmider, M. Mart{\'i}n, J. Mer{\'i}, Extremely non-complex $C(K)$ spaces, \textit{J. Math. Anal. Appl.} \textbf{350} (2009), no. 2, 601--615.

\bibitem{KoSlod} S. Kowalski, Z. S{\l}odkowski, A characterization of multiplicative linear functionals in Banach algebras, \textit{Studia Math.} \textbf{67} (1980), 215--223.

\bibitem{KP1999} A. Kro\'o, I. Pritsker, A sharp version of Mahler's inequality for products of polynomials, \textit{Bull. London Math. Soc.} \textbf{31} (1999), no. 3, 269--278.

\bibitem{Kul84} S.H. Kulkarni, Gleason--Kahane--\.{Z}elazko theorem for real Banach algebras, \textit{J. Math. Phys. Sci.} \textbf{18} (1983/84), Special Issue, S19--S28.

\bibitem{LW} E. Lacey, P. Wojtaszczyk, {Banach lattice structures on separable $L\sb{p}$ spaces}, \textit{Proc. Amer. Math. Soc.} \textbf{54} (1976), 83--89.

\bibitem{LarSou} D.R. Larson, A.R. Sourour, Local derivations and local automorphisms of $B(X)$, \textit{Proc. Sympos. Pure Math.} \textbf{51}, Part 2, Providence, Rhode Island 1990, pp. 187--194.

\bibitem{LeNgWon2009} C.-W. Leung, C.-K. Ng, N.-C. Wong, Geometric unitaries in JB-algebras, \textit{J. Math. Anal. Appl.} \textbf{360} (2009), 491--494.

\bibitem{Li75} B.R. Li, Real C$^*$-algebras (in Chinese), \textit{Acta Math. Sinica} \textbf{18} (1975), 216--218.

\bibitem{LI} B.R. Li, \textit{Real operator algebras}, World Scientific Publishing Co., Inc., River Edge, NJ, 2003.

\bibitem{LiPeWangWang} L. Li, A.M. Peralta, L. Wang, Y.-S. Wang, Weak-$2$-local isometries on uniform algebras and Lipschitz algebras, \textit{Publ. Mat.} \textbf{63} (2019), 241--264.

\bibitem{LT} J. Lindenstrauss, L. Tzafriri, \textit{Classical Banach spaces I}, Sequence spaces. Ergebnisse der Mathematik und ihrer Grenzgebiete, Vol. 92. Springer-Verlag, Berlin-New York, 1977.

\bibitem{LMS1998} A.E. Litvak, V.D. Milman, G. Schechctman, Averages of norms and quasi-norms, \textit{Math. Ann.} \textbf{312} (1998) 95--124.

\bibitem{LOM} V.I. Lomonosov, Invariant subspaces of the family of operators that commute with a completely continuous operator (Russian), \textit{Funkcional. Anal. i Prilo\v{z}en.} \textbf{7}
(1973), no. 3, 55--56.

\bibitem{LOM2000} V.I. Lomonosov, \textit{A counterexample to the Bishop-Phelps theorem in complex spaces}, \textit{Israel J. Math.} \textbf{115}, 25--28 (2000).

\bibitem{LoMarPa} G. L{\'o}pez, M. Mart{\'i}n, R. Pay{\'a}, Real Banach spaces with numerical index 1, \textit{Bull. London Math. Soc.} \textbf{31} (1999), no. 2, 207-212.

\bibitem{LUM} G. Lumer, \textit{Complex methods, and the estimation of operator norms and spectra from real numerical ranges}, \textit{J. Funct. Anal.} \textbf{10} (1972), 482--495.

\bibitem{Malicet} D. Malicet, I. Nourdin, G. Peccati, G. Poly, Squared chaotic random variables: new moment inequalities with applications, \textit{J. Funct. Anal.} \textbf{270} (2016), no. 2, 649--670.


\bibitem{Mank71} P. Mankiewicz, A superreflexive Banach space $X$ with $L(X)$ admitting a homomorphism onto the Banach algebra $C(\beta N)$, \textit{Israel J. Math.} \textbf{65} (1989), no. 1, 1--16.

\bibitem{MankTomcz2003} P. Mankiewicz, N. Tomczak-Jaegermann, Quotients of finite-dimensional Banach spaces; random phenomena. Handbook of the geometry of Banach spaces, Vol. 2, 1201--1246, North-Holland, Amsterdam, 2003.


\bibitem{Marcus1997} M. Marcus, A lower bound for the product of linear forms, \textit{Linear Multilinear Algebra} \textbf{43} (1–3) (1997) 115--120.

\bibitem{MAR} V.A. Markov, On polynomials least deviating from zero in a given interval. With a preface by Serge Bernstein. (\"{U}ber Polynome, die in einem gegebenen Intervalle m\"{o}glichst wenig von Null abweichen.) (German) \textit{Math. Ann.} \textbf{77} (1916), 213--258.

\bibitem{MarPe} J. Mart\'{\i}nez, A.M. Peralta, Separate weak$^*$-continuity of the triple product in dual real JB$^*$-triples, \textit{Math. Z.}
\textbf{234} (2000), 635--646.

\bibitem{MashreRansRans2015} J. Mashreghi, T. Ransford, A Gleason-Kahane-\'{Z}elazko theorem for modules and applications to holomorphic function spaces, \textit{Bull. Lond. Math. Soc.} \textbf{47} (2015), 1014--1020.

\bibitem{MashreRansRans2018} J. Mashreghi, J. Ransford, T. Ransford, A Gleason-Kahane-\'{Z}elazko theorem for the Dirichlet space, \textit{J. Funct. Anal.} \textbf{274} (11) (2018), 3254--3263.

\bibitem{Mathieu89} M. Mathieu, Weakly compact homomorphisms from C$^*$-algebras are of finite rank, \textit{Proc. Amer. Math. Soc.} \textbf{107} (1989), no. 3, 761--762.

\bibitem{Matolcsi2005_1} M. Matolcsi, The libear polarization constant of ${\mathbb R}^n$, \textit{Acta Math. Hungar}, \textbf{108} (1–2) (2005), 129--136.


\bibitem{Matolcsi2005_2} M. Matolcsi, A geometric estimate on the norm of product of functionals, \textit{Linear Algebra Appl.}, \textbf{405} (2005), 304--310.

\bibitem{MitiEldel70} B.S. Mityagin, I.C. Edelstein, Homotopy of linear  groups for two classes of Banach spaces  (Russian), \textit{Funkt. Analiz. Priloz.} \textbf{4}, no. 3, (1970), 61--72.

\bibitem{Molnar19} L. Moln{\'a}r, On 2-local $^*$-automorphisms and 2-local isometries of $B(H)$, \textit{J. Math. Anal. Appl.} \textbf{479} (2019), no. 1, 569-580.


\bibitem{monta} A. Montanaro, Some applications of hypercontractive inequalities in quantum information theory, \textit{J. Math. Phys.} \textbf{53} (2012), no. 12, 122206, 15 pp.

\bibitem{MoriOza2020} M. Mori, N. Ozawa, Mankiewicz's theorem and the Mazur--Ulam property for C$^*$-algebras, \textit{Studia Math.} \textbf{250} (2020), no. 3, 265--281.

\bibitem{MOS3} M.S. Moslehian, A. Kusraev, M. Pliev, Matrix KSGNS construction and a Radon-Nikodym type theorem, \textit{Indag. Math.} (N.S.) \textbf{28} (2017), no. 5, 938--952.

\bibitem{MOS4} M.S. Moslehian, Q. Xu, A. Zamani, Seminorm and numerical radius inequalities of operators in semi-Hilbertian spaces, \textit{Linear Algebra Appl.} \textbf{591} (2020), 299--321.

\bibitem{M1998} G.A. Mu\~noz, Complexifications of polynomials and multilinear maps on real Banach spaces. Function spaces (Pozna\.{n}, 1998), 389--406, Lecture Notes in Pure and Appl. Math., 213, Dekker, New York, 2000.


\bibitem{Munoz2002} G.A. Mu\~noz, Y. Sarantopoulos, Bernstein and Markov-type inequalities for polynomials on real Banach spaces,
\textit{Math. Proc. Cambridge Philos. Soc.} \textbf{133} (2002), no. 3, 515--530.


\bibitem{MUS} G.A.. Mu\~{n}oz, Y. Sarantopoulos, J.B. Seoane-Sep\'{u}lveda, {The real plank problem and some applications}, \textit{Proc. Amer. Math. Soc.} \textbf{138} (2010), no. 7, 2521--2535.


\bibitem{MUN} G.A. Mu\~{n}oz, Y. Sarantopoulos, A. Tonge, {Complexifications of real Banach spaces, polynomials and multilinear maps}, \textit{Studia Math.} \textbf{134} (1999), no. 1, 1--33.

\bibitem{N1950} L. Nachbin, A theorem of the Hahn-Banach type for linear transformations. \textit{Trans. Am.Math.
Soc.} \textbf{68} (1950), 28--46.

\bibitem{NavNav2012} J.C. Navarro-Pascual, M.A. Navarro, Unitary operators in real von Neumann algebras, \textit{J. Math. Anal. Appl.} \textbf{386} (2012), no. 2, 933--938.

\bibitem{NEE} J.M.A.M. van Neerven, The norm of a complex Banach lattice, \textit{Positivity} \textbf{1} (1997), no. 4, 381--390.

\bibitem{NT1986} P. Nevai, V. Totik, Weighted polynomial inequalities, \textit{Constr. Approx.} \textbf{2} (1986) 113--127.


\bibitem{Nguyen2009} T. Nguyen, A lower bound on the radius of analyticity of a power series in a real Banach space, \textit{Studia Math.} \textbf{191} (2009), 171--179.

\bibitem{Oi19} S. Oi, A generalization of the Kowalski-S{\l}odkowski theorem and its application to 2-local maps on function spaces, \textit{J. Aust. Math. Soc.} (2020), 1--26.

\bibitem{Palmer70} T.W. Palmer, Real C$^*$-algebras, \emph{Pacific J. Math.} \textbf{35} (1970), 195--204.

\bibitem{Palmer} T.W. Palmer, \textit{Banach algebras and the general theory of $^*$-Algebras. Vol. I.} Encyclopedia Math. Appl. 49 Cambridge University Press, 1994.

\bibitem{PS2016} M. K. Papadiamantis, Y. Sarantopoulos, Polynomial estimates on real and complex $L_p(\mu)$ spaces, \textit{Studia Math.} \textbf{235} (2016), no. 1, 31--45.

\bibitem{Sarantopoulos2016} M. K. Papadiamantis, Y. Sarantopoulos, Radius of analyticity of a power series on real Banach spaces, \textit{J. Math. Anal. Appl.} \textbf{434} (2016), no. 2, 1281--1289.

\bibitem{PR2004} A. Pappas, S, R\'ev\'esz, Linear polarization constants of Hilbert spaces, \textit{J. Math. Anal. Appl.} \textbf{300} (2004) 129--146.

\bibitem{PaSi} A.L.T. Paterson, A.M. Sinclair, Characterisation of isometries between C$^*$-algebras, \textit{J. London Math. Soc.} \textbf{5} (1972), 755--761.

\bibitem{Paulsen} V.I. Paulsen, \textit{Completely bounded maps and operator algebras}, {Cambridge Studies in Advanced Mathematics}, Vol. 78, Cambridge University Press, 2002.

\bibitem{Ped} G.K. Pedersen, \textit{C$^*$-algebras and their automorphism groups}, {London Mathematical Society Monographs}, Vol. 14, Academic Press, London, 1979.

\bibitem{psseo} D. Pellegrino, J.B. Seoane-Sep\'ulveda, New upper bounds for the constants in the Bohnenblust-Hille inequality, \textit{J. Math. Anal. Appl.} \textbf{386} (2012), no. 1, 300--307.

\bibitem{Pe03} A.M. Peralta, On the axiomatic definition of real JB$^*$-triples, \textit{Math. Nachr.} \textbf{256} (2003), 100--106.

\bibitem{Pe2018} A.M. Peralta, A survey on Tingley's problem for operator algebras, \textit{Acta Sci. Math. (Szeged)} \textbf{84} (2018), 81--123.

\bibitem{Phelps65} R.R. Phelps, Extreme points in function algebras, \textit{Duke Math. J.} \textbf{32} (1965), 267--277.

\bibitem{P1940} R.S. Phillips, On linear transformations. \textit{Trans. Am. Math. Soc.} \textbf{48} (1940), 516--541.

\bibitem{P2012} D. Pinasco, Lower bounds for norms of products of polynomials via Bombieri inequality, \textit{Trans. Amer. Math. Soc.} \textbf{364} (2012), no. 8, 3993--4010.

\bibitem{PisBook} G. Pisier, \textit{Espaces de Banach quantiques: une introduction à la théorie des espaces d'opérateurs.} SMF Journ. Annu., 1994, Soc. Math. France, Paris, 1994.

\bibitem{plichkozagorodnyuk1998} A. Plichko, A. Zagorodnyuk, On automatic continuity and three problems of The Scottish book
concerning the boundedness of polynomial functionals, \textit{J. Math. Anal. Appl.} \textbf{220} (1998), no. 2, 477--494.

\bibitem{RACK} H.-J. Rack, A generalization of an inequality of V. Markov to multivariate polynomials, \textit{J. Approx. Theory}, \textbf{35} (1982),
94--97.

\bibitem{RACK1} H.-J. Rack, A genercasealization of an inequality of V. Markov to multivariate polynomials, II, \textit{J. Approx. Theory}, \textbf{40} (1984), 129--133.

\bibitem{RAN} B. Randrianantoanina, Contractive projections and isometries in sequence spaces, \textit{Rocky Mountain J. Math.} \textbf{28} (1998), no. 1, 323--340.

\bibitem{Read1985} C.J. Read, A solution to the invariant subspace problem on the space $l_1$, \textit{Bull. London Math. Soc.} \textbf{17} (1985), no. 4, 305--317.

\bibitem{REIMER} M. Reimer, On multivariate polynomials of least deviation from zero on the unit cube, \textit{J. Approx. Theory}, \textbf{23} (1978), 65--69.


\bibitem{RS2003} S. R\'ev\'esz, Y. Sarantopoulos, On Markov constants of homogeneous polynomials over real normed spaces, \textit{East J. Approx.} \textbf{9} (2003), no. 3, 277--304.

\bibitem{RS2004} S. R\'ev\'esz, Y. Sarantopoulos, Plank problems, polarization and Chebyshev constants, \textit{J. Korean Math. Soc.} \textbf{41} (2004), 157--174.

\bibitem{Rod2010} A. Rodríguez Palacios, Banach space characterizations of unitaries: a survey, \textit{J. Math. Anal. Appl.} \textbf{369} (2010), no. 1, 168--178.

\bibitem{RoiStern81} M. Roitman, Y. Sternfeld, When is a linear functional multiplicative?, \textit{Trans. Amer. Math. Soc.} \textbf{267} (1981), 111--124.

\bibitem{RUA03} Z.-J. Ruan, On real operator spaces, \textit{Acta Math. Sinica (English Series)} \textbf{19} (2003), 485--496.

\bibitem{RUA} Z.-J. Ruan, Complexifications of real operator spaces, \textit{Illinois J. Math.} \textbf{47} (2003), no. 4, 1047--1062.

\bibitem{RUA2} Z.-J. Ruan, \textit{real operator spaces}, International Workshop on Operator Algebra and Operator Theory (Linfen, 2001), \textit{Acta Math. Sin.} (Engl. Ser.) 19 (2003), no. 3, 485--496.

\bibitem{Rudin73} W. Rudin, \textit{Functional analysis}, McGraw-Hill Series in Higher Mathematics. McGraw-Hill Book Co., New York-D\"{u}sseldorf-Johannesburg, 1973.

\bibitem{RuDye} B. Russo, H.A. Dye, A note on unitary operators in C$^*$-algebras, \textit{Duke Math. J.} \textbf{33} (1966), 413--416.

\bibitem{Ryan2002} R. Ryan, \textit{Introduction to tensor products of Banach spaces}, Springer Monographs in Mathematics. Springer-Verlag London, Ltd., London, 2002. xiv+225 pp.

\bibitem{RT} R.A. Ryan, B. Turett, Geometry of spaces of polynomials, \textit{J. Math. Anal. Appl.} \textbf{221} (1998), no. 2, 698--711.

\bibitem{S} S. Sakai, \textit{C$^*$-algebras and $W^*$-algebras}, Springer Verlag. Berlin, 1971.

\bibitem{Sarant1986} Y. Sarantopoulos, \textit{Polynomials and multilinear mappings in Banach spaces.} Ph. D. dissertation, Brunel University, 1986.

\bibitem{Sarant1987} I. Sarantopoulos, Extremal multilinear forms on Banach spaces, \textit{Proc. Amer. Math. Soc.} \textbf{99} (1987), no. 2, 340--346.

\bibitem{Sarant1987BGMS} Y. Sarantopoulos, Polynomials on certain Banach spaces, \textit{Bull. Greek Math. Soc.} \textbf{28} (1987), 89--102.

\bibitem{Sarantopoulos1991} Y. Sarantopoulos, Bounds on the derivatives of polynomials on Banach spaces, \textit{Math. Proc. Cambridge Philos. Soc.} \textbf{110} (1991), no. 2, 307--312.

\bibitem{ScottishBook} The Scottish Book, Mathematics from the Scottish Caf\'{e}, ed. R. D. Mauldin, Birkh\"{a}user, 1981.

\bibitem{SS} K. Seddighi, M. H. Shirdarreh Haghighi, {Sufficient conditions for a linear functional to be multiplicative}, \textit{Proc. Amer. Math. Soc.} \textbf{129} (2001), no. 8, 2385--2393.

\bibitem{Semrl97} P. \v{S}emrl, Local automorphisms and derivations on $B(H)$, \textit{Proc. Amer. Math. Soc.} \textbf{125} (1997), 2677--2680.

\bibitem{Sinc71} A.M. Sinclair, The norm of a hermitian element in a Banach algebra, \textit{Proc. Amer. Math. Soc.} \textbf{28} (1971), 446--450.

\bibitem{SIR} G. Sirotkin, {A version of the Lomonosov invariant subspace theorem for real Banach spaces}, \textit{Indiana Univ. Math. J.} \textbf{54} (2005), no. 1, 257--262.

\bibitem{Skalyga1} V. I. Skalyga, Bounds on the derivatives of polynomials on centrally symmetric convex bodies (Russian), \textit{Izv. Ross. Akad. Nauk Ser. Mat.} \textbf{69} (2005), no. 3, 179--192; translation in \textit{Izv. Math.} \textbf{69} (2005), no. 3, 607--621.

\bibitem{Skalyga2} V.I. Skalyga, Theorems of V. A. Markov in normed spaces (Russian), \textit{Izv. Ross. Akad. Nauk Ser. Mat.} \textbf{72} (2008), no. 2, 193--222; translation in \textit{Izv. Math.} \textbf{72} (2008), no. 2, 383--412.

\bibitem{S1941} A. Sobczyk, Projection of the space $m$ on its subspace $c_0$. \textit{Bull. Am. Math. Soc.} \textbf{47} (1941), 938--947.

\bibitem{Stacho82} L.L. Stach{\'o}, A projection principle concerning biholomorphic automorphisms, \textit{Acta Sci. Math.} \textbf{44} (1982), 99--124.

\bibitem{Stacho02} L.L. Stach{\'o}, A counterexample concerning contractive projections of real JB$^*$-triples, \textit{Publ. Math. Debrecen} \textbf{58} (2001), no. 1-2, 223--230.

\bibitem{Stinespring} W.F. Stinespring, Positive functions on C$^*$-algebras, \textit{Proc. Amer. Math. Soc.} \textbf{6} (1955), 211--216.

\bibitem{Sto} M. Stone, Applications of the Theory of Boolean rings to general topology, \textit{Trans. Amer. Math. Soc.} \textbf{41} (1937), 375--481.

\bibitem{Sza86} S. Szarek, On the existence and uniqueness of complex structure and spaces with ``few'' operators, \textit{Trans. Amer. Math. Soc.} \textbf{293} (1) (1986), 339--353.

\bibitem{Sza86b} S. Szarek, A superreflexive Banach space which does not admit complex structure, \textit{Proc. Amer. Math. Soc.} \textbf{97} (3) (1986), 437--444.

\bibitem{Tak} M. Takesaki,\newblock {\em Theory of operator algebras I}, \newblock Springer, New York, 2003.

\bibitem{Tarski} A. Tarski, Further remarks about the degree of equivalence of polygons (in Polish), \textit{Odbitka A. Parametru}, \textbf{2} (1932), 310--314.


\bibitem{Taylor1938} A.E. Taylor, Additions to the theory of polynomials in normed linear spaces, \textit{T\^ohoku Math. J.}, \textbf{44} (1938), 302--318.

\bibitem{Topp65} D.M. Topping, Jordan algebras of self-adjoint operators, \textit{Mem. Amer. Math. Soc.} \textbf{53} (1965), 48 pp..

\bibitem{TouBrits2020} C. Tour{\'e}, R. Brits, Multiplicative spectral functionals on $C(X)$, \textit{Bull. Aust. Math. Soc.} \textbf{102} (2020), no. 2, 303--307.

\bibitem{TouSchuBrits2017} C. Tour{\'e}, F. Schulz, R. Brits, Multiplicative maps into the spectrum, \textit{Studia Math.} \textbf{239} (2017), 55--66.

\bibitem{TouSchuBrits2018} C. Tour{\'e}, F. Schulz, R. Brits, Some character generating functions on Banach algebras, \textit{J. Math. Anal. Appl.} \textbf{468} (2018), 704--715.

\bibitem{Visser1946} C. Visser, A generalization of Chebyshev's inequality to polynomials in more than one variable, \textit{Indagationes Math.}, \textbf{8} (1946), 310--311.

\bibitem{Wri77} J.D.M. Wright, Jordan C$^*$-algebras, \textit{Michigan Math. J.} \textbf{24} (1977), 291--302.

\bibitem{YangZhao2014} X. Yang, X. Zhao, On the extension problems of isometric and nonexpansive mappings. In: \textit{Mathematics without boundaries}, Edited by Themistocles M. Rassias and Panos M. Pardalos. 725--748, Springer, New York, 2014.

\bibitem{youngson1978vidav} M.A. Youngson, A Vidav theorem for Banach Jordan algebras, \textit{Math. Proc. Cambridge Philos. Soc.} \textbf{84} (1978), 2, 263--272.

\bibitem{Zagorod2001} A. Zagorodnyuk, The zero-subspaces of symmetric polynomials, \textit{Nonlinear Boundary Problems} \textbf{11} (2001), 224--229.

\bibitem{Ze68} W. \.{Z}elazko, A characterization of multiplicative linear functionals in complex Banach algebras, \textit{Studia Math.} \textbf{30} (1968), 83--85.

\bibitem{Ze94} W. \.{Z}elazko, What is known and what is not known about multiplicative linear functionals. Topological vector spaces, algebras and related areas (Hamilton, ON, 1994), 102--115, Pitman Res. Notes Math. Ser., 316, Longman Sci. Tech., Harlow, 1994.


\bibitem{Zettl} H. Zettl, A characterization of ternary rings of operators, \textit{Adv. in Math.} \textbf{48} (1983), no. 2, 117--143.

\bibitem{ZHA} F. Zhang, \textit{Matrix theory: Basic results and techniques}, Second edition. Universitext. Springer, New York, 2011.

\bibitem{Z1977} M. Zippin, The separable extension problem. \textit{Isr. J. Math.} \textbf{26} (1977), 372--387.

\bibitem{Zyl} G. van Zyl, Complexification of the projective and injective tensor products, \textit{Studia Math.} \textbf{189} (2008), no. 2, 105--112.

\end{thebibliography}

\end{document}